\documentclass[preprint,11pt]{elsarticle}

\usepackage{times}
\usepackage{mathptm}
\usepackage{graphicx}
\usepackage{epstopdf}
\usepackage[top=1in, bottom=1in, left=1in, right=1in]{geometry}
\usepackage{fancyhdr}
\usepackage{paralist}
\usepackage{subfigure}
\usepackage{amsmath}
\usepackage{algorithm}
\usepackage{algcompatible}
\usepackage{color}
\usepackage{amssymb}

\usepackage{ulem} % For striking out lines using \sout

\journal{Journal of Computational Physics}

\begin{document}

\begin{frontmatter}


 \title{Estimating the Trace of the Matrix Inverse by Interpolating from the Diagonal of an Approximate Inverse}
 \author[wm]{Lingfei Wu\corref{cor1}}
 \ead{lwu@email.wm.edu} 
 \author[wm]{Jesse Laeuchli} 
 \author[umn]{Vassilis Kalantzis} 
 \author[wm]{Andreas Stathopoulos\corref{cor1}}
 \ead{andreas@cs.wm.edu}
 \author[up]{Efstratios Gallopoulos}  
 \cortext[cor1]{Corresponding author}
 \address[wm]{Department of Computer Science, College of William and Mary, Williamsburg, VA 23187, United States}
 \address[umn]{Department of Computer Science,  University of Minnesota, Minneapolis, MN 55455, United States}
 \address[up]{Department of Computer Engineering and Informatics, University of Patras, Patras, Greece}

%\title{}


%\author{}
%
%\address{}

\begin{abstract}
A number of applications require the computation of the trace of a matrix 
  that is implicitly available through a function.
A common example of a function is the inverse of a large, sparse matrix, 
  which is the focus of this paper.
When the evaluation of the function is expensive, 
  the task is computationally challenging because the standard approach 
  is based on a Monte Carlo method which converges slowly.
We present a different approach that exploits the pattern correlation,
  if present, between the diagonal of the inverse of the matrix and 
  the diagonal of some approximate inverse that can be computed inexpensively. 
We leverage various sampling and fitting techniques to fit the diagonal of 
  the approximation to the diagonal of the inverse. 
Depending on the quality of the approximate inverse, our method may serve 
  as a standalone kernel for providing a fast trace estimate with a 
  small number of samples. 
Furthermore, the method can be used as a variance reduction method 
  for Monte Carlo in some cases.
This is decided dynamically by our algorithm.
An extensive set of experiments with various technique combinations 
  on several matrices from some real applications demonstrate the 
  potential of our method.
\end{abstract}

\begin{keyword}
Matrix trace, Monte Carlo method, variance reduction, preconditioner, fitting, interpolation


\end{keyword}

\end{frontmatter}


%---------------------------------------------------------------
\section{Introduction}
Computing the trace of a matrix $A$ that is given explicitly is a straightforward operation. 
However, for numerous applications we need to compute the trace of a matrix that is given 
  implicitly by its action on a vector $x$, i.e., $Ax$. 
Specifically, many applications are interested in computing the trace of a function of 
  a matrix $F(A)$. 
Examples include estimating parameters in image restoration using the generalized cross-validation approach \cite{hutchinson1990stochastic}, exploring the inverse covariance matrix in uncertainty quantification \cite{bekas2009low,kalantzis2013accelerating}, computing observables in lattice quantum chromodynamics (LQCD) \cite{stathopoulos2013hierarchical}, or counting triangles in large graphs \cite{avron2010counting}. 
The matrix $A$ is large, and often sparse, so its action $F(A)x$ is typically computed 
  through iterative methods.
Because it is challenging to compute $F(A)$ explicitly, the Monte Carlo (MC) approach has become 
the standard method for computing the trace by averaging samples of the bilinear form $x^TF(A)x$ \cite{avron2011randomized,bai1996some}. 
The main purpose of this paper is to develop practical numerical techniques to address the computation of the trace of the inverse of a large, sparse matrix. But our technique can also be adapted to other functions such as the 
trace of the logarithm (yielding the determinant) or the trace of the 
matrix exponential. 

For small size problems, computing $A^{-1}$ through a dense or sparse LDU 
decomposition is the most efficient and accurate approach \cite{duff1986direct, xia2015fast}.
This works well for discretizations of differential operators in low dimensions but becomes intractable in high dimensional discretizations. For larger size problems, domain decomposition and divide and conquer strategies are more tractable but still expensive \cite{tang2011domain}.
In many cases, however, a low accuracy approximation is sufficient. Numerous methods have been presented to address this need for estimating the trace of the inverse of symmetric positive definite matrices through Gaussian bilinear forms \cite{bai1996some, guo2000computing}, modified moments \cite{meurant2009estimates, brezinski2012moments}, 
and Monte Carlo (MC) techniques \cite{bai1996some, guo2000computing, wong2004computing, bekas2007estimator, meurant2009estimates, avron2011randomized, chen2016}. 

MC methods for computing the trace of a matrix are based on the structure of the
Hutchinson method \cite{hutchinson1990stochastic}, which iteratively computes an 
average of matrix bilinear forms with random vectors.
Variants of MC estimators are mainly analyzed and compared based on the variance of one sample \cite{avron2011randomized,roosta2014improved}, which depends on the choice of the selected random vectors.
 %and the structure of $A^{-1}$. One of the ways to accelerate the convergence of the MC procedure is to reduce the variance of the estimator by different choices of random vectors from various distributions \cite{hutchinson1990stochastic, avron2011randomized}, by Hadamard vectors as columns of a Hadamard matrix for MC  \cite{wong2004computing,bekas2007estimator}, by linear combinations of columns of the probing matrix \cite{tang2012probing}, and by combining Hadamard vectors with the probing technique \cite{stathopoulos2013hierarchical}. 
For real matrices, choosing random vectors having each element $\pm 1$ with equal 
  probability is known to minimize variance over all other choices of random vectors \cite{hutchinson1990stochastic,avron2011randomized} and therefore has been widely used in many applications. 
For complex matrices, the same result holds for vectors with $\pm 1, \pm i$ elements.
In \cite{avron2011randomized}, Avron and Toledo analyze the quality of trace estimators through three different metrics such as trace variance, $(\epsilon,\delta)$-approximation of the trace, and the number of random bits for different choices of random vectors. In \cite{roosta2014improved}, Khorasani and Ascher improve the bounds of $(\epsilon,\delta)$-approximation for the Hutchinson, Gaussian and unit vector estimators. 

There has been a number of efforts to combine the Hutchinson method with 
well-designed vectors based on the structure of the matrix 
\cite{wong2004computing,bekas2007estimator,tang2012probing,stathopoulos2013hierarchical}.
In \cite{bekas2007estimator}, the authors use columns of the 
Hadamard matrix, rather than random vectors, to systematically capture 
certain diagonals of the matrix. Then, the MC iteration achieves the 
required accuracy by continuously annihilating more diagonals with more 
Hadamard vectors. However, the location of the nonzeros, or of the large 
elements of $A^{-1}$, often does not coincide with the diagonals annihilated 
by the Hadamard vectors. In \cite{tang2012probing}, graph coloring and 
probing vectors are used to identify and exploit special structures, such 
as bandedness or decaying properties in the elements of $A^{-1}$, to 
annihilate the error contribution from the largest elements. However, if 
the error for the chosen number of colors is large, all work has to be 
discarded and the probing procedure repeated until the accuracy is satisfied.
In \cite{stathopoulos2013hierarchical}, the authors introduce hierarchical 
probing on lattices to avoid the previous problems and achieve the required 
accuracy in an incremental way. For all these approaches, the approximation
error comes from non-zero, off-diagonal elements that have not been 
annihilated yet.
A different MC method samples elements of the main diagonal of the $A^{-1}$ 
  to estimate the trace.
However, its variance depends on the variance of the diagonal which could 
  be much larger than the variance of the Hutchinson method 
  \cite{avron2011randomized}.
This paper focuses on this type of method with the extra assumption that 
  an approximation to the main diagonal of $A^{-1}$ is available.

Our motivation for focusing only on the main diagonal is that
  the trace of $A^{-1}$ is simply a summation of a discrete, 
  1-D signal of either the eigenvalues or the diagonal elements of $A^{-1}$.
Although we cannot compute all the diagonal elements, 
  we may have an approximation to the whole signal from the diagonal of 
  an approximation of $A^{-1}$ (e.g., of a preconditioner).
If the two diagonals have sufficiently correlated patterns, 
  fitting methods can be used to refine the approximation both for 
  the diagonal and the trace.
Therefore, the proposed method may serve as a standalone kernel for 
  providing a good trace estimate with a small number of samples.
But it can also be viewed as a preprocessing method for stochastic 
 variance reduction for MC in cases where the variance reduces sufficiently.
This can be monitored dynamically by our method. 

We present several techniques that improve the robustness of our method 
  and implement dynamic error monitoring capabilities.
Our extensive experiments show that we typically obtain trace estimates with 
much better accuracy than other competing methods, and in some cases the 
  variance is sufficiently reduced to allow for further improvements 
  through an MC method.

%Our goal is to provide an alternative that is competitive with existing approaches when the pattern of the diagonal can be captured well by some approximation.  
%The reminder of this paper is structured as follows. In section 2, we review the Hutchinson method and unit vector estimator, discuss different means of computing approximation, analyze the deflation techniques on reducing stochastic variance as well as comparison of different MC methods. Section 3 presents a number of sampling strategies and illustrates two fitting models to exploit the pattern correlation between the diagonals of the approximation and the diagonals of the matrix inverse. In section 4, we propose a dynamic evaluation framework for estimating the variance of different MC methods and monitor the trace error from fitting stage. In section 5, we give some numerical experiments to demonstrate the effectiveness of the proposed method. Conclusions and further thoughts are gathered in section 8.

%---------------------------------------------------------------
\section{Preliminaries}
We denote by $\|.\|$ the 2-norm of a vector or a matrix, by $N$ the order of $A$, 
by $Z$ an approximation of $A$, 
by $D$ the diagonal elements of $A^{-1}$, 
by $M$ the diagonal elements of $Z^{-1}$, 
by $Tr(F(A))$ the trace of the matrix $F(A)$, 
and by extension, $Tr(D)$ the sum of the elements of the vector $D$,
by $T_{e_i}(F(A))$ the MC trace estimator of $F(A)$ using unit (orthocanonical) vectors, 
by $T_{Z_2}(F(A))$ the MC trace estimator of $F(A)$ using Rademacher vectors, 
by $diag(.)$ the diagonal operator of a matrix, 
and by $Var(.)$ the variance operator of a random variable or a vector.

%---------------------------------------------------------------
\subsection{Hutchinson trace estimator and unit vector estimator}
The standard MC method to estimate the trace of the matrix inverse is due to Hutchinson \cite{hutchinson1990stochastic}.
It estimates the $Tr(A^{-1})$ by averaging $s$ bilinear forms with 
random vectors $z_j \in Z_2^N = \{z(i) = \pm 1 \ \text{with probability} \ 0.5\}$,
\begin{equation} \label{eq: Hutchinson trace of invA}
T_{Z_2}(A^{-1}) = \dfrac{1}{s} \sum_{j=1}^s z_j^{T}A^{-1}z_j.
\end{equation}
The variance of this method is given by 
\begin{equation} \label{eq: Hutchinson variance of invA}
Var(T_{Z_2}(A^{-1})) = \dfrac{2}{s} \|A^{-1}\|_{F}^{2} - \dfrac{2}{s} \sum_{i=1}^{N}\|D_{i}\|^{2} . 
\end{equation}
The variance of this trace estimator is proven to be minimum over  all vectors with real entries \cite{hutchinson1990stochastic}. The confidence interval of a MC method reduces as $O(\sqrt{Var(T_{Z_2}(A^{-1}))})$ for the given matrix. 

The unit vector estimator uniformly samples $s$ vectors from the orthocanonical 
  basis $\{e_1, \ldots, e_N\}$
  \cite{avron2011randomized},
\begin{equation} \label{eq: unit vector trace of invA}
T_{e_i}(A^{-1}) = \dfrac{N}{s} \sum_{j=1}^s e_{i_j}^{T}A^{-1}e_{i_j},
\end{equation}
where $i_j$ are the random indices. The variance of the unit vector estimator is given by
\begin{equation} \label{eq: unit vector variance of invA}
Var(T_{e_i}(A^{-1})) = \dfrac{N^2}{s} Var(D).
\end{equation}
The variance of the Hutchinson method depends on the magnitude of the off-diagonal elements. It converges in one step for diagonal matrices and rapidly if $A^{-1}$ is highly diagonal dominant.
On the other hand, the variance of the unit vector estimator depends only on the variance of the diagonal elements. It converges in one step if the diagonal elements are all the same and rapidly if the diagonal elements are similar. Thus, the method of choice depends on the particular matrix. 

%---------------------------------------------------------------
\subsection{Reducing stochastic variance through matrix approximations}
\label{sec:approximations}
Given an approximation $Z \approx A$, for which $Z^{-1}$ and
  $Tr(Z^{-1})$ are easily computable, we can decompose 
\begin{equation} \label{eq: approximation}
%Tr(A^{-1}) = Tr(Z^{-1}) + Tr(Z^{-1}(ZA^{-1} - I)) = Tr(Z^{-1}) + Tr(E),
Tr(A^{-1}) = Tr(Z^{-1}) + Tr(E),
\end{equation}
where $E = A^{-1} - Z^{-1}$.
We hope that by applying the MC methods on $E$ instead on $A^{-1}$, 
  the variance of the underlying trace estimator,
  in (\ref{eq: Hutchinson variance of invA}) or (\ref{eq: unit vector variance of invA}), can be reduced, thereby accelerating the convergence of MC. 
Among many ways to obtain a $Z$, we focus on the following two.

The first approach is when $Z^{-1} = (LU)^{-1}$, where the $L,\ U$ matrices stem from an incomplete LU (ILU) factorization of $A$, one of the most commonly used preconditioners. 
If the ILU is sufficiently accurate, then $M = diag(Z^{-1})$ may be a good approximation to $D$.
To obtain the vector $M$ without computing the entire $Z^{-1}$, we can use an 
  algorithm described in \cite{erisman1975computing}. 
This algorithm requires the computation of only those entries $Z_{ij}^{-1}$ for which 
$L_{ij}$ or $U_{ij} \neq 0$. If the $L, U$ factors are sufficiently sparse or structured,
this computation can be performed efficiently 
(see \cite{journals/toms/LinYMLYE11} for an example in the symmetric case).

%Another approach is to update the ILU decomposition by a low-rank approximation of the error matrix $\bar{E} = A-LU$. The low-rank approximation can be derived by a Lanczos (Arnoldi for nonsymmetric) procedure. If $\bar{E}$ is nearly a low-rank matrix, then such an update can be used as a preconditioner to accelerate  convergence when solving linear system and also update the approximation $Z$. However, the diagonal of $Z^{-1}$ may not get close to that of $A^{-1}$ even many Lanczos steps since the error, which depends on the smallest magnitude eigenpairs of $\bar{E}$ (or the smallest singular triplets of $\bar{E}$), will not be distributed only along the diagonal. The updated preconditioner and the corresponding updated approximation of the diagonal of $A^{-1}$ can be applied through the Woodburry formula \cite{higham2002accuracy}, 
%\begin{equation*}(\hat{A} + VTV^{\top})^{-1} = \hat{A}^{-1} - \hat{A}^{-1}V(T^{-1} + V^{\top}\hat{A}V)^{-1}V^{\top}\hat{A}
% \end{equation*}
%where $\hat{A} = LU$ and $VTV^{\top}$ is a low-rank approximation of $\bar{E}$. We have to mention that we did not explore this technique in our current research. 

The second approach is a low rank approximation $Z^{-1} = V \Sigma^{-1} U^T$, where $\Sigma$ is a diagonal matrix with the $k \ll N$ smallest singular values of $A$,
and $U$ and $V$ are the corresponding left and right singular vectors.
Eigenvalues and eigenvectors can be used instead but we do not
  consider it in this paper.
This subspace can be obtained directly by an iterative eigensolver \cite{Wu2015PHSVD, wu2016svds} 
  as a preprocessing step. 
The cost of this procedure is relatively small since the singular space 
  is not needed in high accuracy.
Alternatively, the space can be approximated by methods such as eigCG \cite{stathopoulos2010computing} or eigBiCG \cite{abdel2014extending} while solving linear systems of equations during the MC method. 
This incremental approach adds only minimal overhead to the iterative linear 
  solver but it cannot compute as many and as good quality singular vectors 
  as the first approach.
The quality of the approximation of $A^{-1}$ by $Z^{-1}$ depends on the 
  separation of the computed singular space. Therefore, depending on the matrix, 
  more singular triplets may be needed for a good low rank approximation.
On the other hand, this space can also be used to deflate 
  and thus accelerate subsequent linear systems.
Once the singular space is computed, each diagonal element can 
  be obtained with a single inner product of short vectors, 
  $M_i = V(i,:)^TU(i,:)/\sigma_i$, and thus it is computationally 
  inexpensive.
Finally, the incremental SVD approach requires some special algorithmic 
  attention during our algorithm, which will be pointed out later.

A computationally inexpensive, albeit less accurate approach for computing an approximation $M$ is based on variational bounds on the entries of $A^{-1}$
\cite{bai1996some,robinson1992variational}. 
Upper and lower bounds on the $i$-th diagonal entry $ A^{-1}_{ii}$ are derived inexpensively since they only depend on estimates of the smallest and largest algebraic eigenvalues, $\lambda_1,\lambda_N$, and the entries of $A$. 
The bounds apply to both symmetric and unsymmetric matrices. For the case of a real symmetric $A$, we have \cite{robinson1992variational},
% \begin{equation} \label{eq: lower bound}
% \dfrac{1}{\lambda_N} \leq (A^{-1})_{ii} \leq \dfrac{1}{\lambda_1} 
% \end{equation}
\begin{equation} \label{eq:RobinWathen1}
 \dfrac{1}{\lambda_N} + \dfrac{(\lambda_N-A_{ii})^2}{\lambda_N(\lambda_N A_{ii}-s_{ii})} \leq (A^{-1})_{ii} \leq \dfrac{1}{\lambda_1} - 
 \dfrac{(A_{ii}-\lambda_1)^2}{\lambda_1(s_{ii} - \lambda_1 A_{ii})},
 \end{equation}
where $s_{ij} = \sum_{k=1}^N A_{ik} A_{kj}$. However the bounds in (\ref{eq:RobinWathen1}) 
will not be sharp especially the upper bound \cite{meurant2009estimates} and the error in the approximation can be large. 
 
\begin{figure}[h!]
  \centering
  \subfigure[ $M$ does not captures pattern of $D$]
  {\includegraphics[width=0.35\textwidth]{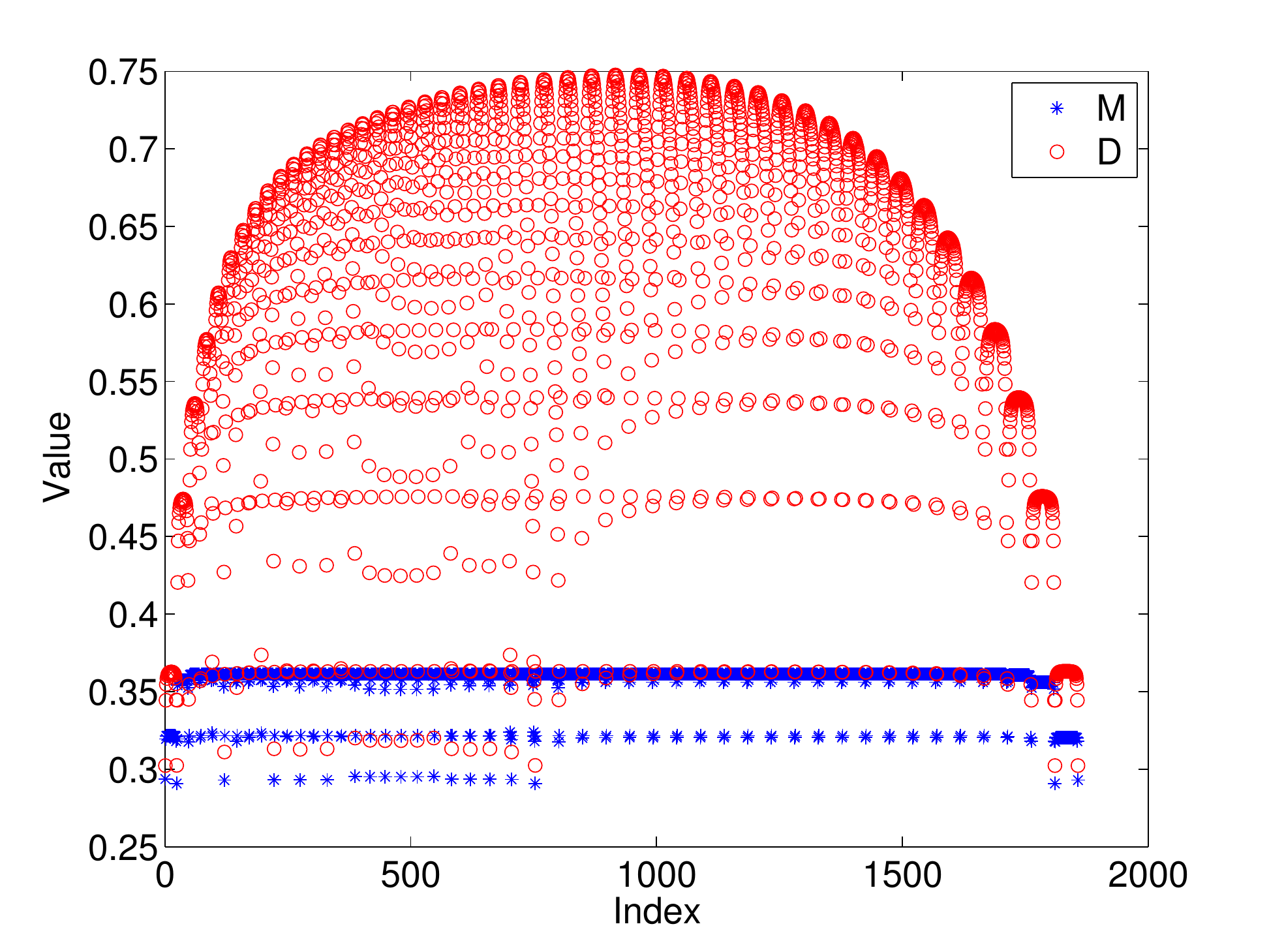}\label{fig:motivation_ILU_bad}}
  \subfigure[$M$ quite close to $D$ and captures pattern]{\includegraphics[width=0.35\textwidth]{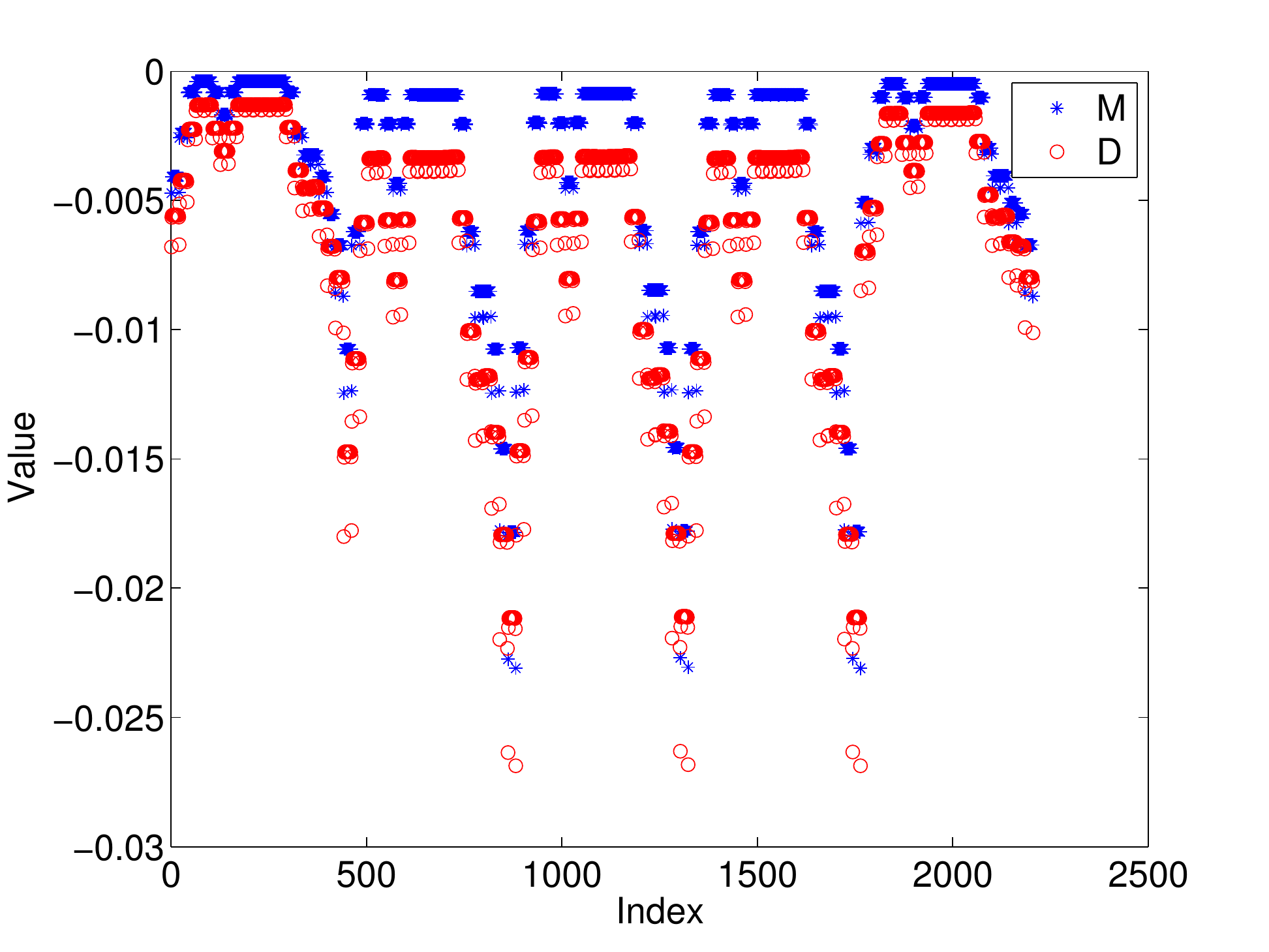}\label{fig:motivation_ILU_good}}   
  \caption{
  The pattern correlation between the diagonals of $A^{-1}$ and its approximation $Z^{-1}$ computed by ILU(0) on matrices (a) delsq50 and (b) orsreg2205. delsq50 is created in MATLAB by {\tt delsq(numgrid('S',50)).}}
  \label{fig:motivation_ILU}
\end{figure}

\begin{figure}
  \centering
  \subfigure[$M$ not close to $D$ but captures pattern]{\includegraphics[width=0.35\textwidth]{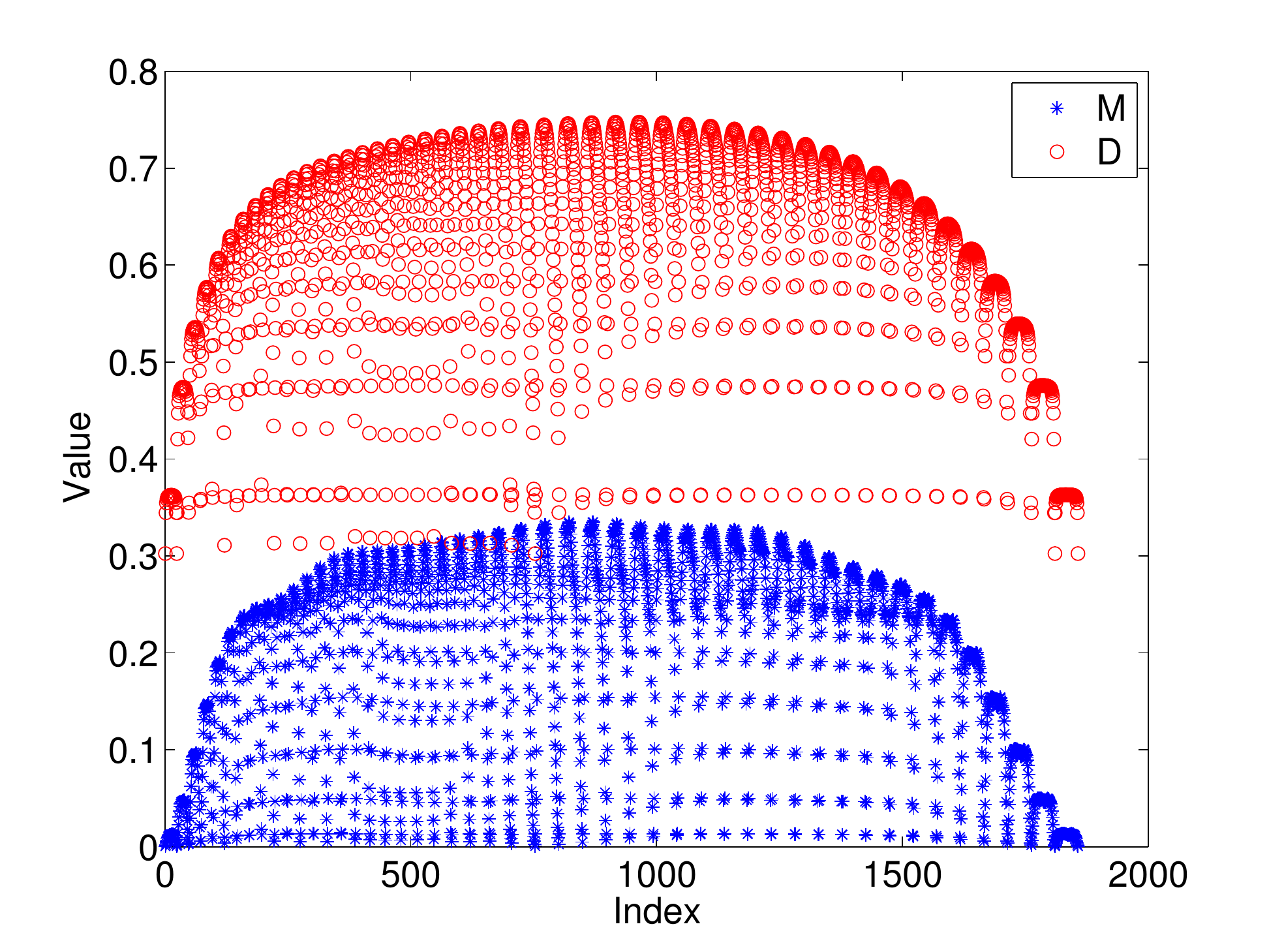}\label{fig:motivation_SVD_bad}}                
  \subfigure[$M$ not close to $D$ but captures pattern]{\includegraphics[width=0.35\textwidth]{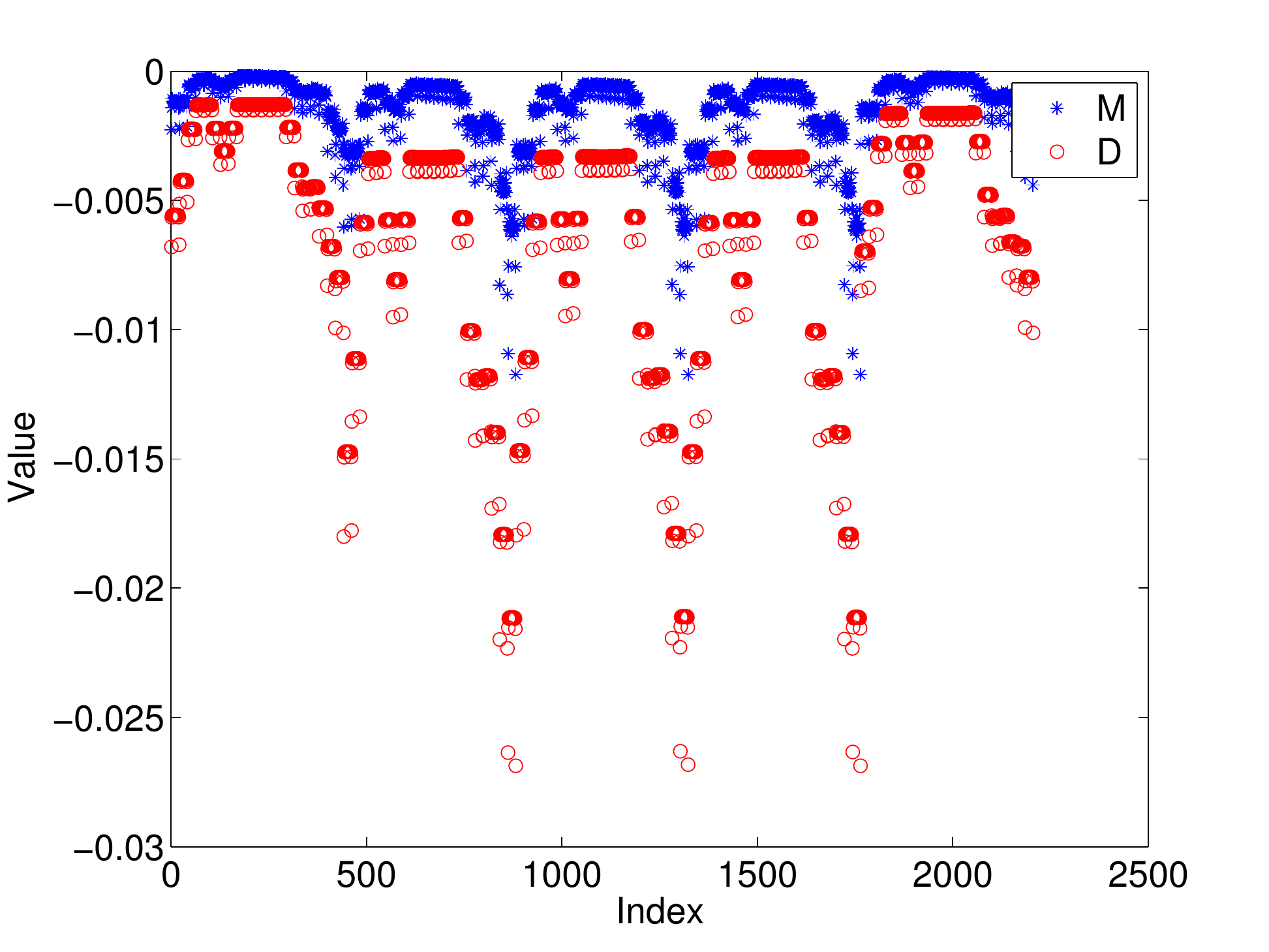}\label{fig:motivation_SVD_good}}   
  \caption{
  The pattern correlation between the diagonals of $A^{-1}$ and its SVD approximation $Z^{-1}$ computed from the 20 smallest singular triplets of $A$ on matrices (a) delsq50 and (b) orsreg2205. delsq50 is created in MATLAB by {\tt delsq(numgrid('S',50)).}}
  \label{fig:motivation_SVD}
\end{figure}

In general, unless $Z^{-1}$ is highly accurate, we do not expect $Tr(M)$ 
  to be close to $Tr(A^{-1})$.
However, the patterns of $M$ and $D$ often show some correlation.
We demonstrate this for two example matrices, delsq50 and orsreg2205, in Figures 
\ref{fig:motivation_ILU} and \ref{fig:motivation_SVD} using ILU and SVD respectively.
For matrix orsreg2205, both the ILU and SVD approaches return an approximate diagonal $M$ which captures the pattern of $D$ very well,
with the $M$ returned by ILU being slightly better than the one from SVD. 
For matrix delsq50, SVD clearly captures the pattern of $D$ better than what ILU does. 
As in preconditioning for linear systems of equations, 
  the appropriate approximation technique depends on the given matrix.
%We have found that for more diagonally dominant matrices, 
%  inexpensive approaches that approximate $M$ based on variational bounds 
%  on the entries of $A^{-1}$ \cite{bai1996some,robinson1992variational} 
%  could be equally useful.
%Here, we focus only on ILU and SVD based approximations.

%---------------------------------------------------------------
\subsection{Comparison of different MC methods and discussion on importance sampling}
Based on 
(\ref{eq: Hutchinson variance of invA}--\ref{eq: unit vector variance of invA}),
we express the variance of the trace estimators $T_{Z_2}(E)$ and $T_{e_i}(E)$ as follows:
%by replacing the matrix $A$ with $E$.
\begin{equation} \label{eq: Hutchinson variance of E}
Var(T_{Z_2}(E)) = \dfrac{2}{s} \|E\|_{F}^{2} - \dfrac{2}{s} \sum_{i=1}^{N}\|diag(E)\|^{2} , 
\end{equation}
\begin{equation} 
	\label{eq: unit vector variance of E}
Var(T_{e_i}(E)) = \dfrac{N^2}{s} Var(diag(E)).
\end{equation}
Figures \ref{fig:motivation_ILU} and \ref{fig:motivation_SVD} show there is potential 
  for the variances of $T_{Z_2}(E)$ or $T_{e_i}(E)$ to be smaller than those
  of $T_{Z_2}(A^{-1})$ and $T_{e_i}(A^{-1})$.
However, for a given matrix, we must gauge which MC method would be better, and 
  whether the variances need further improvement.

The estimator $T_{e_i}(E)$ has the interesting property that if $M = D + c$, where $c$
  is a constant, then its variance in (\ref{eq: unit vector variance of E}) is zero 
  and we obtain the correct trace in one step.
Although we cannot expect this in practice, it means that the shift observed 
  between $M$ and $D$ in Figure \ref{fig:motivation_SVD_bad} should not affect 
  the effectiveness of $T_{e_i}(E)$. 

On the other hand, $T_{e_i}(E)$ fails to identify correlations of the form $M=cD$.
For such cases, importance sampling is preferred, where $M$ plays the role of a new 
  distribution simulating the distribution of $D$.
Assume that both $D$ and $M$ have been shifted by the same shift so that $M_i > 0$ if $D_i > 0,\ i=1,\ldots,N$.
To transform $M$ into a probability mass function, let 
	$G =  \tfrac{1}{Tr(M)} M$.
To obtain an estimator of the trace of $D$ with importance sampling, we replace
  the uniform sampling of $D_i$ values with sampling with probability $G_i$
  \cite{mackay1998introduction}. 
Then, instead of (\ref{eq: unit vector trace of invA}),
  the importance sampling estimator is:
\begin{equation}
T_{IS}(D) = \dfrac{N}{s} \sum_{j=1}^s D_{i_j} \dfrac{\tfrac{1}{N}}{G_{i_j}} = 
	\dfrac{Tr(M)}{s} \sum_{j=1}^s \dfrac{D_{i_j}}{M_{i_j}}.
\end{equation}
When $M=cD$, the variance of $T_{IS}(D)$ is zero and it finds the trace 
  in one step.
However, it completely fails to identify shift correlations.
In general $D$ and $M$ may have a more complex relationship that 
neither $T_{IS}(D)$ or $T_{e_i}(E)$ can capture. 
This motivates our idea to explore general fitting models to 
approximate $D$. 

%---------------------------------------------------------------
\section{Approximating the trace of a matrix inverse} 
We seek to construct a function $f$, such that $D \approx f(M)$. 
Then we can decompose 
\begin{equation} \label{eq:refine approximation}
Tr(A^{-1}) =  Tr(D - f(M)) + Tr(f(M)).
\end{equation}
$Tr(f(M))$ is trivially computed for a given $f$.
A key difference between the approaches in (\ref{eq:refine approximation}) and (\ref{eq: approximation}) is that 
it is easier to find a fitting of the two vectors $M$ and $D$ 
  if a strong pattern correlation exists between them than to fit all corresponding elements in 
 matrices $A^{-1}$ and $Z^{-1}$.
If $Tr(f(M))$ is a good approximation to $Tr(A^{-1})$ and its accuracy can be evaluated easily, we can directly use this quantity; Otherwise, we can apply the unit vector MC estimator to compute $Tr(E_{fit}) = Tr(D - f(M))$, provided that its variance
\begin{equation} \label{eq:unit vector variance of E_fit}
Var(T_{e_i}(E_{fit})) = \dfrac{N^2}{s} Var(E_{fit})
\end{equation}
is smaller than the variances in (\ref{eq: Hutchinson variance of invA}), (\ref{eq: unit vector trace of invA}), (\ref{eq: Hutchinson variance of E}) and (\ref{eq: unit vector variance of E}).

\begin{algorithm}[h!]
\caption{Basic algorithm for approximating $Tr(A^{-1})$ \label{alg:main idea}}
\begin{algorithmic}[1]
\STATEx ${\bf Input:}\ A\in \mathbb{R}^{N\times N}$
\STATEx ${\bf Output:}\ Tr(A^{-1})$ estimation and $\ Z\in \mathbb{R}^{N \times N}$
%\STATEx \% $Z$ and $Z^{-1}$ are given in their implicit forms
%\STATEx \% $Z^{-1}$ also acts as a preconditioner when solving linear systems
\STATE Compute $M=diag(Z^{-1})$, where $Z^{-1}$ is an approximation to $A^{-1}$
	 (Section \ref{sec:approximations})
\STATE Compute fitting sample $S_{fit}$, a set of $k$ indices
	 (Section \ref{sec:pointIdentification})
\STATE Solve linear systems $D_i = e_{i}^{T} A^{-1} e_{i},\ \forall i \in S_{fit}$ (Section \ref{sec:approximations})
\STATE Obtain a fitting model $f(M) \approx D$ by fitting $f(M(S_{fit}))$ to $D(S_{fit})$
	 (Section \ref{section:two fitting models})
\STATE Compute refined trace approximation $T_{e_i}(E_{fit})$ 
	 using (\ref{eq: unit vector trace of invA})
\STATE Estimate the relative trace error and, if needed, the variances for different MC methods
	 (Section \ref{sec:DynamicEvaluation})
\end{algorithmic}
\end{algorithm} 

The basic description of the proposed estimator is outlined in Algorithm \ref{alg:main idea}.
First, our method computes an approximation $M$ of $D$ using one of the methods discussed in the previous section. 
If that method is based on a preconditioner or a low rank approximation, that 
	preconditioner $Z^{-1}$ could also be used to speed up the solution of linear systems
        in step 3.
Second, it finds a fitting sample $S_{fit}$, a set of indices that should capture the important distribution characteristics of $D$. 
Since we have no information about $D$, Section \ref{sec:pointIdentification} discusses
  how to tackle this task by considering the distribution of $M$.
Third, it computes the values of $D(S_{fit})$ by solving the corresponding 
  linear systems using a preconditioned iterative solver.
Since this is the computational bottleneck, the goal is to obtain good accuracy with far 
  fewer fitting points than the number of vectors needed in MC.
Fourth, it computes a fitting model that has sufficient predictive power to improve the diagonal approximation. This critical task is discussed in Section \ref{section:two fitting models}.
Finally, since there are no a-posteriori bounds on the relative error of the trace, 
  we use a combination of statistical approaches and heuristics to estimate it incrementally
  at every step.
If the estimated error is not sufficiently small, our method can be followed by an MC method.
For this reason, the algorithm also estimates dynamically the variances of the two 
  different MC estimators (\ref{eq: Hutchinson variance of E}) and 
  (\ref{eq: unit vector variance of E}) so that the one with the smallest variance is picked.
The dynamic evaluation of error and variances is discussed in 
  Section \ref{sec:DynamicEvaluation}.

%---------------------------------------------------------------
\subsection{Point Identification Algorithm}
\label{sec:pointIdentification}
We need to identify a set of indices $S_{fit}$ based on which we can compute a
  function $f$ that fits $f(M(S_{fit}))$ to $D(S_{fit})$ so that 
  $|Tr(f(M)) - Tr(D)|$ is minimized.
It is helpful to view $Tr(D)$ 
as an integral of some hypothetical one dimensional function 
  which takes the values $D(i) = D_i$ on the discrete points $i=1,\ldots,n$.
Our goal is to approximate $\int_0^N D(x)dx$ with far fewer points than $n$. 
Because it is one dimensional, 
%it may be surprising that Monte Carlo is the standard approach and not numerical integration. 
it may be surprising that the Monte Carlo approach can be
    advantageous than numerical integration.
However, the effectiveness of any numerical quadrature rule relies on the 
  smoothness of the function to be integrated, specifically on the magnitude
  of its higher order derivatives.
In our case, our hypothetical function may have no smoothness or bounded 
  derivatives.
More practically, we deal with a set of discrete data points in $D$ for which 
  smoothness must be defined carefully. 
A typical definition of smoothness for discrete data is based on the Lipschitz
  continuity, i.e., $|D_i -D_j| < c|i-j|$ \cite{discreteSmoothness}.
The lower the constant $c$, the smoother the set of data.
Hence, for a continuous function, the larger the magnitude of its first 
  derivative, the less smooth its discretized points are.
In the context of integrating an arbitrary $D$, the first order divided 
  difference (i.e., the discretized first derivative) $D_i-D_{i+1}$ may be 
  arbitrarily large in magnitude and therefore numerical integration may not 
  work any better than random averaging.

Even for matrices that model physically smooth phenomena (e.g., in PDEs), 
  the matrix may be given in an ordering that does not preserve physical locality.

Consider a sorted permutation of the diagonal, $\hat{D} = sort(D)$.
Obviously $Tr(D) = Tr(\hat D)$, but $\hat D$ is monotonic and, 
   based on the definition of Lipschitz discrete smoothness,
   it is maximally smooth, i.e., 
\begin{equation}
   \label{eq:sortedSmoothness}
	|\hat D(i) - \hat D(j)| \leq \Delta |i-j|,
\end{equation}
  for the smallest possible $\Delta \in \mathbb{R}_+$ 
  among all permutations of $D$\footnote{
   It is easy to see that $\Delta = \max_{i=2,N} |\hat D(i)-\hat D(i-1)|
	= \hat D(i_0)-\hat D(i_0-1)$ is the smallest value that satisfies 
	\eqref{eq:sortedSmoothness} for the sorted array. 
	Assume there is a different permutation $D'$ that satisfies 
	\eqref{eq:sortedSmoothness} with $\delta < \Delta$.
	If $D'(k) = \hat D(i_0)$, then 
	$\delta \geq |D'(k)-D'(k-1)|\geq \hat D(i_0)-\hat D(i_0-1) = \Delta$,
	which contradicts the assumption.
}.
Monotonicity implies that, in the absence of any additional information 
  about the data, a simple trapezoidal rule minimizes the worst case 
  integration error \cite{Kiefer_optimality_1957}.
In addition, if we are allowed to choose the integration points sequentially, based on the points computed so far, 
then a much better average case error can be 
  obtained \cite{Novak_IntegrMonotoneF,Sukharev_seqOpt_Integr}. 
On the other hand, if bounds are known on the discrete smoothness of $\hat D$, 
  better worst case error bounds can be established.
Since $D$, however, is not available, we turn to its approximation $M$.

A close pattern correlation between $M$ and $D$ means that the elements 
  of $M$ should have a similar distribution as those of $D$, or that 
  $\hat M = sort(M)$ should be similar to $\hat D$.
Let us assume for the moment, that the index that sorts $D$ to $\hat D$
  sorts also $M$ to $\hat M$. 
In other words, we assume a complete correlation between the monotonic 
  orderings of $M$ and $D$ even though their values may differ.
Then, we can work on the surrogate model $\hat M$ for which we
  can afford to identify the best quadrature points that yield the
  smallest error in $Tr(\hat M)$.
These will be the ideal points over all permutations of $M$ for computing 
  the fitting function $f$.

Specifically, we need to select indices that capture the important 
  distribution changes in $\hat M$. 
For example, identifying minimum and maximum elements of $M$ sets the 
  range of approximation for $f$ and avoids extrapolation.
We also look for entries in $\hat M$ that deviate highly from their neighbors 
  (where the first order divided difference of $\hat M$ 
  and hopefully of $\hat D$ has a large value).
This strategy has three advantages.
First, between such indices the data is smooth so the integral $Tr(\hat M)$ 
  should be captured well by the trapezoidal rule.
Second, and more important, we can obtain a more accurate fitting 
  function $f$ in a piecewise manner in intervals where $\hat D$ has
  similar behavior.
Third, in this way we avoid picking points with the same value 
  $\hat M_i = \hat M_j$ which can create problems in the fitting function.

The proposed index selection method is shown in 
  Algorithm \ref{alg:point identification algorithm based on the trapezoidal rule}. 
Initially the set of sampled indices, $\hat S_{fit}$, 
  includes the indices of the extrema of $\hat M$, 1 and $N$.
Then, for every interval $(i,j)$, with $i,j$ consecutive indices in $\hat S_{fit}$, 
  we find the index $t$ from $i+1$ to $j-1$ that minimizes the 
  trapezoidal rule error for computing $Tr(\hat M(i:j))$:
\begin{equation} \label{eq:trapezoidal rule}
\underset{{i < t < j,\ t \in \mathbb{Z}}} {\mathrm{argmin}} (|(\hat M(i) - \hat M(j))*(i-j) - (\hat M(i) - \hat M(t))*(i-t) - (\hat M(t) - \hat M(j))*(t-j)|).
\end{equation}
The process continues until it reaches a maximum number of sampling points 
  or until the maximum error over all intervals decreases by a fixed value, 
  say $0.001$. 
The value of this threshold depends on how close the patterns of $M$ and $D$
  are correlated, a question we had postponed and we discuss next.

The points selected by the previous algorithm minimize the error 
  for computing $Tr(M)$ with only $|\hat S_{fit}|$ points. 
How good are these points for fitting $\hat M$ to $\hat D$ and thus 
  obtaining a small error in $Tr(D)$?
Consider the $M$ and $D$ as discrete functions from 
  $\left\{1,\ldots, n\right\} \rightarrow \mathbb{R}$ and assume they are 
  bijective in their ranges.
Let $J, G: \left\{1,\ldots, n\right\} \rightarrow \left\{1,\ldots, n\right\}$
  be the two permutation functions that sort 
  $M$ to $\hat M$ and $D$ to $\hat D$, i.e., 
  $\hat M(i) = M (J(i))$ and
  $\hat D(i) = D (G(i))$.
%$$ \hat M = M \circ J\ \mbox{ and }\ \hat D = D \circ G.$$
Our method relies on the assumption that the distribution of $\hat M$ 
  captures the distribution of $\hat D$, or that there exists
  a smooth function $f$ that can be approximated well with a low 
  degree polynomial that fits $\hat D(i) = f(\hat M(i))$.
For any $i$, there are two indices $m=G(i), k=J(i)$ with $m\neq k$ in general
  that satisfy, 
  $$D(m) = D(G(i)) = \hat D(i) = f(\hat M(i))
  = f( M(J(i)) ) = f(M(k)).$$
Our algorithm first picks an index $i$ for $\hat M$ and 
  derives the original index $k=J(i)$ in $M$.
Since the permutation $G$ is unknown, we do not know the $D(m)$ that 
  should be matched with $M(k)$.
Thus, the algorithm makes the simplifying assumption that $G=J$ and computes
  $D(k)$ instead.
However, this corresponds to a different $\hat D(j)$, where $G(j) = k$.
Therefore, using (\ref{eq:sortedSmoothness}), we can bound the error 
  for this mismatch by 
\begin{equation}
	|\hat D(i) - \hat D(j)| \leq \Delta |i-j| = 
	\Delta | J^{-1}(k) - G^{-1}(k) |. 
	\label{eq:bound}
\end{equation}
This implies that as long as the permutations $J$ and $G$ are 
  locally similar, or in other words they do not shuffle the same index 
  of $M$ and $D$ too far from each other, the error is small 
  and therefore the fitting should work well.
On the other hand, if $M$ is a random permutation of $D$, there is no good fitting,
  even though $\hat M = \hat D$.
This implies that the traces of $f(M)$ and $D$ may be very close 
but not their individual elements. It also means that
 its variance of $E_{fit}$ may not always be small (see the examples in the 
  following section as well as in 
  Section \ref{sec:Monitoring Relative Trace Error}).

In practice, since we do not know the permutation $G$, it is possible
 that some flat area of $\hat M$ is associated with important changes in
 $\hat D$.
To alleviate the effect of a possible local pattern mismatch, 
  we instead empirically insert the midpoint of the current largest 
  interval every $k$ samples (we choose $k=5$ in lines 11--14).
Returning to the choice of threshold in the algorithm, we see
  that going below a threshold helps the accuracy of $Tr(M)$ but 
  not necessarily of $Tr(D)$. 
The reason is that the greedy point selection strategy 
   becomes less effective in determining the local mismatching patterns.
Therefore, we terminate searching for a new index based on 
  (\ref{eq:trapezoidal rule}) if the maximum error is less than 0.001.
We found this threshold to be sufficient in our experiments.
If the maximum required samples
    have not been generated, we continue by simply bisecting the
     largest intervals until $maxPts$ is reached (in lines 16--18).
On exit, we compute $S_{fit}$ which maps $\hat S_{fit}$ to the original unsorted ordering. 

%Finally, we have to mention that if the indices of the elements of $\hat{M}$ and $D$ are not matched, it will not help reduce $Var(T_{e_i}(E_{fit}))$ but have no effect on estimating the trace. We will discuss this issue in details in the following Section \ref{sec:Monitoring Relative Trace Error}.

\begin{algorithm}[h!]
\caption{Point identification algorithm based on the trapezoidal rule}
\begin{algorithmic}[1]
\STATEx ${\bf Input:}\ M \in \mathbb{R}^N$ and $maxPts$
\STATEx ${\bf Output:}\ \hat S_{fit}$: desired sampling index set
	\STATE [$\hat{M}$, $J$] = sort($M$)
	\STATE $numSamples=1$, 
		$initErr = tempErr = |Tr(\hat M)-(\hat{M}(1)+\hat{M}(N))\tfrac{N-1}{2}|$
	\STATE Add $1, N$ into $\hat S_{fit}$ and push interval $(1,N)$ and its $tempErr$ in the queue $Q$ including interval and error
	\WHILE{ $numSamples < maxPts$ and $tempErr > 0.001 * initErr$}
		\STATE Pop interval $(L,R)$ with largest error from $Q$
		\FOR{$t = L+1 : R-1$}
			\STATE Use (\ref{eq:trapezoidal rule}) to find a bisecting index $t$ 
		\ENDFOR
		\STATE Add $t$ into $\hat S_{fit}$, $numSamples = numSamples + 1$
			\STATE Push intervals $(L,t)$ and $(t,R)$ each with their corresponding $tempErr$ in $Q$
 		\IF{$numSamples$ is a multiple of 5}
			\STATE Insert one midpoint index into largest interval in $\hat S_{fit}$, $numSamples = numSamples + 1$
			\STATE Insert its corresponding left and right intervals in $Q$
		\ENDIF
	\ENDWHILE
	\WHILE {$numSamples < maxPts$}
	\STATE Insert middle index of the largest interval into $S$,  $numSamples = numSamples + 1$
	\ENDWHILE
	\STATE Return $S_{fit}$ indices in original ordering,
		such that $\hat S_{fit} = J(S_{fit})$
\end{algorithmic}
\label{alg:point identification algorithm based on the trapezoidal rule}
\end{algorithm}
 
\begin{figure}[h!]
  \centering
  \subfigure[ Distribution of elements of $M$ and $D$]{\includegraphics[width=0.45\textwidth]{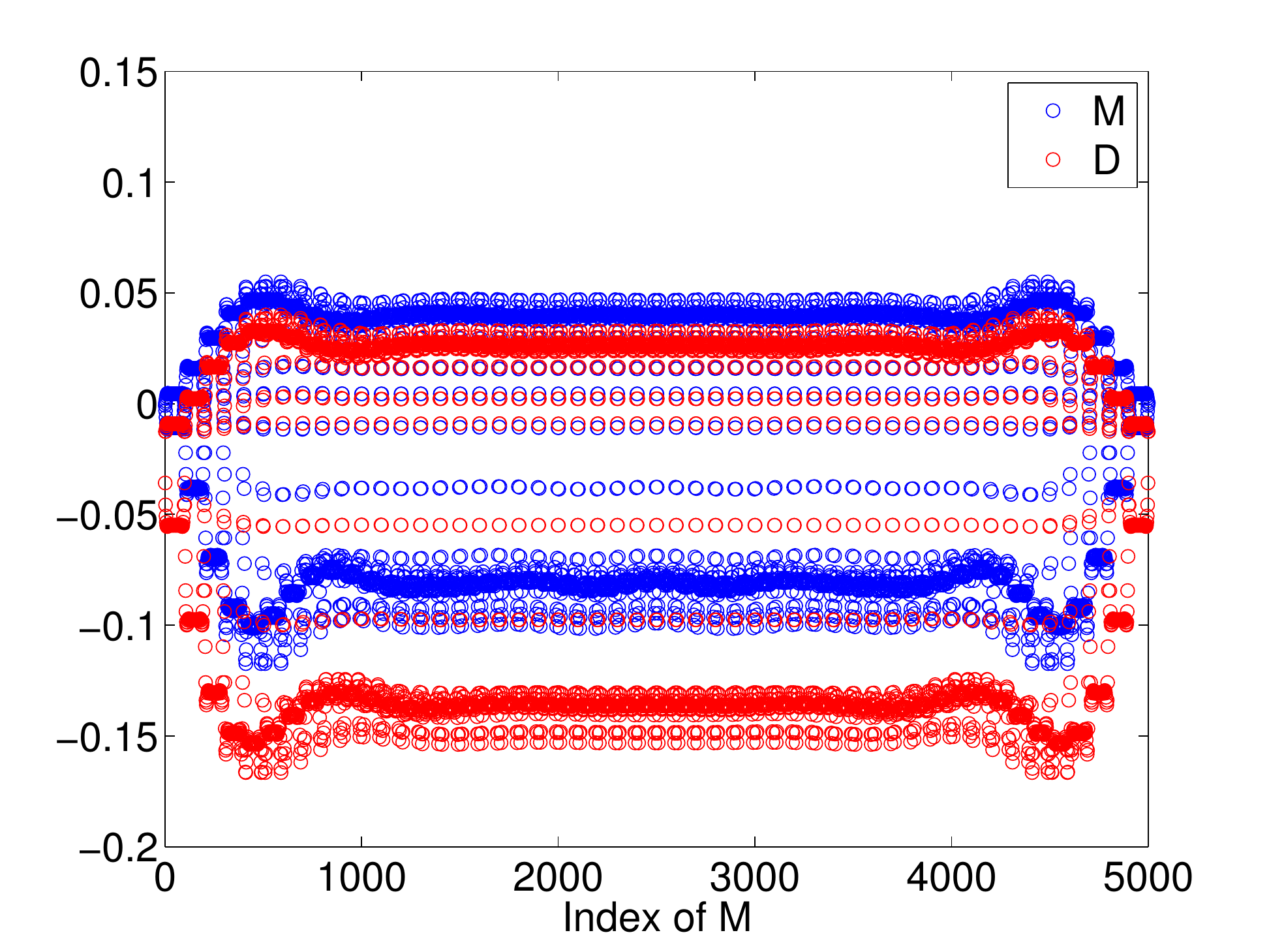}\label{fig:RDB5000 with M and D}}                
  \subfigure[Sample points for fitting]{\includegraphics[width=0.45\textwidth]{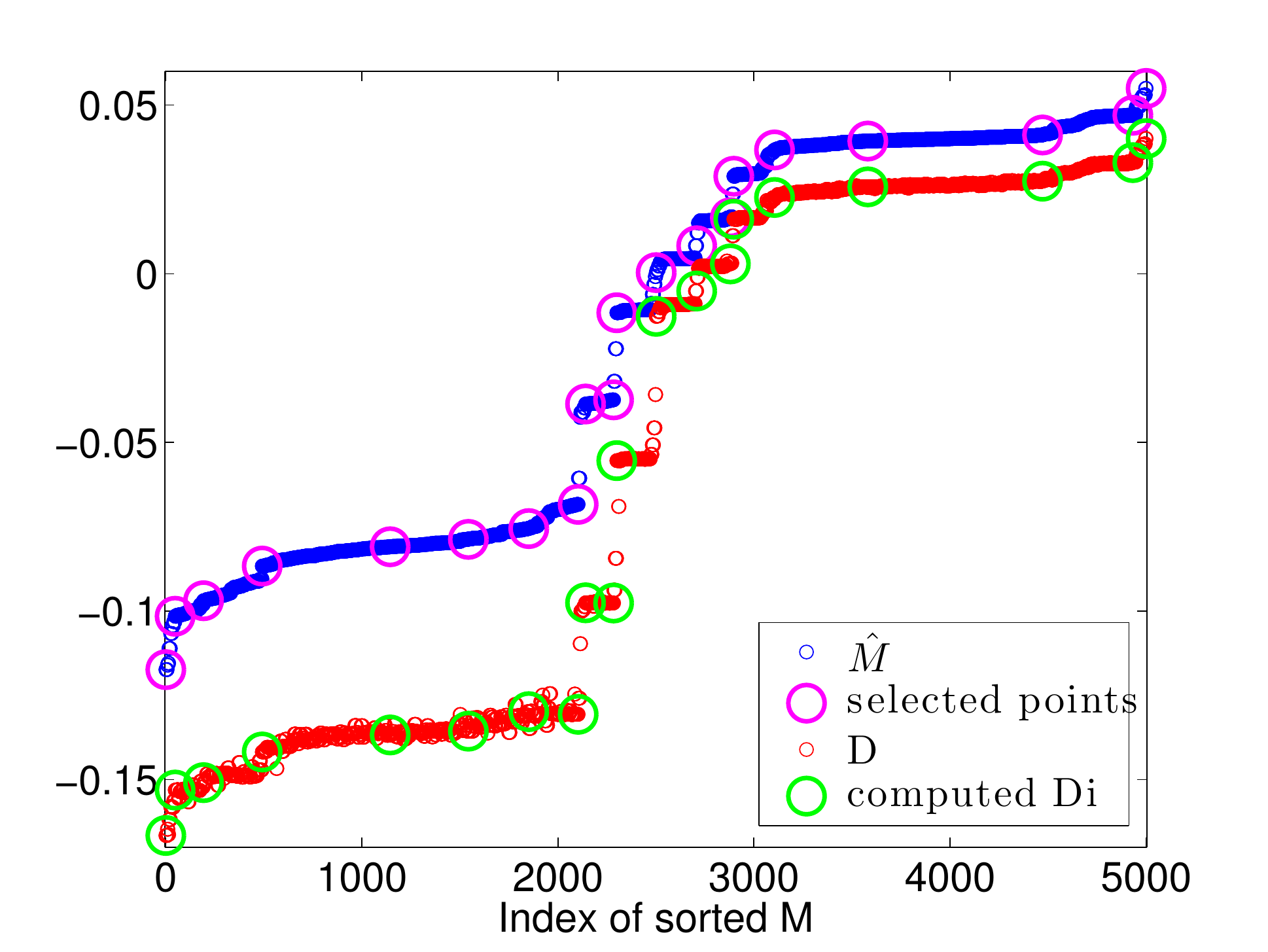}\label{fig:RDB5000 with sample points}}   
  \caption{A typical example to show our sampling strategy based on the pattern correlation of $M$ and $D$. $M$ is computed by 20 singular vectors of matrix RDB5000. In the right figure, the magenta and green circles denote the sample points associated with the sampling indices $S_{fit}$ in $\hat{M}$ and $D$.}
  \label{fig:RDB5000-Sampling}
\end{figure}

A typical sampling result produced by our method is shown in Figure \ref{fig:RDB5000-Sampling}.
The graph on the left shows the diagonals $M$ and $D$ plotted in their 
  original order.
The pattern correlation between them is clear, but it is less straightforward how to 
  pick the fitting points. 
Figure \ref{fig:RDB5000 with sample points} shows both diagonals plotted with 
the order of indices that sorts $\hat M$, i.e., $M(J), D(J)$.
The points picked by our algorithm based on $\hat M$ correspond to an 
  almost monotonic sequence of points in $D(J)$. The associated indices can then be used to 
  identify a suitable set of entries in $D$ to perform the numerical integration.

%---------------------------------------------------------------
\subsection{Two Fitting Models}
\label{section:two fitting models}

Next we construct a fitting model that minimizes 
  $\|f(M(S_{fit})) - D(S_{fit})\|$.
The fitting model must have sufficient predictive power with only a 
  small number of points and it should avoid oscillating behavior 
since we work on the monotonic sequence $\hat M$ which we assume correlates well with $\hat D$. 
Therefore, we consider a linear model and a piecewise polynomial model. 

\begin{figure}[!h]
  \centering
  \subfigure[Linear LS model in original order]{\includegraphics[width=0.35\textwidth]{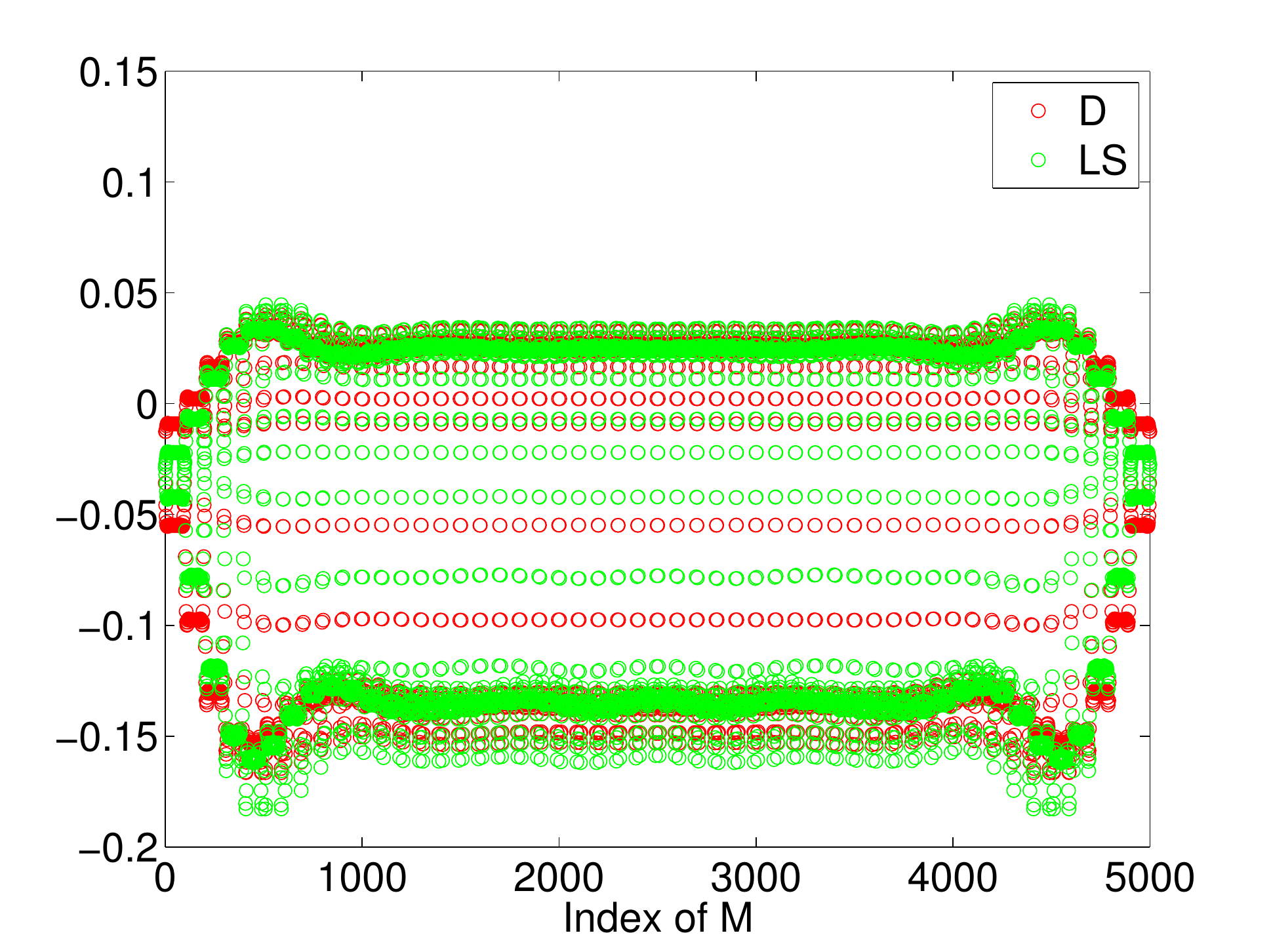}\label{fig:Linear LS model in original order}}                
  \subfigure[PCHIP model in sorted order]{\includegraphics[width=0.35\textwidth]{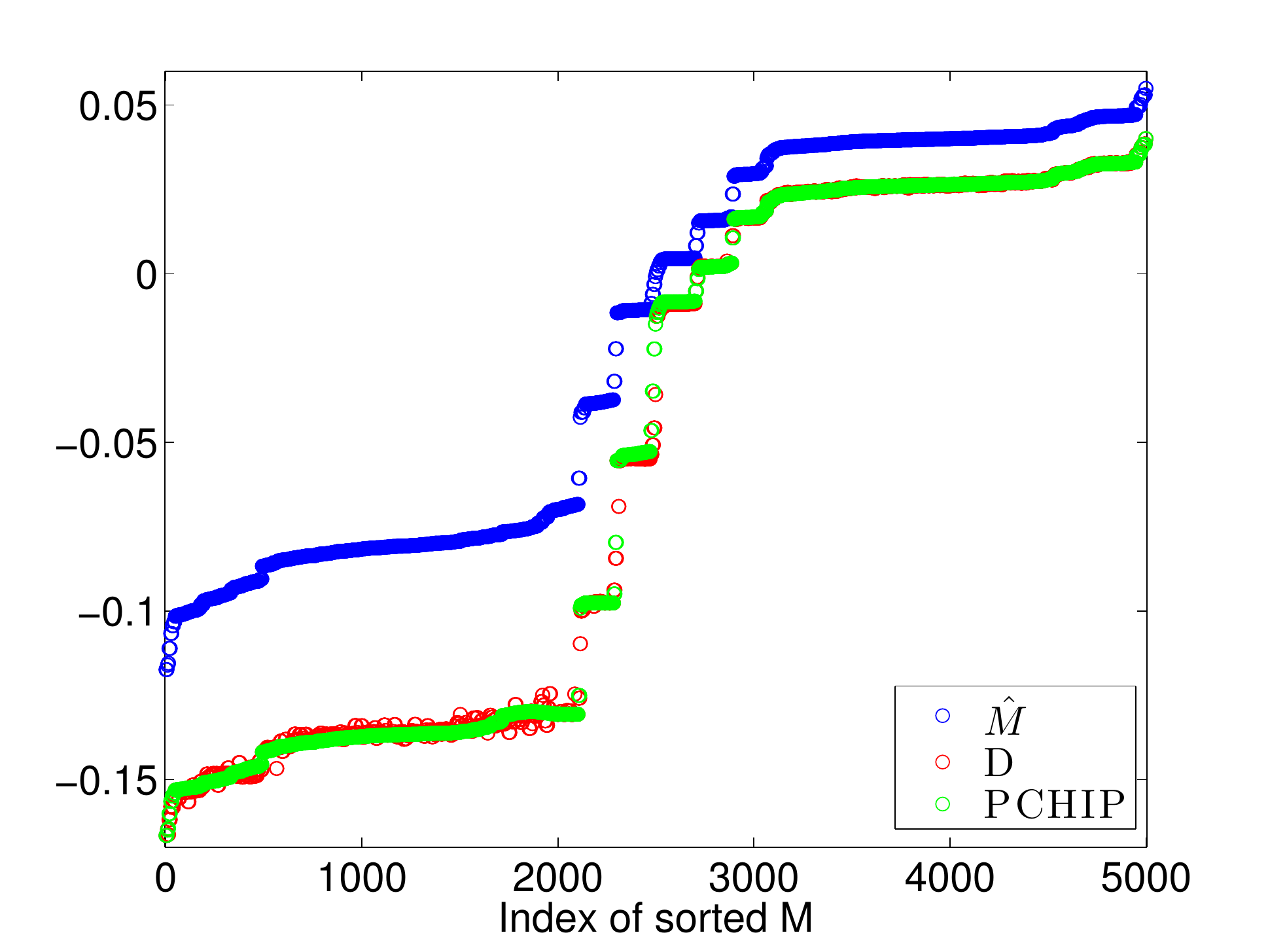}\label{fig:PCHIP model in sorted order}}   
  \caption{Fitting results of the matrix RDB5000 in original order and sorted order with linear LS model and PCHIP model.}
  \label{fig:RDB5000-Fitting results in original order and sorted order}
\end{figure}

The MC methods in (\ref{eq: approximation}) and (\ref{eq:refine approximation}) can resolve the trace when $D = M + c$ while importance sampling can resolve the trace when $D = cM$. 
To combine these, we first use a linear model, $ y = bM + c.  $
We determine the parameters $b, c$ by a least squares fitting,
$
{\mathrm{argmin}}_{b, c \in \mathbb{R}}  {\|D(S_{fit}) - (bM(S_{fit}) + c) \|_2}.
$
The linear model may be simple but avoids the large oscillations of higher 
  degree polynomials, and in many cases it is quite effective in improving the 
  accuracy of the trace estimation and reducing the variance of the diagonal elements 
  of $E_{fit}$.
The linear fitting algorithm is described in Algorithm \ref{alg:LS fitting}. 
Figure \ref{fig:Linear LS model in original order} shows 
  the fitting result on the example matrix of the previous section, 
  in the original order of $D$.

\begin{algorithm} 
\caption{Linear least squares fitting model for approximating $Tr(A^{-1})$.}
\begin{algorithmic}[1]
\STATEx ${\bf Input:}\ A\in \mathbb{R}^{N\times N}$
\STATEx ${\bf Output:}\ Tr(A^{-1})$ estimation: $T_f$
\STATE Compute $M$ using ILU or Eigendecompostion or SVD on $A$
\STATE Call Algorithm \ref{alg:point identification algorithm based on the trapezoidal rule} to compute sample set $S_{fit}$.
\STATE Find $[b, c] = {\mathrm{argmin}} {\|D( S_{fit}) - (bM( S_{fit}) + c) \|_2}$.
\STATE Compute trace approximation $T_f = \sum_{i=1}^N (b*M(i) + c)$
\end{algorithmic}
\label{alg:LS fitting}
\end{algorithm}

The linear model preserves the shape of $M$, and therefore relies exclusively on the quality of $M$. 
To take advantage of our premise that the distribution $\hat M$ approximates well the distribution of $\hat D$,
  our next fitting model is the Piecewise Cubic Hermite Spline Interpolation (PCHIP). 
It was proposed in \cite{fritsch1980monotone} to construct a visually pleasing monotone piecewise cubic interpolant to monotone data. The PCHIP interpolant is only affected locally by changes in the data and, most importantly, it preserves the shape of the data and respects monotonicity.
Therefore, we work on $\hat M$ and the indices 
  $\hat S_{fit} = [1=s_1, s_2, \ldots, s_{k-1}, s_k=N]$ 
  which are given in an order such that
  $\alpha = \hat{M}({s_1}) \leq \hat{M}({s_2}) \leq \cdots \leq \hat{M}({s_k}) = \beta $ 
  is a partition of the interval $[\alpha,\beta]$. 
An index $s_i$ corresponds to the index $J^{-1}(s_i)$ in the original ordering of $M$, 
  where $J$ is from 
  Algorithm \ref{alg:point identification algorithm based on the trapezoidal rule} and $J^{-1}$ denotes the inverse mapping of the sorted list.
Thus, for each $s_i$ we compute $D(J^{-1}(s_i)), i=1,\ldots ,k$, and we
  use PCHIP to construct a piecewise cubic function such that,
\begin{equation}
p(\hat{M}(s_i)) = D(J^{-1}(s_i)), \quad i = 1, 2, \ldots, k.
\end{equation}
Notice that $p(x)$ will be monotone in the subintervals 
  where the fitting points $D(J^{-1}(s_i))$ are also monotone. 
  Therefore, as long as $M$ is close to $D$ in the sense of (\ref{eq:bound}),
  integration of $p(x)$ will be very accurate.

%\begin{equation}
%p(x) = f_i H_1(x) + f_{i+1} H_2(x) + d_i H_3(x) + d_{i+1} H_4(x) .
%\end{equation}
%where $H_i(x)$ are the cubic Hermite basis functions, and the values $d_i$ and $d_{i+1}$ are the derivatives $p^{'}(\hat{M}_i)$ and $p^{'}(\hat{M}_{i+1})$ respectively. 

The PCHIP model is given in Algorithm \ref{alg:PCHIP fitting}.
The first two steps are the same as in Algorithm \ref{alg:LS fitting}. 
In step 3, we apply the function {\bf unique} to remove the duplicate elements 
  of $\hat M(\hat S_{fit}))$ to produce a sequence of unique values as required by PCHIP. 
This yields a subset of the indices, $\hat S_{fit}'$, which is mapped to original 
  indices as $I=J^{-1}(\hat S_{fit}')$ to be used in PCHIP.

\begin{algorithm}[h!]
\caption{PCHIP fitting model for approximating $Tr(A^{-1})$.}
\begin{algorithmic}[1]
\setcounter{ALG@line}{2}
\STATEx ${\bf Input:}\ A\in \mathbb{R}^{N\times N}$
\STATEx ${\bf Output:}\ Tr(A^{-1})$ estimation: $T_f$
	\STATEx \% Steps 1 and 2 are the same as in Algorithm \ref{alg:LS fitting}
	\STATE Remove duplicates: 
	$\hat S_{fit}' = \bold{unique}(\hat M(\hat S_{fit}))$, $I =J^{-1}(\hat S_{fit}')$
	\STATE Apply PCHIP to fit $p({M(I)}) = D(I)$ and obtain a polynomial $p(M) \approx D$
	\STATE Compute trace approximation $T_f = \sum_{i=1}^N p(M(i))$
\end{algorithmic}
\label{alg:PCHIP fitting}
\end{algorithm}

Figure \ref{fig:PCHIP model in sorted order} shows the result of fitting 
  on the RDB5000 example of the previous figures.
Table \ref{ta:Comparision of trace estimation, relative trace error and variances between two fitting models} compares the relative error for the trace,
as well as the variances of $D$, $D - M$ and $D - p(M)$ when using the linear LS and PCHIP models on two test matrices, OLM5000 and KUU.
In both cases, the two models can provide a trace estimate of surprising relative accuracy of the order of {\tt 1e-2} with only 20 fitting points. 
In addition, for matrix OLM5000, the standard deviation of MC on $D - p(M)$ is reduced 
   by a factor of 10 compared to that of MC on $D$ (a speedup of 100 in terms of samples). 
However, for matrix KUU, the standard deviation of MC on $D - p(M)$ does not improve much over MC on $D$.
This reflects the discussion in Section \ref{sec:pointIdentification}, and
  suggests that our method can serve as a standalone kernel for 
  giving a fast trace estimate.
In addition, in some cases the method can reduce the variance significantly
  to accelerate a second stage MC.

\begin{table}[htbp]
\caption{Comparing the trace relative error and variances in (\ref{eq: unit vector variance of invA}), (\ref{eq: unit vector variance of E}) and (\ref{eq:unit vector variance of E_fit}) for matrices OLM5000 and KUU, using the linear LS and PCHIP models with 20 fitting points each. The traces of matrix
OLM5000 and KUU are $Tr(A^{-1})=-5.0848e+02$ and $3.6187e+03$, 
respectively.}

\centering
\small
\begin{center}
	\begin{tabular}{|r|rrr|rrr|}
	\hline
	\multicolumn{1}{|r}{Matrix}
	 & \multicolumn{3}{c}{{\tt OLM5000}} 
	 & \multicolumn{3}{c|}{{\tt KUU}} \\ \hline
 	Model & M & LS & PCHIP & M & LS & PCHIP  \\ \hline
	%{\footnotesize $Trace$ } 
	%& \Red{-5.0848e+02} & -5.0848e+02 & -5.0848e+02 & \Red{3.6187e+03} & 3.6187e+03 & 3.6187e+03\\ \hline 
	%{\footnotesize $TraceEst$ } 
	%& \Red{-1.7033e+01} & -4.9713e+02 & -5.0122e+02 & \Red{1.5084e+03} & 3.5745e+03 & 3.5663e+03\\ \hline
	{\footnotesize $Relative\ error$ } 
	& 9.6650e-01 & 2.2320e-02 & 1.4288e-02 & 5.8316e-01 & 1.2207e-02 & 1.4469e-02\\ \hline	  	
	{\footnotesize $Std(T_{e_i}(A^{-1}))$ }
	& 1.1425e+02 & 1.1425e+02 & 1.1425e+02 & 2.8731e+02 & 2.8731e+02 & 2.8731e+02\\ \hline	
	{\footnotesize $Std(T_{e_i}(D-M))$ }
	& 1.1007e+02 & 1.1007e+02 & 1.1007e+02 & 1.4999e+02 & 1.4999e+02 & 1.4999e+02\\ \hline
	{\footnotesize $Std(T_{e_i}(D - p(M)))$ }
	& -- & 2.0332e+01 & 1.7252e+01 & -- & 1.6137e+02 & 1.6289e+02\\ \hline	
	\end{tabular}
\end{center}
\label{ta:Comparision of trace estimation, relative trace error and variances between two fitting models}
\end{table}

%---------------------------------------------------------------
\section{Dynamic Evaluation of Variance and Relative Trace Error}
\label{sec:DynamicEvaluation}
Since there are no a-posteriori bounds for the accuracy of our results,
  we develop methods that use the information from the solution of
  the linear systems to incrementally estimate the trace error and 
  the variances of the resulting approximations.
This approach is also useful when $M$ is updated with more left and right 
  eigenvectors or singular vectors obtained from the solution of additional 
  linear systems.

%---------------------------------------------------------------
\subsection{Dynamic Variance Evaluation}
To decide which MC method we should use after the fitting stage or even whether it 
  is beneficial to use the fitting process for variance reduction, we 
  monitor incrementally the variances
$Var(T_{e_i}(A^{-1}))$, $Var(T_{e_i}(E_{fit}))$, 
$Var(T_{Z_2}(A^{-1}))$,  and $Var(T_{Z_2}(E))$,
with the aid of the cross-validation technique \cite{harrell2001regression}. 

Our training set is the fitting sample set $D(S_{fit}),$ 
  while our test set $D(S_{mc})$ is a small random set which is independent of 
  the fitting sample set.
If we want to combine our method with MC, eventually more samples need to be computed, 
  and thus we can pre-compute a certain number of them as the test set $D(S_{mc})$.
We have used the holdout method \cite{arlot2010survey}, a single train-and-test experiment for some data splitting strategy since the fitting sample set is fixed. 

To compute $D(S_{fit})$ or $D(S_{mc})$ , a linear system with multiple right hand sides is solved as follows:
\begin{equation} 
 A_{ii}^{-1} = e_{i}^{T}x_{i},\ \ \ Ax_{i} = e_{i}, \ \forall i \in S_{fit} \cup S_{mc}.
\end{equation}
The computed column vectors $x_{i}$ can be used to estimate the Frobenius norm of both $A^{-1}$ and $E = A^{-1} - Z^{-1}$ \cite{gudmundsson1995small,kenney1998statistical}. Then $Var(T_{Z_2}(A^{-1}))$ and $Var(T_{Z_2}(E))$ can be estimated as follows:
\begin{equation} \label{eq:estimation of Hutchinson variance of invA}
Var(T_{Z_2}(A^{-1})) \approx \dfrac{2N}{s^2} \sum_{i}{(\|x_i\|^{2} - |D_i|^{2}}),  \ \forall i \in S_{fit} \cup S_{mc},
\end{equation}
\begin{equation} \label{eq:estimation of Hutchinson variance of E}
Var(T_{Z_2}(E)) \approx \dfrac{2N}{s^2} \sum_{i}{(\|E(:,i)\|^{2} - |E(i,i)|^{2}}),  \ \forall i \in S_{fit} \cup S_{mc} , 
\end{equation}
where $E(:,i) = Ee_i$. Simultaneously, based on the sampled diagonal elements $A_{ii}^{-1}$, we can also update the evaluation of $Var(T_{e_i}(E_{fit}))$, $Var(T_{e_i}(E))$ and $Var(T_{e_i}(A^{-1}))$. Here we only show the computation of the unbiased variance estimation for $Var(T_{e_i}(E_{fit}))$ by:
\begin{equation} \label{eq:estimation of unit vector variance of E_fit}
Var(T_{e_i}(E_{fit})) \approx \dfrac{N^2}{s-1} Var(E_{fit}(S_{mc})) . 
\end{equation}

Note that $S_{fit}$ should not be used for estimating the variance of the unit vector MC estimator since these sample points are exact roots of the PCHIP function. 

\begin{algorithm}[h!]
\caption{Dynamic variance evaluation algorithm for estimating variances of different MC methods \label{alg:Dynamic variance evaluation}}
\begin{algorithmic}[1]
\STATEx ${\bf Input:}\ A\in \mathbb{R}^{N\times N}$
\STATEx ${\bf Output:}\ Tr(A^{-1})$ estimation and variances estimations of various MC methods
\STATE Initialize $maxPts$, $S_{fit}$, $S_{mc}$
\IF{$M$ is computed using ILU on A }
\STATE Compute the approximation $M$ 
\ENDIF 
\STATE Generate random index set $S_{mc}$ without replacement and compute MC samples $D(S_{mc})$
\FOR{$i=5:1:maxPts$ }
\IF{$M$ is computed using Eigendecompostion or SVD on A }
\STATE Update $M$ with $2*i$ number of left and right eigenpairs or singular triplets
\ENDIF 
\STATE Call Algorithm \ref{alg:point identification algorithm based on the trapezoidal rule} 
  to find more indices so that $S_{fit}$ has $i$ fitting points
\STATE Call Algorithms \ref{alg:LS fitting} or \ref{alg:PCHIP fitting} to update approximation of $Tr(A^{-1})$ 
\STATE Estimate variances of different MC methods based on (\ref{eq:estimation of Hutchinson variance of invA}), (\ref{eq:estimation of Hutchinson variance of E}) and (\ref{eq:estimation of unit vector variance of E_fit})
\ENDFOR
\end{algorithmic}
\end{algorithm}

We implement the dynamic variance evaluation scheme in Algorithm \ref{alg:Dynamic variance evaluation}. 
In lines 2-4 and 7-9, the approximation $M$ can be computed using an ILU factorization at the beginning 
of the procedure or be updated 
with increasing number of singular triplets or eigenpairs. Note that if $M$ is obtained by a partial eigendecomposition or SVD, 
the updated $M$ is different in two consecutive steps, and thus
Algorithm \ref{alg:point identification algorithm based on the trapezoidal rule}
will return a slightly different index set $S_{fit}$ which may not be incremental.
This may not provide a consistent improvement of the relative trace error 
 during the fitting progress.
In line 10 of the algorithm, we force the points to be incremental between 
  steps $i$ and $i+1$ as follows; we generate the entire set
  $S_{fit}^{(i+1)}$ and remove the indices that lie the closest to the previous 
  index set $S_{fit}^{(i)}$.
The remaining index set is incorporated into $S_{fit}^{(i)}$. 
This simple scheme works quite well experimentally.  

\begin{figure} [h!]
  \centering
  \subfigure[Dynamic variance estimation with ILU]{\includegraphics[width=0.35\textwidth]{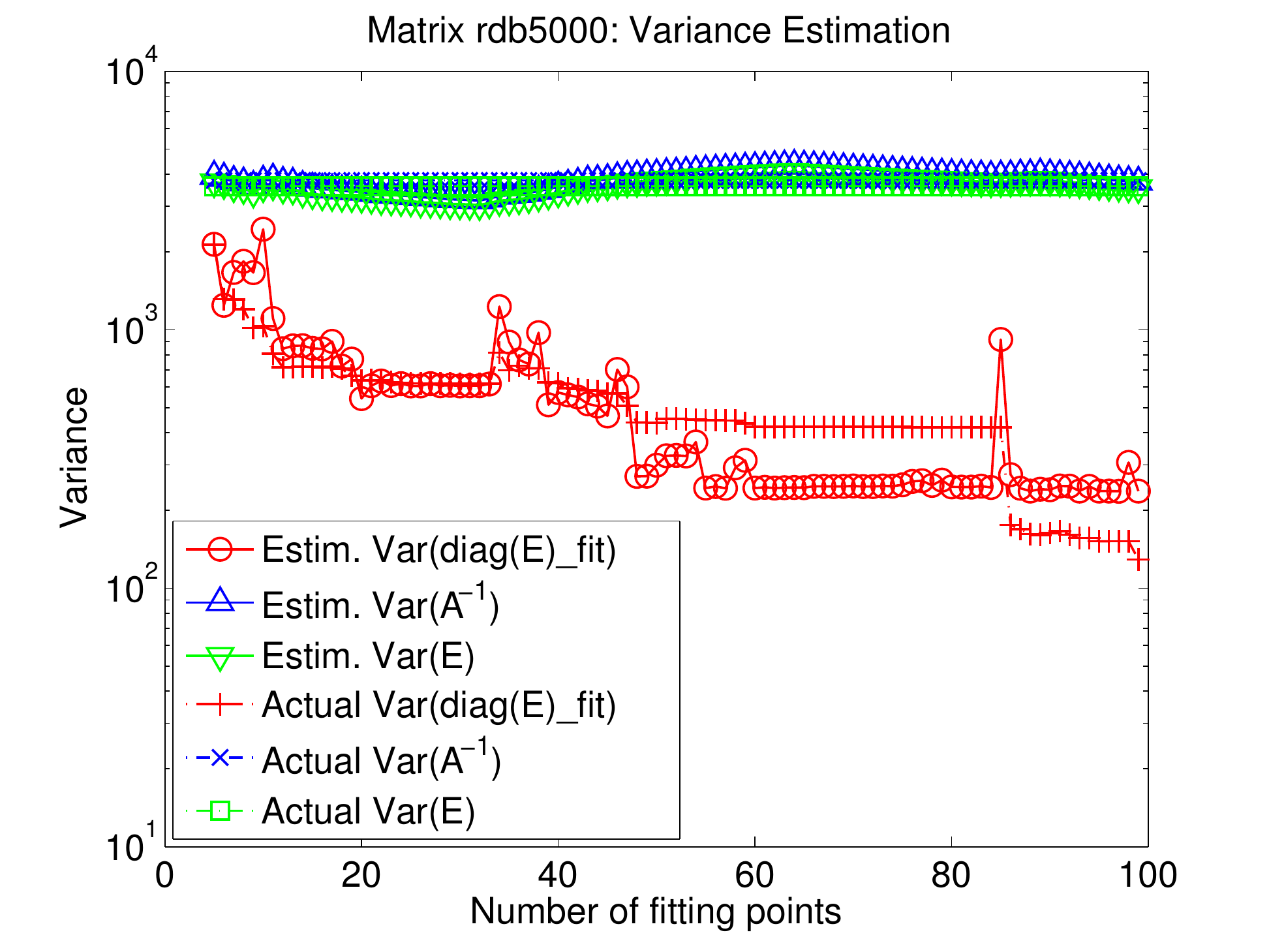}\label{fig:estimate variance with ILU on rdb5000}}                
  \subfigure[Dynamic variance estimation with SVD]{\includegraphics[width=0.35\textwidth]{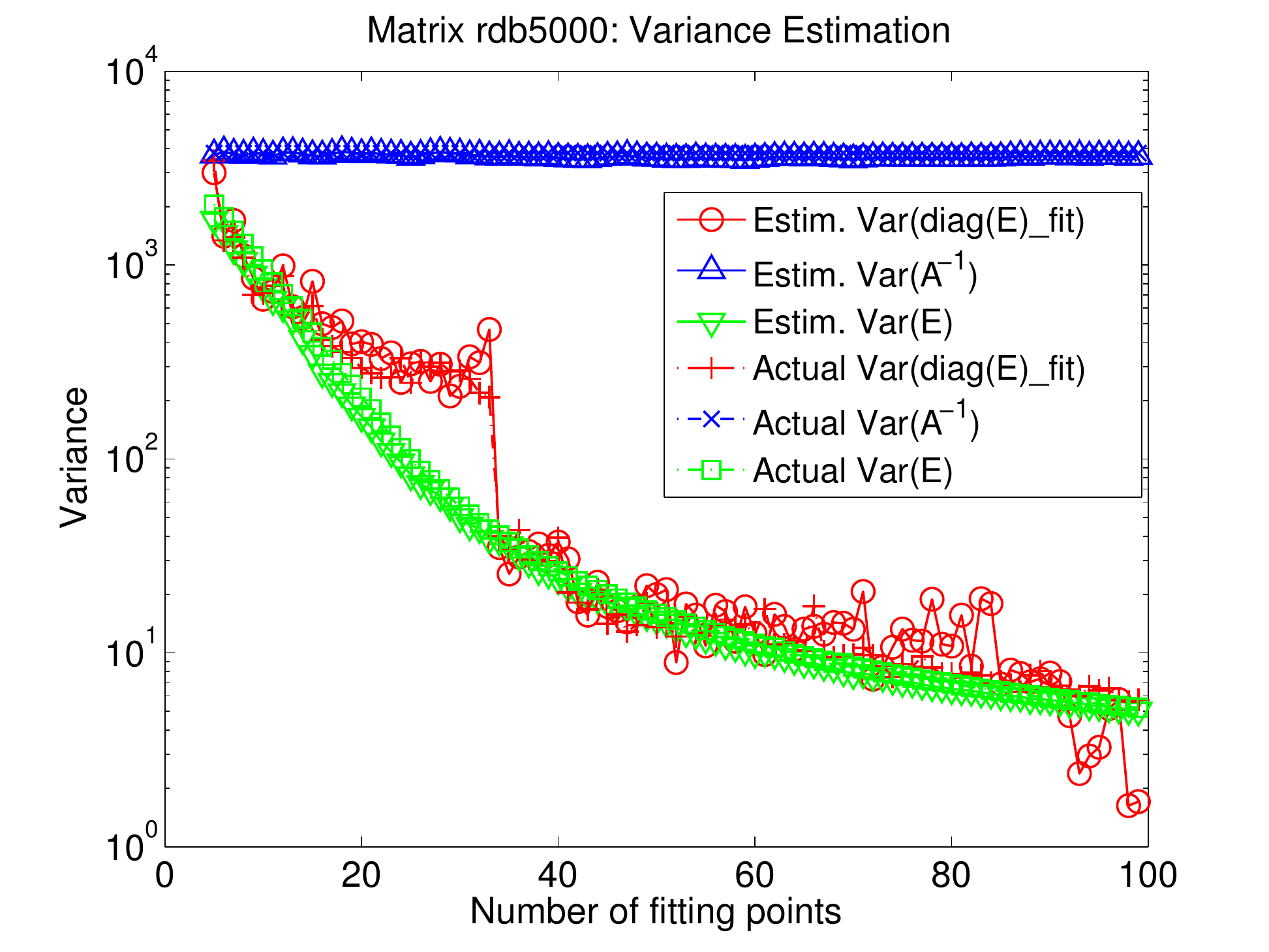}\label{fig:estimate variance with SVD on rdb5000}}   
  \caption{Comparing estimated variances and actual variances of unit vector on $E_{fit}$ (denoted as $Var(diag(E)\_{fit})$) and Rademacher vector on $A^{-1}$ (denoted as $Var(A^{-1})$) and $E$ (denoted as $Var(E)$) of the matrix RDB5000 with ILU and SVD respectively.}
  \label{fig:RDB5000-Dynamic variance estimation results with ILU and SVD}
\end{figure}

Figure \ref{fig:RDB5000-Dynamic variance estimation results with ILU and SVD} 
  shows how the actual and the estimated variances for three MC methods match 
  for the test matrix RDB5000.
In addition, the relative difference between different MC methods becomes clear after 
  only a few points which facilitates not only the proper choice of MC method but 
  also an early decision to stop if further fitting is not beneficial.
In the numerical experiments section we show that these results are 
typical for matrices from a wide variety of applications.

%---------------------------------------------------------------
\subsection{Monitoring Relative Trace Error}
\label{sec:Monitoring Relative Trace Error}
\begin{figure} 
  \centering
  \subfigure[Comparing $D$ and $p(M)$ in order of $\hat{M}$]{\includegraphics[width=0.35\textwidth]{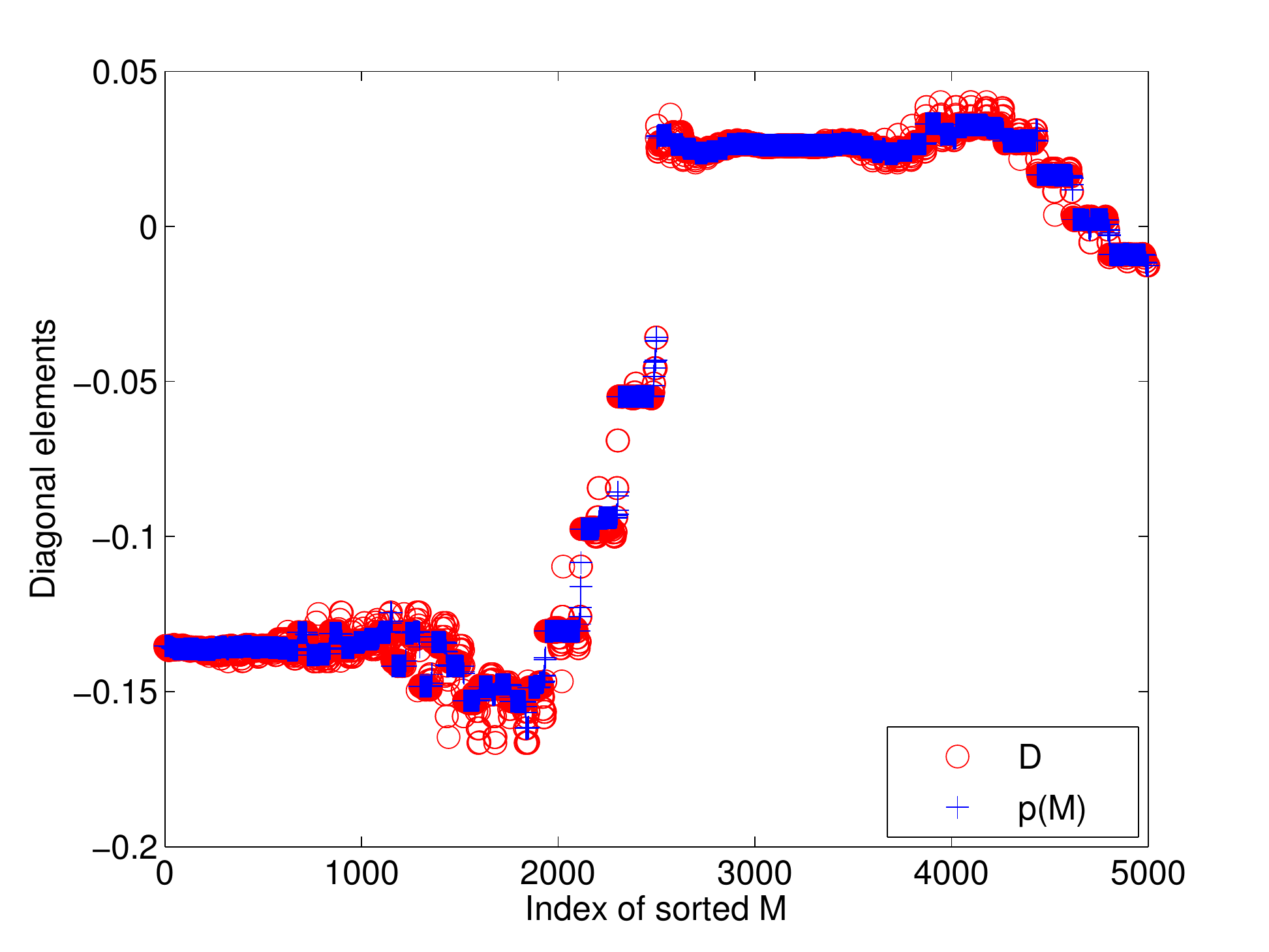}\label{fig:compare diagonals of D and p(M) in original order}}                
  \subfigure[Comparing sorted $D$ and sorted $p(M)$]{\includegraphics[width=0.35\textwidth]{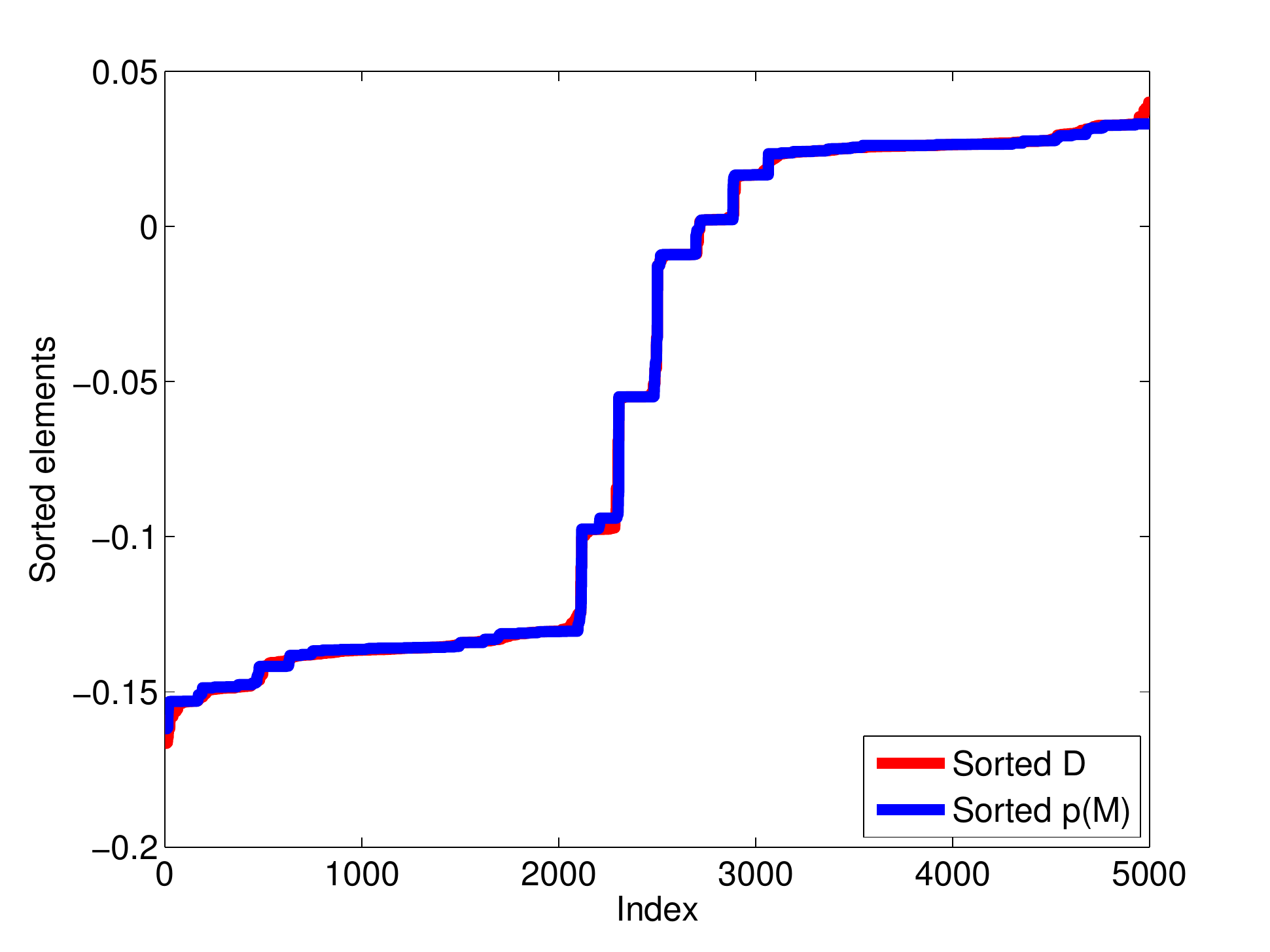}\label{fig:compare diagonals of D and p(M) in sorted order}}   
  \caption{Compare $D$ and $p(M)$ of the matrix RDB5000 in different order where $M$ is computed by ILU.}
  \label{fig:RDB5000-compare the signal of D and p(M)}
\end{figure}

As discussed in Section \ref{sec:pointIdentification} and showed in 
  Table \ref{ta:Comparision of trace estimation, relative trace error and variances between two fitting models},
the variance of MC on the vector $E_{fit}$ may not be reduced if 
the sorting permutations of $\hat M$ and $\hat D$ are dissimilar.
Even in such cases, however, our fitting method can provide a good 
  trace estimation with only a small number of samples. 
We further investigate this by comparing the elements of 
  $D$ and $p(M)$ with different orderings.
Figure \ref{fig:compare diagonals of D and p(M) in original order} shows the elements of $p(M(J))$ and $D(J)$, i.e., with respect to the order of $\hat{M}$. It illustrates that although $p(M)$ captures the pattern of $D$, the order of 
  its elements does not correspond exactly to that of $D$; hence the small reduction in 
  $Var(T_{e_i}(E_{fit}))$. However, Figure \ref{fig:compare diagonals of D and p(M) in sorted order} reveals that the distributions of the sorted $p(M)$ and the sorted $D$ almost coincide; hence, the two integrals $Tr(p(M))$ and $Tr(D)$ are very close.

The main obstacle for using our method as a standalone kernel for trace 
  estimation
  is that there is no known way to measure or bound the relative trace error.
Resorting to the confidence interval computed by the variance of a MC estimator is 
  pessimistic as our results show (see also \cite{avron2011randomized, roosta2014improved}).
We note the similarity to the much smaller error obtained in the average case 
  of integrating monotonic functions with adaptive bilinear forms versus the worst 
  case known bounds \cite{Kiefer_optimality_1957,Novak_IntegrMonotoneF}.

Motivated by these previous research results, we show how to 
 develop practical criteria to monitor the relative trace error in
  our Algorithms. Suppose at each step of the fitting process in Algorithm \ref{alg:Dynamic variance evaluation} we collect a sequence of trace estimates $T_i, i \in [1,maxPts]$. Consider the trace estimations in two successive steps,
\begin{equation} \label{eq:relative error}
\dfrac{|T_i - T_{i+1}|}{|T_{i+1}|} = \dfrac{|(T_i - Tr(D)) - (T_{i+1} - Tr(D))|}{|T_{i+1}|} = \dfrac{|E_i - E_{i+1}|}{|T_{i+1}|} \leq \dfrac{2max(|E_i|, |E_{i+1}|)}{|T_{i+1}|} \approx \dfrac{2max(|E_i|, |E_{i+1}|)}{|Tr(D)|}.
\end{equation} 
As long as $T_i$ converges with more fitting points, the relative difference of two 
  successive trace estimations can serve as an approximation to the relative error. 
However, when the global pattern provided by $M$ and $p(M)$ is not fully 
  matched to that of $D$, convergence of $T_i$ stagnates until enough points 
  have been added to resolve the various local patterns.
To determine whether the current relative trace error estimation can be trusted,
  we present our second heuristic by considering the error bound of our fitting models. 

When approximating $f(\hat M)$ with $p(\hat M)$, the PCHIP Hermite cubic splines with 
  $k$ points on the interval $[\alpha, \beta]$, 
  the bound on the error $E(\hat{M}) = f(\hat{M}) - p(\hat{M})$  is
 given by \cite{lucas1974error},
\begin{equation} \label{eq:error bound}
|E(\hat{M})| \leq \dfrac{1}{384} h^4 \|f^{(4)}\|_{\infty,[\alpha,\beta]} \ ,
\end{equation}
where $h = \tfrac{\beta - \alpha}{k}$, and
$\|f^{(4)}\|_{\infty,[\alpha,\beta]}$ denotes the maximum value of the 
  fourth derivative of $f$ in the entire interval $[\alpha, \beta]$.
Since $\|f^{(4)}\|_{\infty,[\alpha,\beta]}/384$ is a constant, 
  in two successive fitting steps we have,
\begin{equation} \label{eq:absolute error}
\dfrac{|E_{i}(\hat{M})|}{|E_{i+1}(\hat{M})|} \approx \dfrac{h_{i}^4 }{h_{i+1}^4} = \dfrac{((\beta - \alpha)/k_{i})^4 }{((\beta - \alpha)/k_{i+1})^4} = (\dfrac{k_{i+1}}{k_{i}})^4,
\end{equation}
where $k_{i}$ and $k_{i+1}$ are the number of fitting points in two consecutive steps. 
We can use (\ref{eq:absolute error}) to estimate the maximum possible improvement between two consecutive trace errors. If the $i+1$ trace error estimate reduces over the $i$-th estimate
by a factor of more than $(k_{i}/k_{i+1})^4$, we do not trust it.
 
\begin{algorithm} 
\caption{Dynamic relative trace error evaluation algorithm during the fitting process} 
\label{alg:monitoring trace error}
\begin{algorithmic}[1]
\STATEx ${\bf Input:}\ TraceFit$ from Alg. \ref{alg:Dynamic variance evaluation}
\STATEx ${\bf Output:}\ TraceErr$ estimation 
\STATEx \% $TraceErr(i)$ is defined only for $i>5$ and we assume TraceErr(6) is well defined
\IF{$i == 6$}
	\STATE $TraceErr(6) = |TraceFit(6) - TraceFit(5)|/|TraceFit(5)|$	
\ENDIF
\IF{$i > 6$}
	\STATE $TempTraceErr = |TraceFit(i) - TraceFit(i-1)|/|TraceFit(i)|$	
\IF{$(TempTraceErr/TraceErr(i - 1) \geq ((i-1)/i)^4)$}
	\STATE $TraceErr(i) = TempTraceErr$
\ELSE
	\STATE $TraceErr(i) = TraceErr(i - 1) * ((i-1)/i)^{9/4}$
\ENDIF
\ENDIF
\end{algorithmic}
\end{algorithm}

One caveat is that (\ref{eq:error bound}) may not be tight since the same bound
  holds for each subinterval $[\hat{M}(s_j),\hat{M}(s_{j+1})]$.
This means that a high derivative in one subinterval might dominate the bound 
  in (\ref{eq:error bound}) but should not affect the error in other intervals. 
Therefore, convergence might not be fully dictated by (\ref{eq:absolute error}).
In practice, we found that the improvement ratio is between 
  $O(k_i/k_{i+1})$ and $O((k_i/k_{i+1})^4)$. 
Therefore, if the current relative trace error estimate is determined not to be trusted, 
  we may instead use $|E_{i+1}(\hat{M})| \approx (k_i/k_{i+1})^k |E_1(\hat{M})|, \ k \in [1,4]$. 
The choice of $k$ depends on the quality of $M$. 
In our experiments, we use the geometric mean of the four rates, yielding $k=9/4$.
Recall that the corresponding ratio in MC is $O(\sqrt{k_i/k_{i+1}})$, which is much slower than the proposed trace estimation method. 

\begin{figure}[hbtp]
  \centering
  \subfigure[Monitoring relative trace error with ILU]{\includegraphics[width=0.35\textwidth]{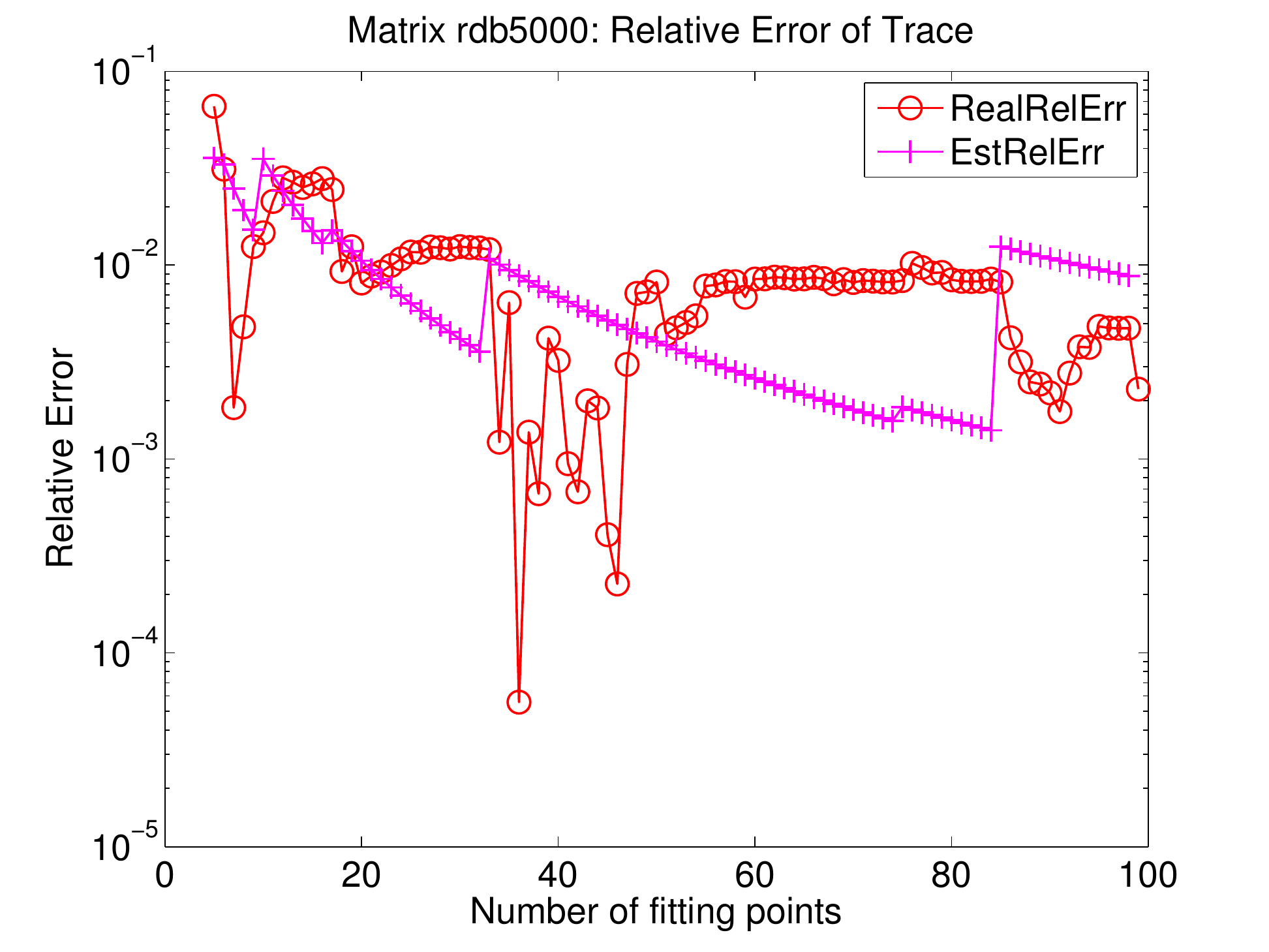}\label{fig:estimate relative trace error with ILU on rdb5000}}                
  \subfigure[Monitoring relative trace error with SVD]{\includegraphics[width=0.35\textwidth]{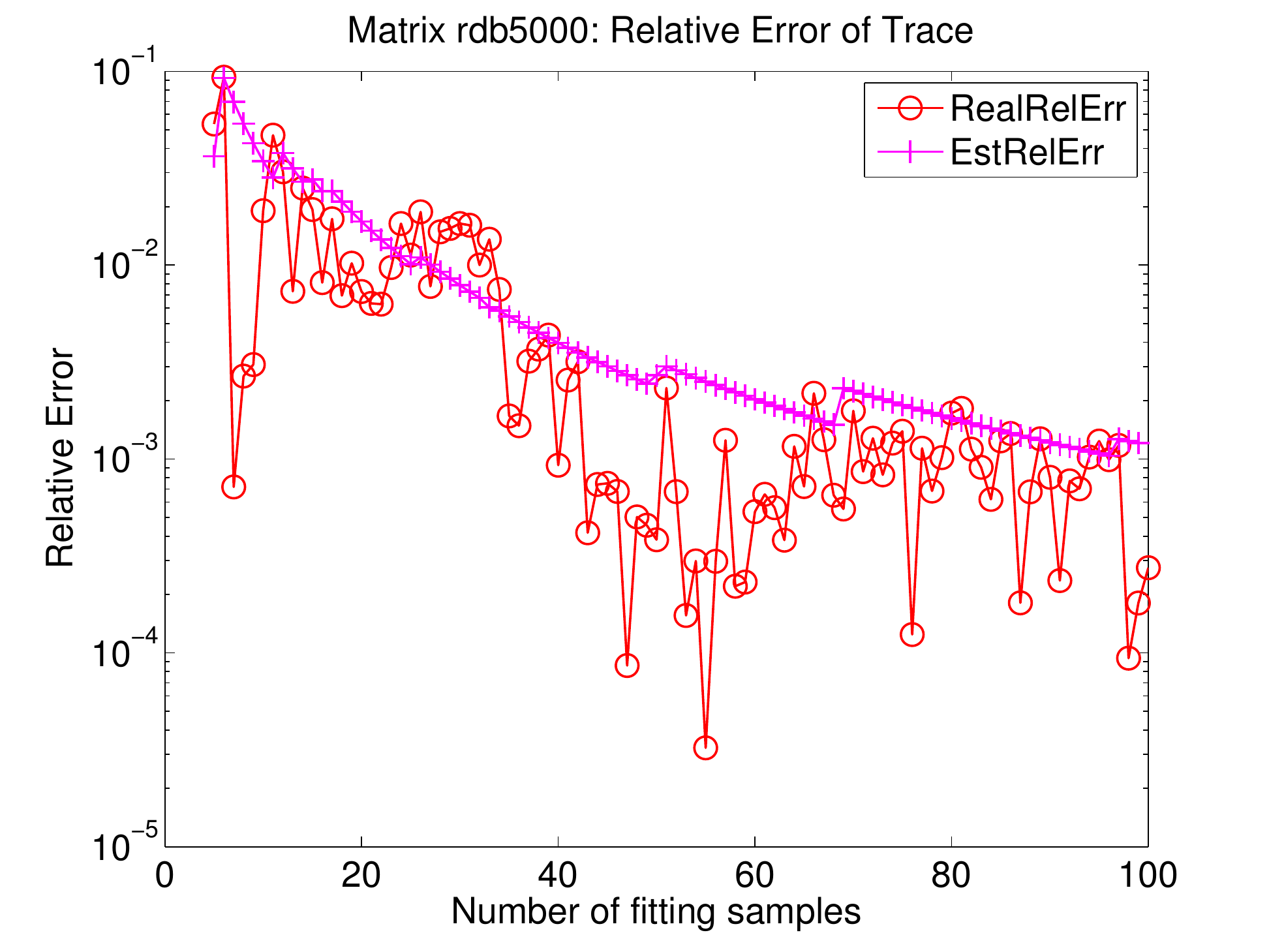}\label{fig:estimate relative trace error with SVD on rdb5000}}   
  \caption{Two examples of monitoring relative trace error of the matrix RDB5000 with ILU and SVD respectively.}
  \label{fig:RDB5000-Monitoring relative trace error with ILU and SVD}
\end{figure}

Algorithm \ref{alg:monitoring trace error} combines the two heuristics in 
  (\ref{eq:relative error}) and (\ref{eq:absolute error}) to dynamically monitor the 
  relative trace error.
It is called after step 12 of Algorithm \ref{alg:Dynamic variance evaluation}.
%In lines 2-4, we first compute the temporary relative trace error after we obtain the first two trace estimates. The idea behind lines 5 to 9 is that the temporary relative trace error below than the best improvement may be not trusted. We compute the geometric mean for estimating the current relative trace error if necessary. 
Figure \ref{fig:RDB5000-Monitoring relative trace error with ILU and SVD} shows two examples
  of how our dynamic method provides reasonable estimates of the relative trace error.

%---------------------------------------------------------------
\section{Numerical Experiments}
We run experiments on matrices that are sufficiently large to avoid misleading results 
due to sampling from small spaces but can still be inverted to obtain the exact trace error.
We select matrices RDB5000 and cfd1 from the University of Florida sparse matrix collection \cite{davis2011university} and generate three test matrices from applications that appear in \cite{bai1996some}.
The Heatflow160 matrix is from the discretization of the linear heat flow problem using the simplest implicit finite difference method.
The matrix Poisson150 is from 5-point central difference discretization of the 2D Poisson's equation on a square mesh. 
The VFH6 matrix is from the transverse vibration of a Vicsek fractal that is constructed self-similarly.
We also use matrix matb5 which is a discretization of the Wilson Dirac operator on a $8^4$ uniform lattice with 12 degrees of freedom at each node, using a mass near to critical. 

Table \ref{ta: test matrices} lists these matrices along with some of their basic properties.
All experiments are conducted using MATLAB 2013a. 
The number of fitting points increases as $s=5:100$.
The approximation $M$ is computed by ILU with parameters 
  {\tt type = ilutp} and {\tt droptol = 1E-2}, or 
  as a low rank approximation of $2s$ smallest singular vectors with accuracy {\tt 1E-6} (twice the number of 
  fitting points at each step),
  or by the bounds on the diagonal. 

\begin{table}[htbp]
\centering
\caption{Basic information of the test matrices }
\label{ta: test matrices}
\small
\begin{center} 
    \begin{tabular}{ |c|c|c|c|c|}
    \hline
    Matrix 		& Order & nnz(A) & $\kappa(A)$ & Application \\ \hline
%    OLM5000		& 5000  & 19996    &  	1		  & computational fluid   \\ \hline
    RDB5000 		& 5000  & 29600    &  	1.7E3	  & computational fluid  \\ \hline
    cfd1 		& 70656 & 1825580  &    1.8E7     & computational fluid  \\ \hline
%    SiNa       	& 5743  & 198787 &  1 		  & quantum chemistry \\ \hline
%    KUU 			& 7102  & 340200 &  1		  & structural problem \\ \hline
%    Wathen100 	& 30401  & 471601  & 8.2E3	  & random 2D/3D \\ \hline
    Heatflow160	& 25600  & 127360  & 2.6E0	  & linear heat flow   \\ \hline
    Poisson150 	& 22500  & 111900  & 1.3E4	  & computational fluid   \\ \hline
    VFH6       	& 15625  & 46873   & 7.2E1 	  & vicsek fractal \\ \hline
    matb5       	& 49152  & 2359296 & 8.2E4 	  & lattice QCD \\ \hline
    \end{tabular}
\end{center}  
\end{table}

%---------------------------------------------------------------
\subsection{Effectiveness of the fitting models}
\begin{figure}
  \centering
  \subfigure[Diagonal elements in (1,2000) ]{\includegraphics[width=0.32\textwidth]{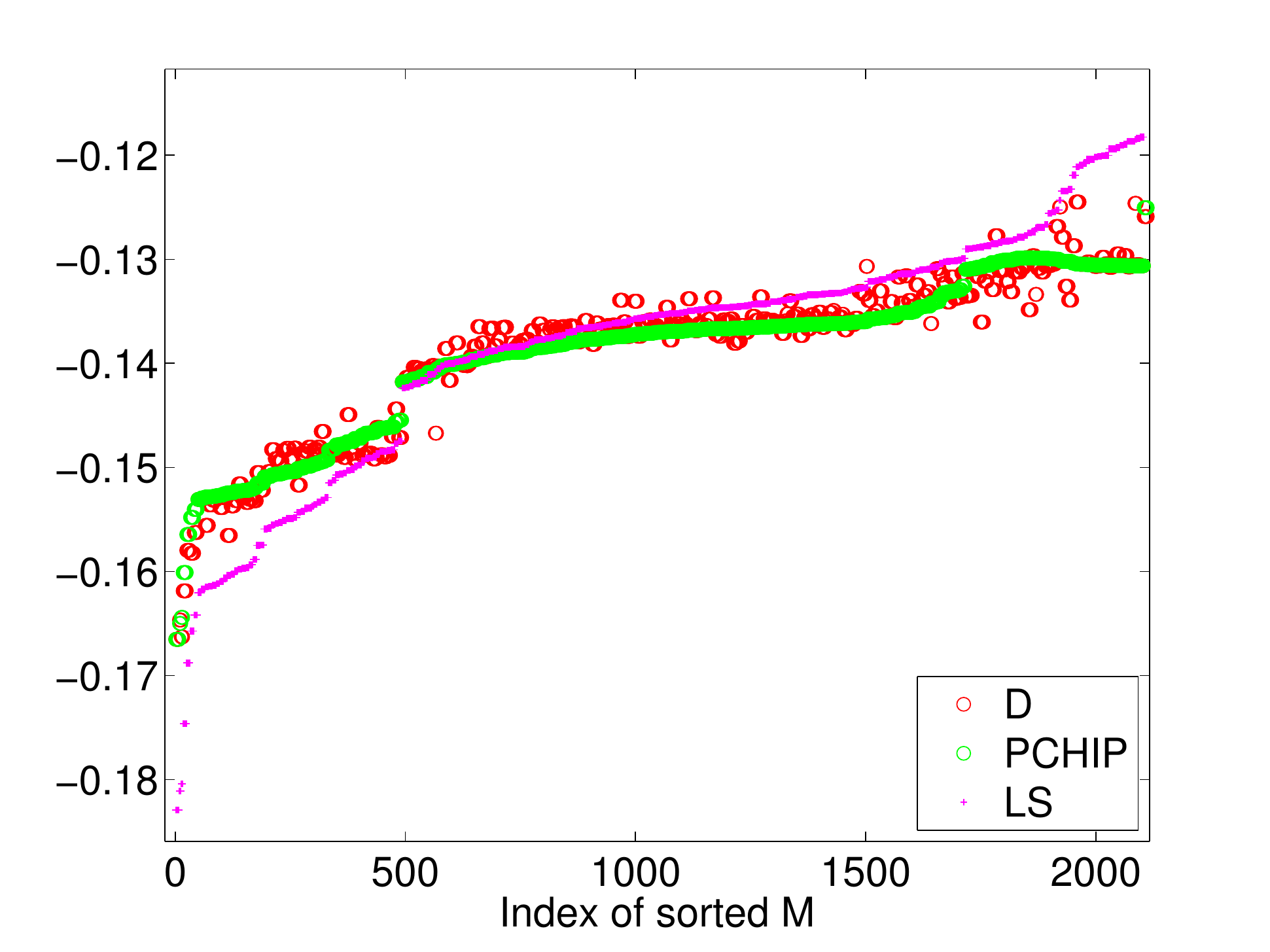}\label{fig:Results in diagonal index range [1,2000]}}                
  \subfigure[Diagonal elements in (2000,3000)]{\includegraphics[width=0.32\textwidth]{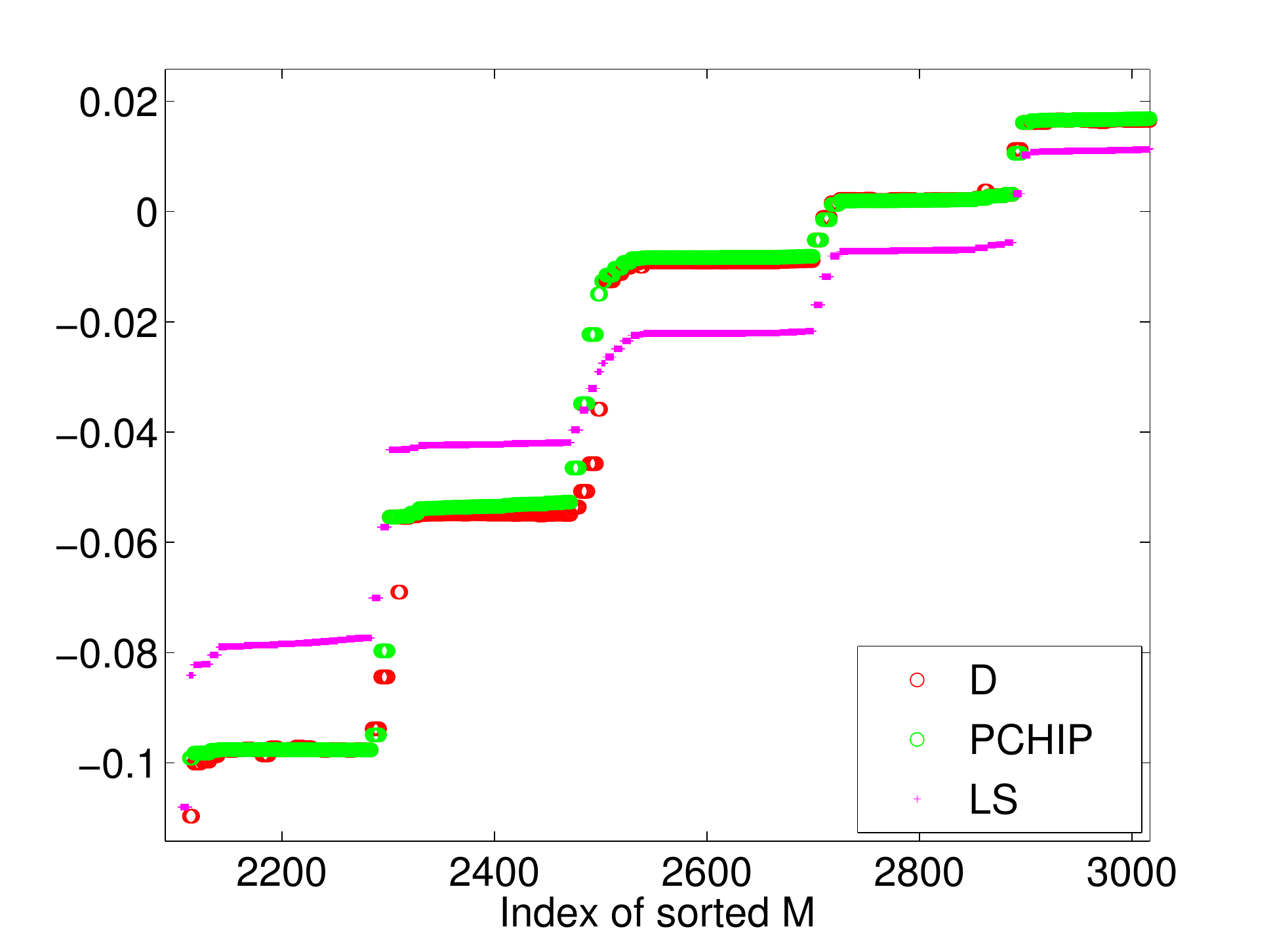}\label{fig:Results in diagonal index range [2000,3000]}}   
  \subfigure[Diagonal elements in (3000,5000)]{\includegraphics[width=0.32\textwidth]{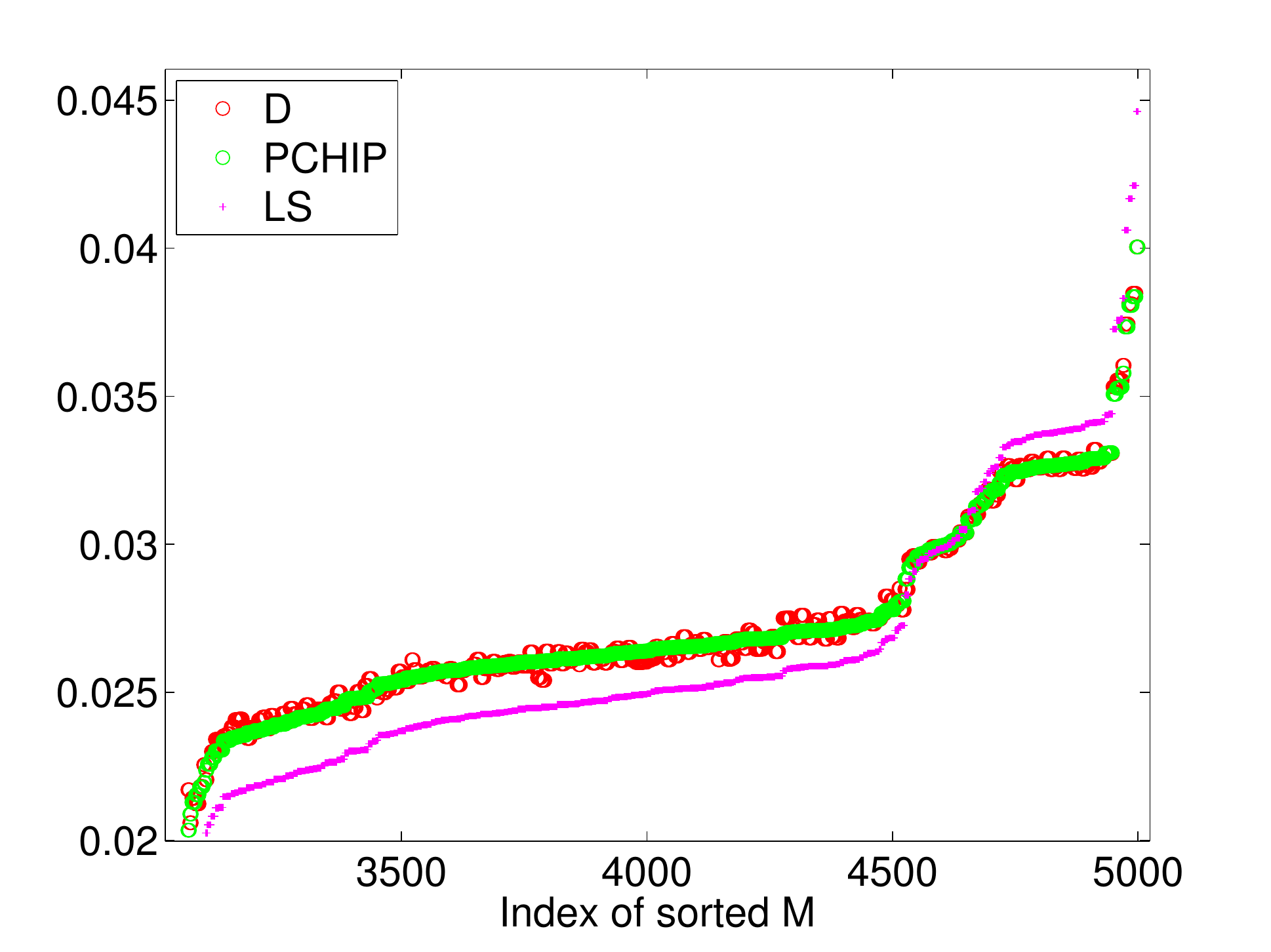}\label{fig:Results in diagonal index range [3000,5000]}}   
  \caption{Comparing the linear LS model with the PCHIP model on RDB5000 matrix with $M$ from SVD.}
  \label{fig:RDB5000-Comparision of two fitting models}
\end{figure}

In Figure \ref{fig:RDB5000-Comparision of two fitting models}, we divide the diagonal elements of the matrix RDB5000 into three contiguous sets and zoom in the details. We see that despite a
 good $M$, the linear LS model cannot scale the entire $M$ onto $D$. The more flexible
 piecewise approach of PCHIP results in a much better fit.

\begin{figure}
  \centering  
  \subfigure[Matrix Heatflow160]
  {\includegraphics[width=0.32\textwidth]{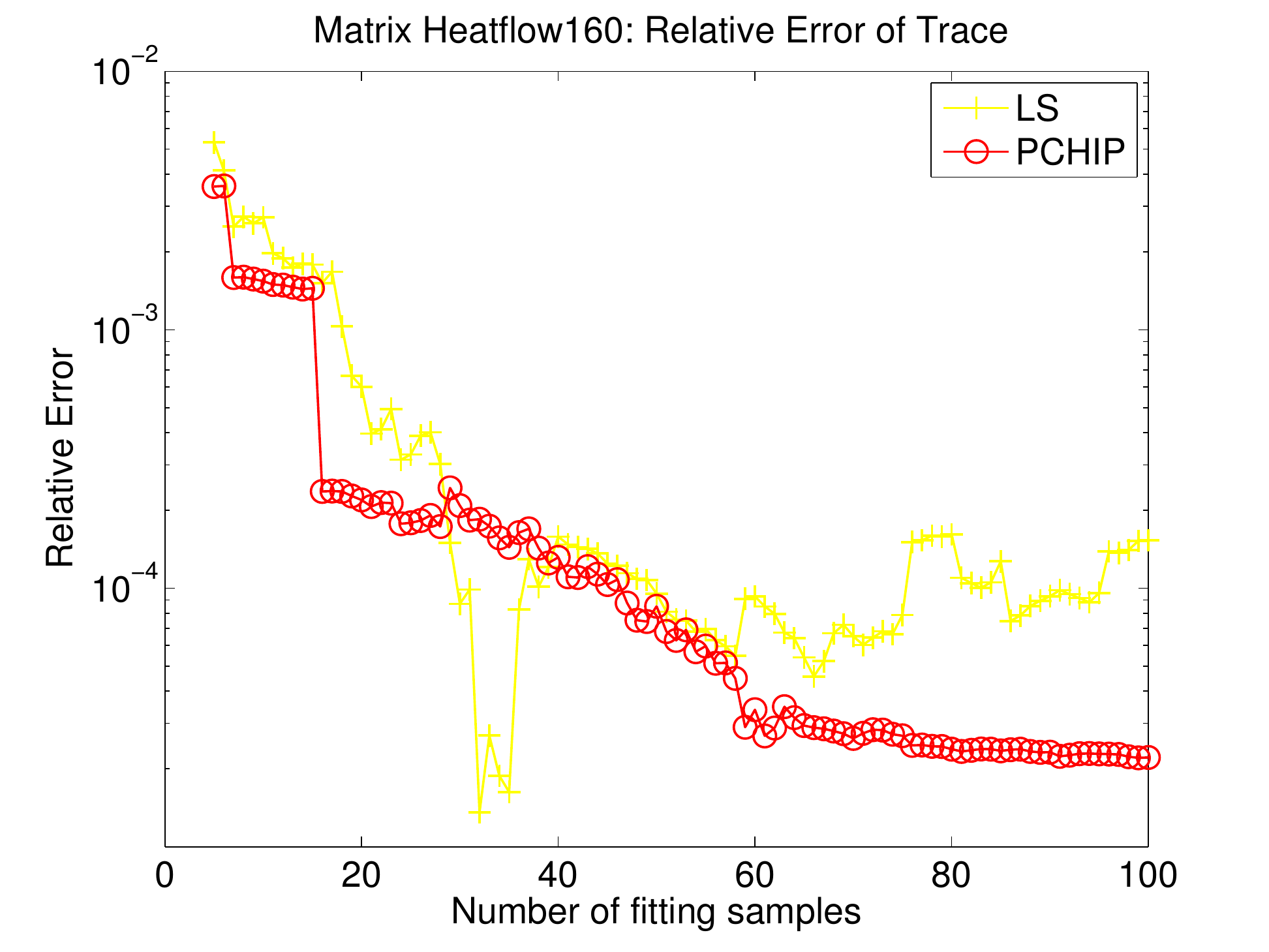}\label{fig:trace estimation on Matrix Heatflow160}} 
  \subfigure[Matrix Poisson150]
  {\includegraphics[width=0.32\textwidth]{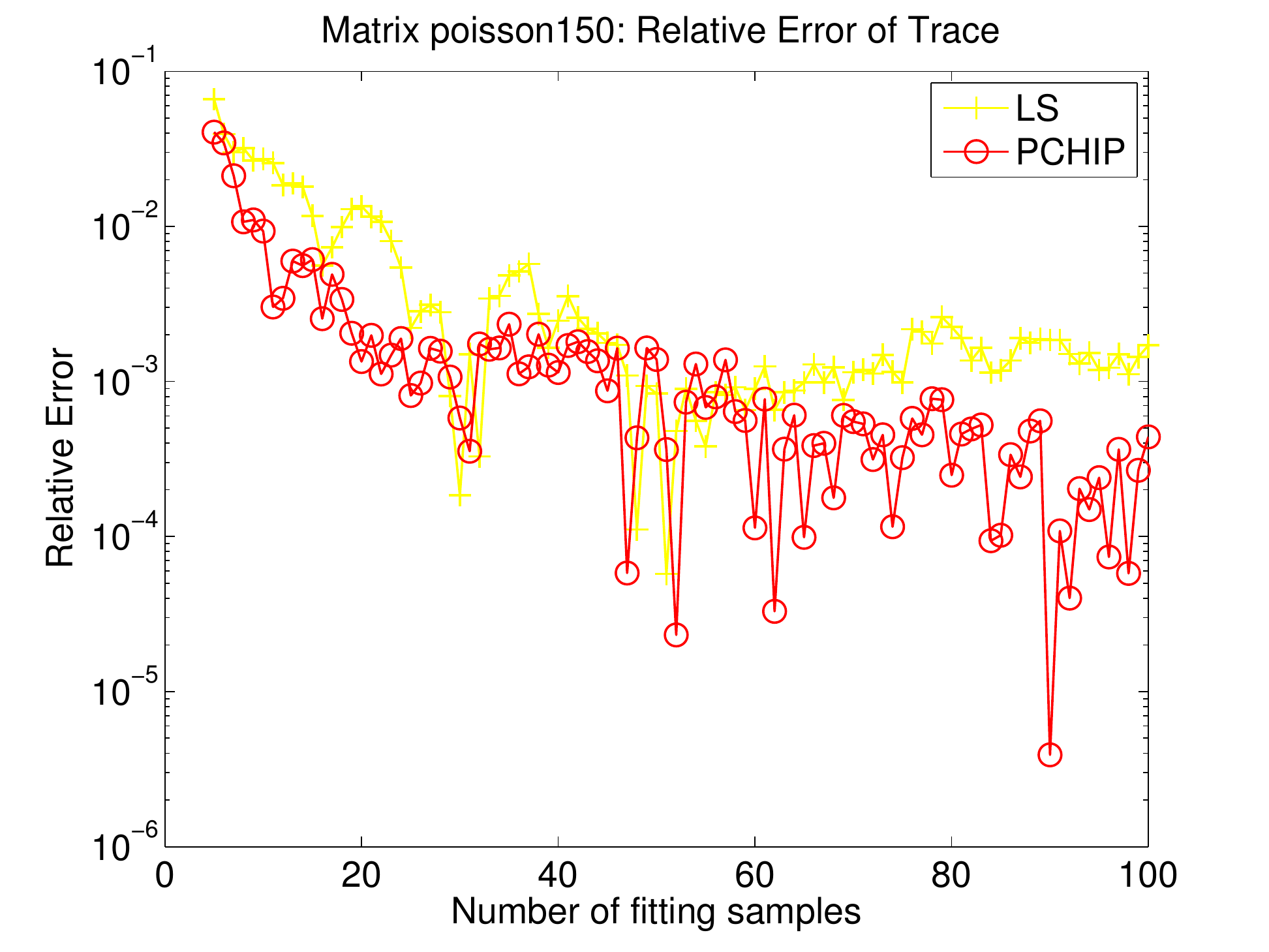}\label{fig:trace estimation on Matrix poisson150}} 
  \subfigure[Matrix VFH6]
  {\includegraphics[width=0.32\textwidth]{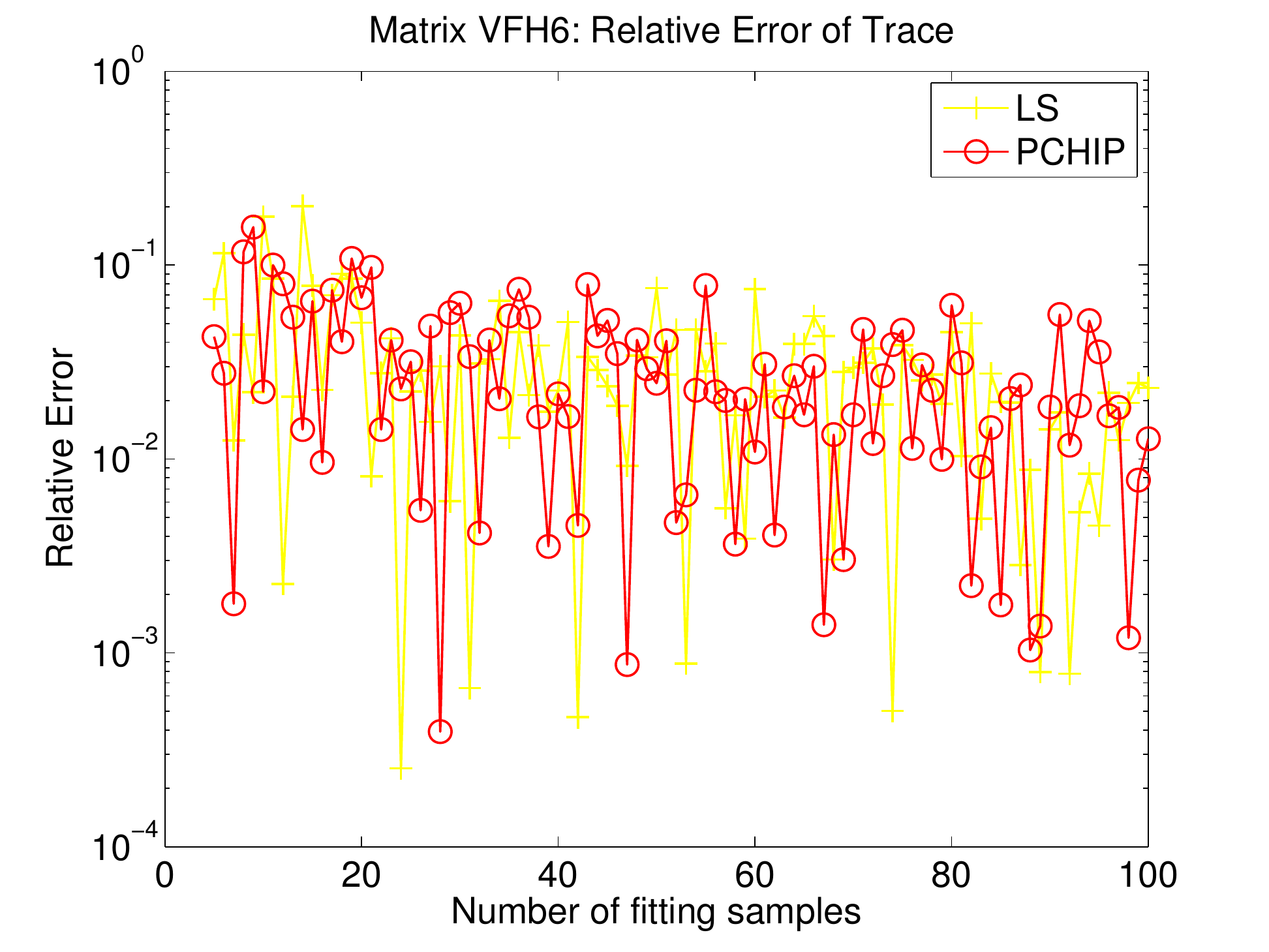}\label{fig:trace estimation on Matrix VFH6}} 
  \caption{Comparing relative trace error between the LS model and the PCHIP model in three typical cases with SVD.}
  \label{fig:Relative trace error in various cases}
\end{figure}

\begin{figure} 
  \centering
  \subfigure[Matrix Heatflow160]
  {\includegraphics[width=0.32\textwidth]{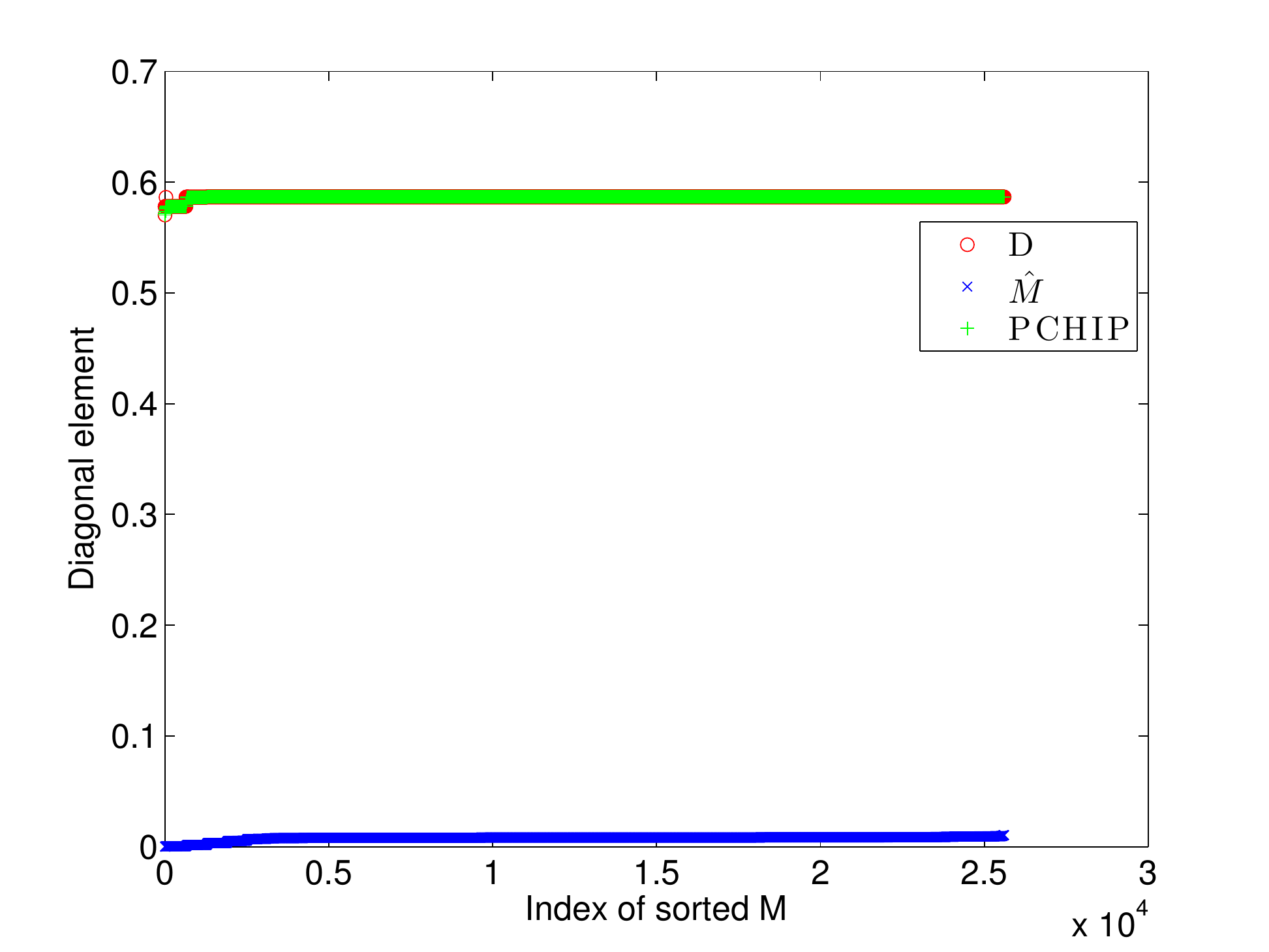}\label{fig:PCHIP on Matrix Heatflow160}} 
  \subfigure[Matrix Poisson150]
  {\includegraphics[width=0.32\textwidth]{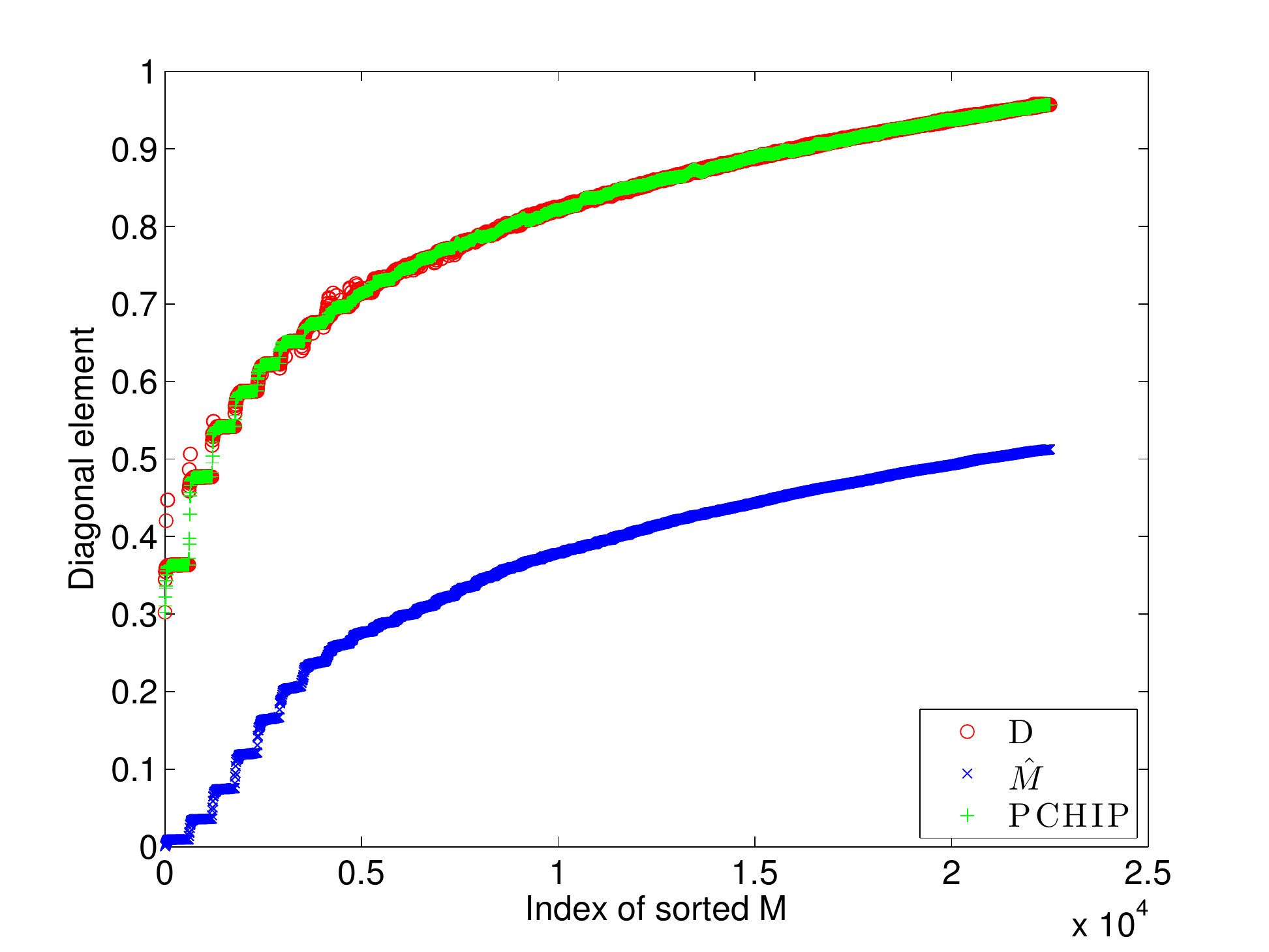}\label{fig:PCHIP on Matrix poisson150}}                
  \subfigure[Matrix VFH6]
  {\includegraphics[width=0.32\textwidth]{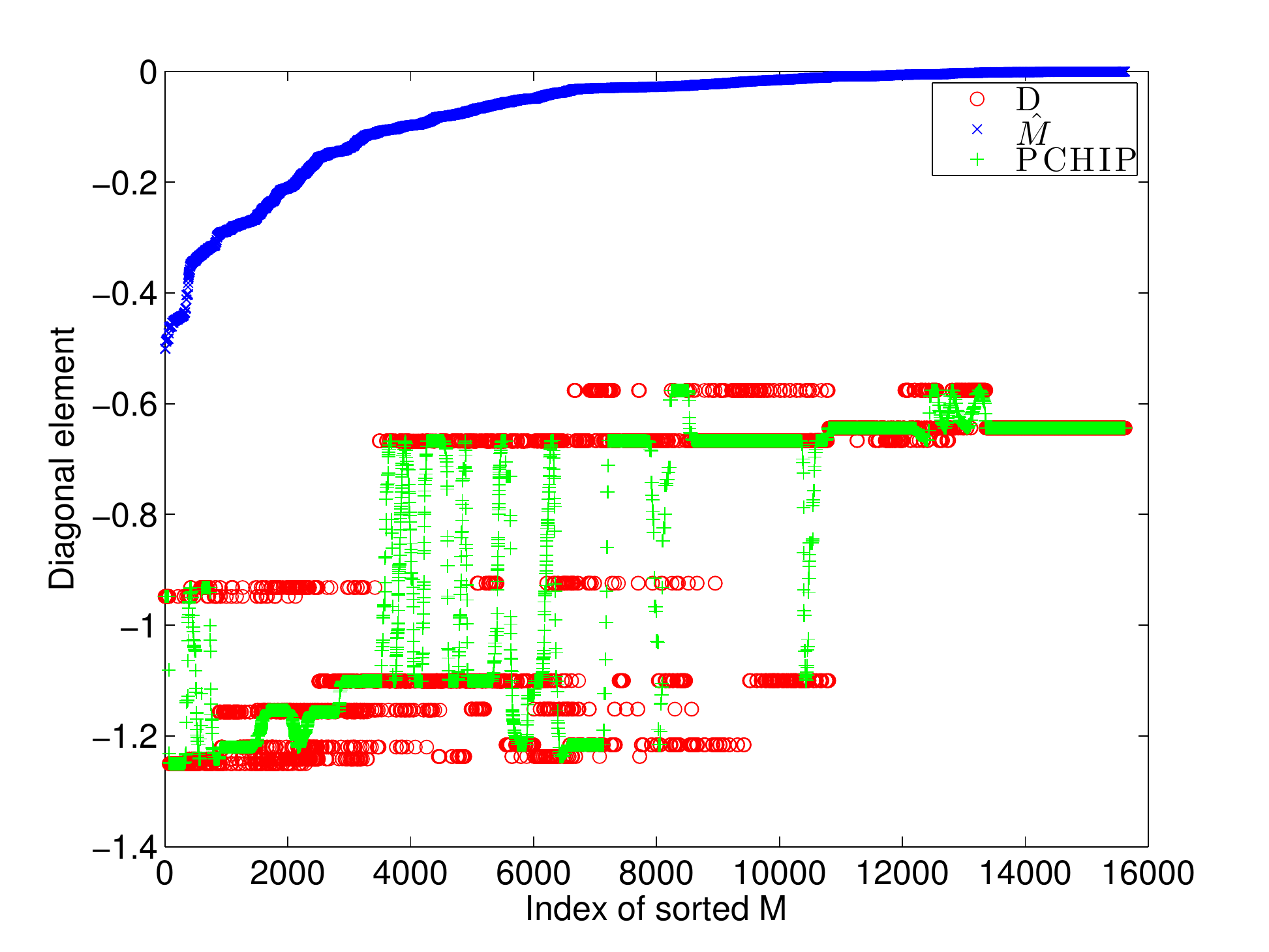}\label{fig:PCHIP on Matrix VFH6}}   
  \caption{Fitting results of three typical cases with the PCHIP model and 100 fitting points using SVD.}
  \label{fig:Different fitting results of three typical cases with the PCHIP model}
\end{figure}

In Figure \ref{fig:Relative trace error in various cases}, we look at three matrices 
  with $M$ generated using the SVD.
The PCHIP model typically has smaller relative trace error than the LS model. 
We also see that as more fitting points are sampled, the relative trace error of both models 
  decreases significantly at early stages and slowly after a certain point.
This relates to the quality of $M$, not of the model.
Typically $M$ will approximate the global pattern of $D$ and the two can be matched well
  with only a few fitting points.
But if the local patterns of $M$ and $D$ differ, a large number of fitting points will
  be required.
This can be seen in Figure \ref{fig:Different fitting results of three typical cases with the PCHIP model}.
Using the permutation of $\hat M$ for each case, we plot $\hat M$, $D$,
 and $p(\hat M)$ using 100 fitting points.
For the Heatflow160 matrix, $p(M)$ approximates $D$ well everywhere
  except for the small leftmost part of the plot, which allows 
  the relative error to reach below $10^{-4}$
  before convergence slows down (Figure \ref{fig:Relative trace error in various cases}).
The behavior is similar for the Poisson150.
The issue is more pronounced on matrix VFH6, where $M$ and $p(M)$ capture the average location 
  of $D$ but completely miss the local pattern, which is reflected by a very slowly 
  improving error in Figure \ref{fig:Relative trace error in various cases}.

We mention that the irregularity of the relative trace errors in 
  Figure \ref{fig:Relative trace error in various cases} relates to the
  variability of successive updates of $M$ and of the sampling indices,
  especially when $M$ is of lower quality.

Figure \ref{fig:comparison of actual variance in various cases} demonstrates that the PCHIP model has smaller actual variance for MC on $E_{fit}$ than the LS model.
Therefore, we only consider the PCHIP model in the rest of experiments.

%---------------------------------------------------------------
\subsection{Comparison between the fitting model and different MC methods}
We address the question of whether the number of matrix inversions we spend
  on computing the fitting could have been used more efficiently in an MC method, 
  specifically the Hutchinson method on $A^{-1}$ and the Hutchinson method 
  on $E$.
In Table \ref{ta:Comparison among the fitting model and different MC methods} 
  we compare the relative trace error of the PCHIP model with 20 fitting points 
  against the relative errors of the two MC methods as computed explicitly from
  their respective standard deviations in  (\ref{eq: Hutchinson variance of invA}) 
  and (\ref{eq: Hutchinson variance of E}), with $s=20$, divided by the actual 
  trace of $D$.
%(\ref{eq:unit vector variance of E_fit}).

When $M$ approximates $D$ sufficiently well, the trace from the fitted
  diagonal is better, and for the ILU approximations far better, than if we just 
  use the Hutchinson method on $A^{-1}$ (the first column of results). 
Although MC on $E$ exploits the ILU or SVD approximation of the entire matrix 
  (not just the diagonal that our method uses), we see that it does not always 
  improve on MC on $A^{-1}$, and in some cases (cfd1 with ILU) it is far worse.
In contrast, our diagonal fitting typically improves on MC on $E$.
The last column shows that even with an inexpensive diagonal approximation
  we obtain a similar or better error than MC on $A^{-1}$.
The only exception is the matrix VFH6 where, as we saw earlier, $M$ cannot 
  capture the pattern of $D$. 
Even then, its error is close to the errors from the MC methods and, 
  as we show next, the best method can be identified dynamically with only 
  a small number of samples.

%\begin{table}[htbp]
%\caption{Comparing the accuracy of the fitting model and various MC methods in (\ref{eq: Hutchinson variance of invA}), (\ref{eq: Hutchinson variance of E}) and (\ref{eq:unit vector variance of E_fit}) with only 20 fitting points.}
%\centering

%\begin{center}
%    \begin{tabular}{ |c|c|cc|c|cc|cc|}
%    \hline
%    	\multicolumn{1}{|r}{}
%	 & \multicolumn{1}{c}{{\tt Trace}} 
%	 & \multicolumn{2}{c}{{\tt RelErr}} 
%	 & \multicolumn{1}{c}{{\tt $Std(T_{Z_2}(A^{-1}))$}}
%	 & \multicolumn{2}{c}{{\tt $Std(T_{Z_2}(E))$}}
%	 & \multicolumn{2}{c|}{{\tt $Std(T_{e_i}(E_{fit}))$}} \\ \hline
%Matrix 	        & --  & ILU   & SVD 	&ILU or SVD	 &ILU	 & SVD 	&ILU 	& SVD    \\ \hline
%    RDB5000 	& 2.7E2 &8.1E-3	&9.7E-3 & 1.4E1  &1.3E1  &3.2E0	&5.6E0 	& 3.8E0  \\ \hline
%    Wathen100 	& 1.4E3 &1.7E-4 &2.0E-2 & 1.4E0  &1.7E-1 &1.0E0	&3.4E0  & 3.0E2   \\ \hline
%    Heatflow160	& 1.5E4 &1.6E-7 &2.2E-4 & 7.3E0  &6.0E-1 &7.3E0	&1.3E-2 & 4.9E0   \\ \hline
%    Poisson150 	& 1.8E4 &2.3E-3 &1.3E-3 & 4.7E2  &4.6E2  &7.7E1	&2.0E2  & 6.3E1   \\ \hline
%    VFH6       	& 1.3E4 &6.8E-5 &6.8E-2 & 4.2E1  &4.0E-1 &4.1E1	&7.4E-1 & 9.3E2 \\ \hline
%    \end{tabular}
%\end{center}
%\label{ta:Comparison among the fitting model and different MC methods}
%\end{table} 

\begin{table}[htbp]
\caption{Relative trace error from our PCHIP model and from the MC method 
  on $A^{-1}$ and on $E$ (computed explicitly as the standard deviation with $s=20$ from
  (\ref{eq: Hutchinson variance of invA}) and (\ref{eq: Hutchinson variance of E}) 
  divided by the actual trace).}
%and (\ref{eq:unit vector variance of E_fit}) with only 20 fitting points.}
\centering
\small
\begin{center}
    \begin{tabular}{ |c|c|cc|cc|c|}
    \hline
    	\multicolumn{1}{|r}{} &
	 & \multicolumn{2}{c|}{{\tt ILU}} 
	 & \multicolumn{2}{c|}{{\tt SVD}} 
	 & \multicolumn{1}{c|}{{\tt Bounds}} \\ \hline
    Matrix 	& $T_{Z_2}(A^{-1})$  &PCHIP  & $T_{Z_2}(E)$ & PCHIP  & $T_{Z_2}(E)$  &	PCHIP  \\ \hline

    RDB5000 	&5.2E-2 & 8.1E-3 &4.8E-2  &4.1E-3 &1.2E-2 & 5.3E-2 \\ \hline
    cfd1 	&1.3E-1 & 2.8E-2 & 8.2E+2  &8.8E-3 &1.8E-2 & 2.6E-2 \\ \hline
    Heatflow160	&4.9E-4 & 1.6E-7 &4.0E-5  &2.0E-4 &4.9E-4 &3.5E-4 \\ \hline
    Poisson150 	&2.6E-2 & 2.3E-3 &2.5E-2  &1.4E-3 &4.3E-3 & 8.3E-3 \\ \hline
    VFH6       	&3.2E-3 & 6.8E-5 &3.1E-5  &1.0E-2 &3.2E-3 & 6.0E-2 \\ \hline

    \end{tabular}
\end{center}
\label{ta:Comparison among the fitting model and different MC methods}
\end{table}

\begin{figure}[htbp]
  \centering
  \subfigure[Matrix Heatflow160]
  {\includegraphics[width=0.32\textwidth]{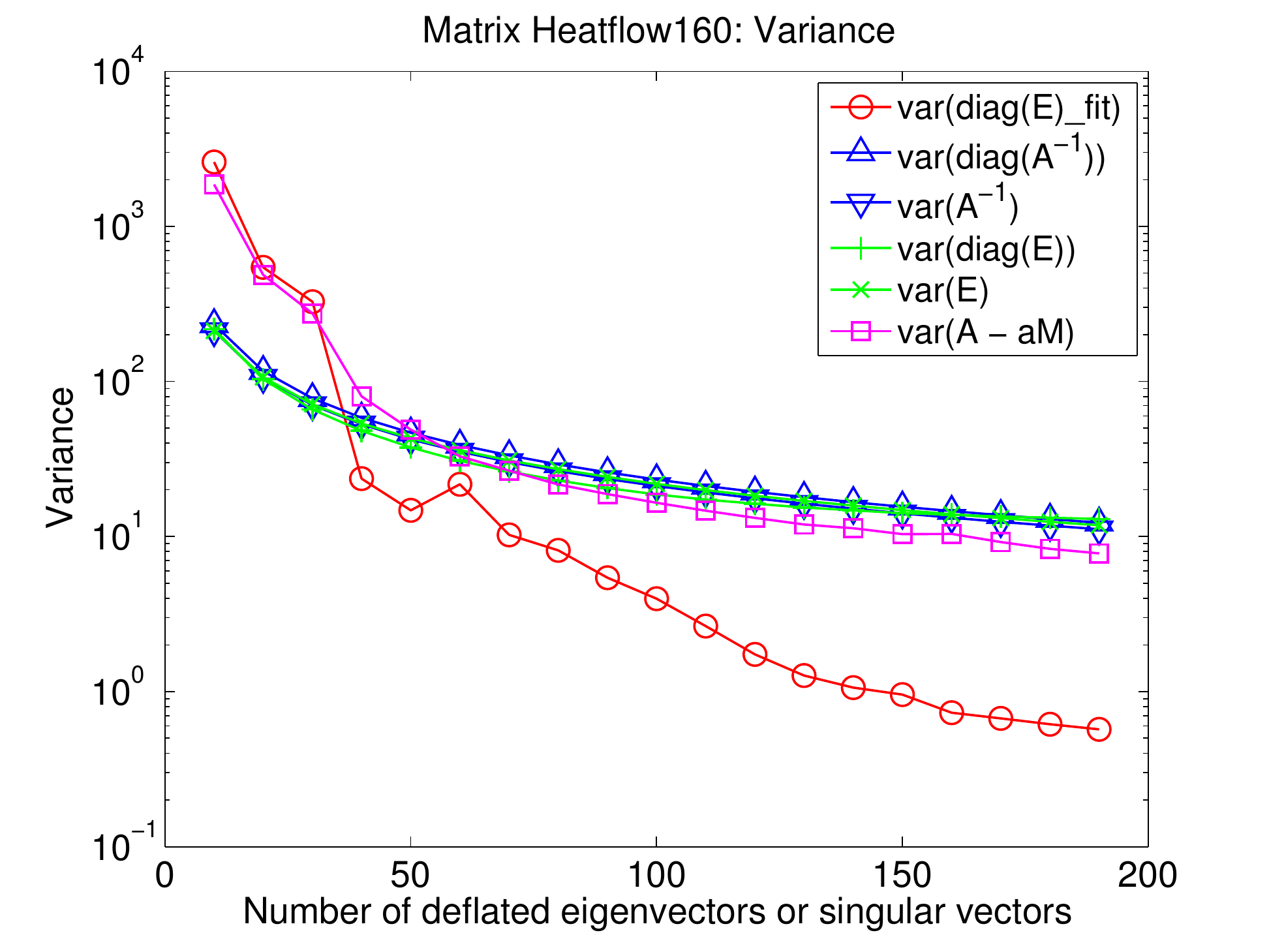}\label{fig:comparison of actual variance on Matrix Heatflow160}}  
  \subfigure[Matrix Poisson150]
  {\includegraphics[width=0.32\textwidth]{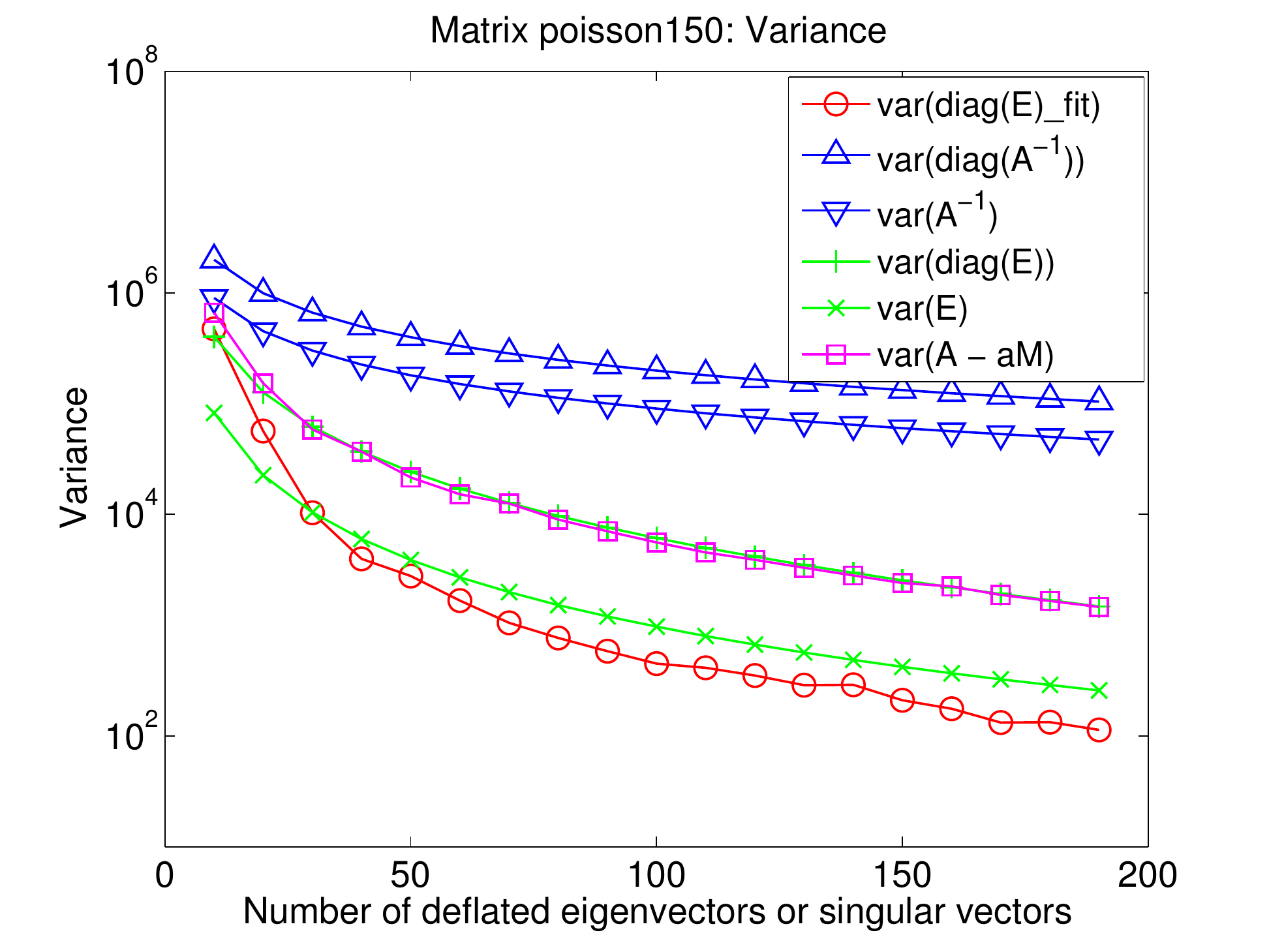}\label{fig:comparison of actual variance on Matrix poisson150}}                
  \subfigure[Matrix VFH6]
  {\includegraphics[width=0.32\textwidth]{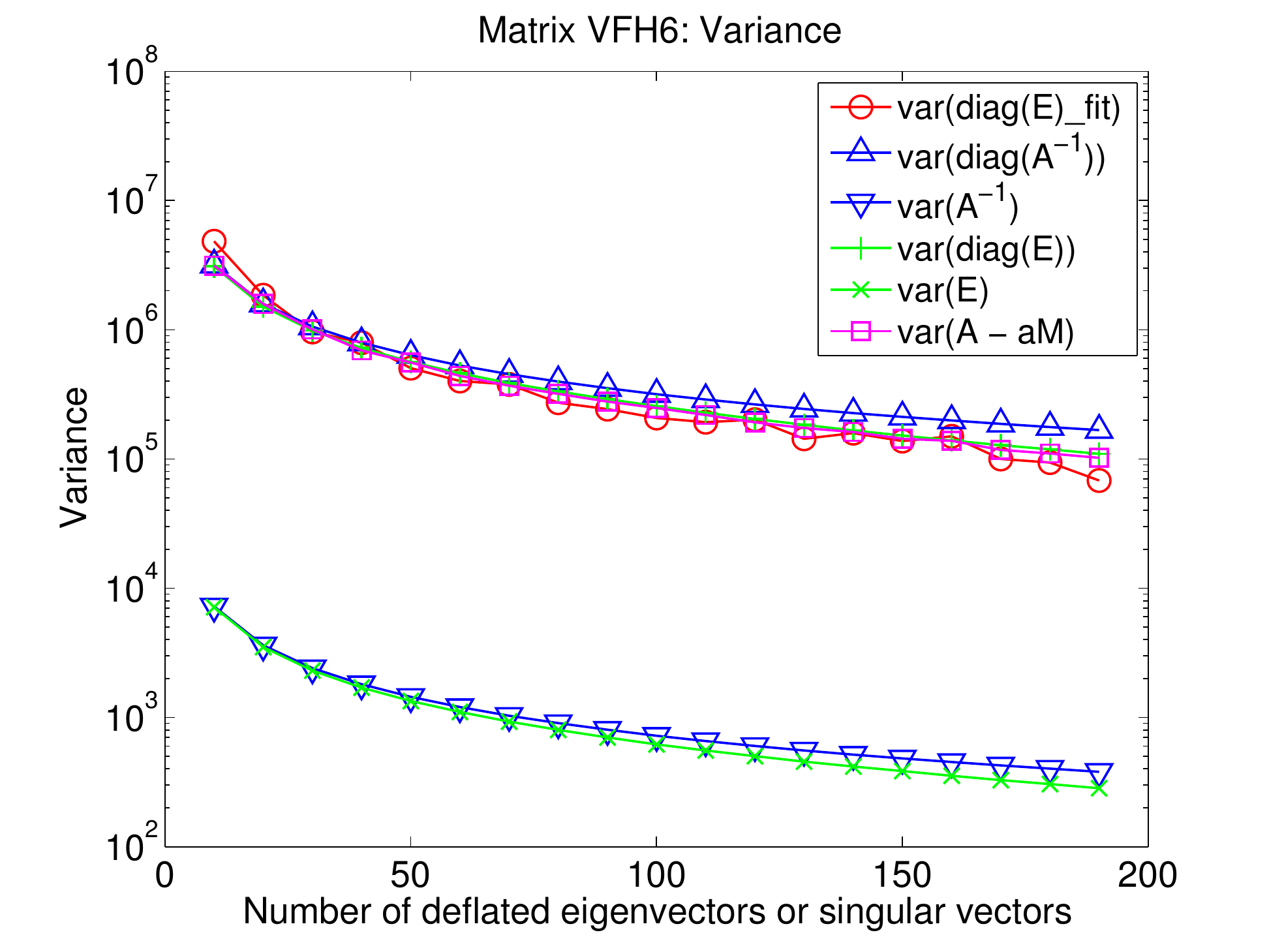}\label{fig:comparison of actual variance on Matrix VFH6}}   
  \caption{Comparing actual variances of different MC methods in three typical cases with SVD.}
  \label{fig:comparison of actual variance in various cases}
\end{figure}

If the user requires better trace accuracy than our fitting technique provides,
  we explore the performance of the diagonal fitting as a variance reduction for 
  MC with unit vectors.
We compute the actual values of $Var(T_{e_i}(E_{fit}))$, $Var(T_{e_i}(A^{-1}))$, $Var(T_{e_i}(E))$, $Var(T_{Z_2}(A^{-1}))$ and $Var(T_{Z_2}(E))$ for every step $s=5:100$ 
and show results for three matrices in 
Figure \ref{fig:comparison of actual variance in various cases}.
Note that the low rank approximation uses $2s$ singular vectors.
As before, for Heatflow160, MC on $E_{fit}$ performs much better than other MC methods, achieving about two orders reduction in variance. For Poisson150, MC on $E_{fit}$ is slightly better compared to the Hutchinson method on $E$. In contrast, for VFH6, MC with unit vectors do not perform well regardless of the diagonal.

\subsection{Dynamic evaluation of variance and relative trace error}
The above results emphasize the importance of being able to assess quickly 
  and accurately the relative differences between the variances of different 
  methods as well as the trace error, so that we can decide whether to continue 
  with fitting or which MC method to switch to.
First we show the effectiveness of the dynamic variance evaluation algorithm 
  for our fitting MC method on $(D - p(M))$ with unit vectors, 
  and on $A^{-1}$ and $E$ with Rademacher vectors. 
Then, we evaluate our algorithm for estimating the relative trace error during 
  the fitting process. 

\begin{figure}
  \centering
  \subfigure[Matrix Heatflow160]
  {\includegraphics[width=0.32\textwidth]{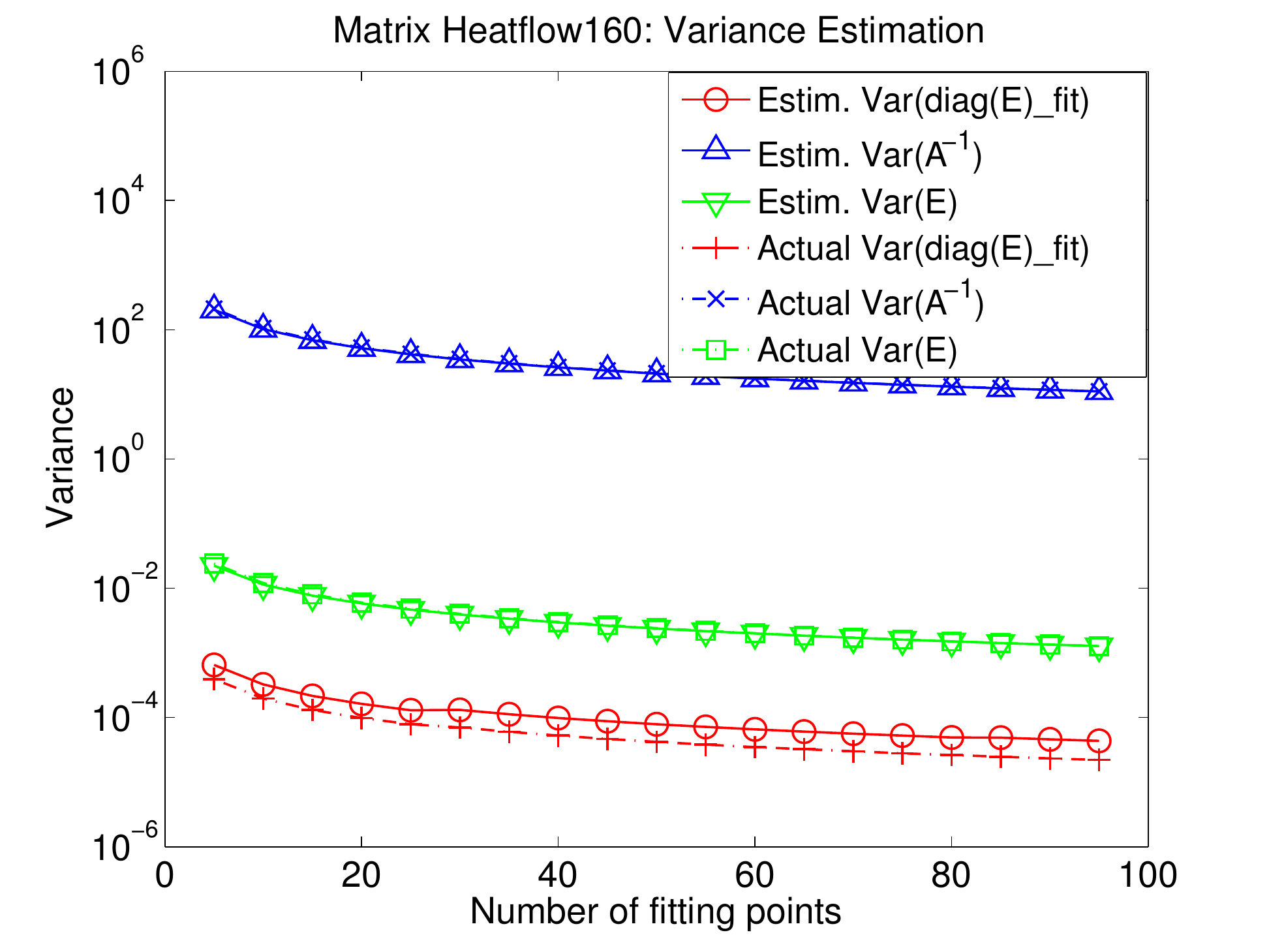}\label{fig:comparison of estimated variances and actual variances on Matrix Heatflow160ILU}} 
  \subfigure[Matrix Poisson150]
  {\includegraphics[width=0.32\textwidth]{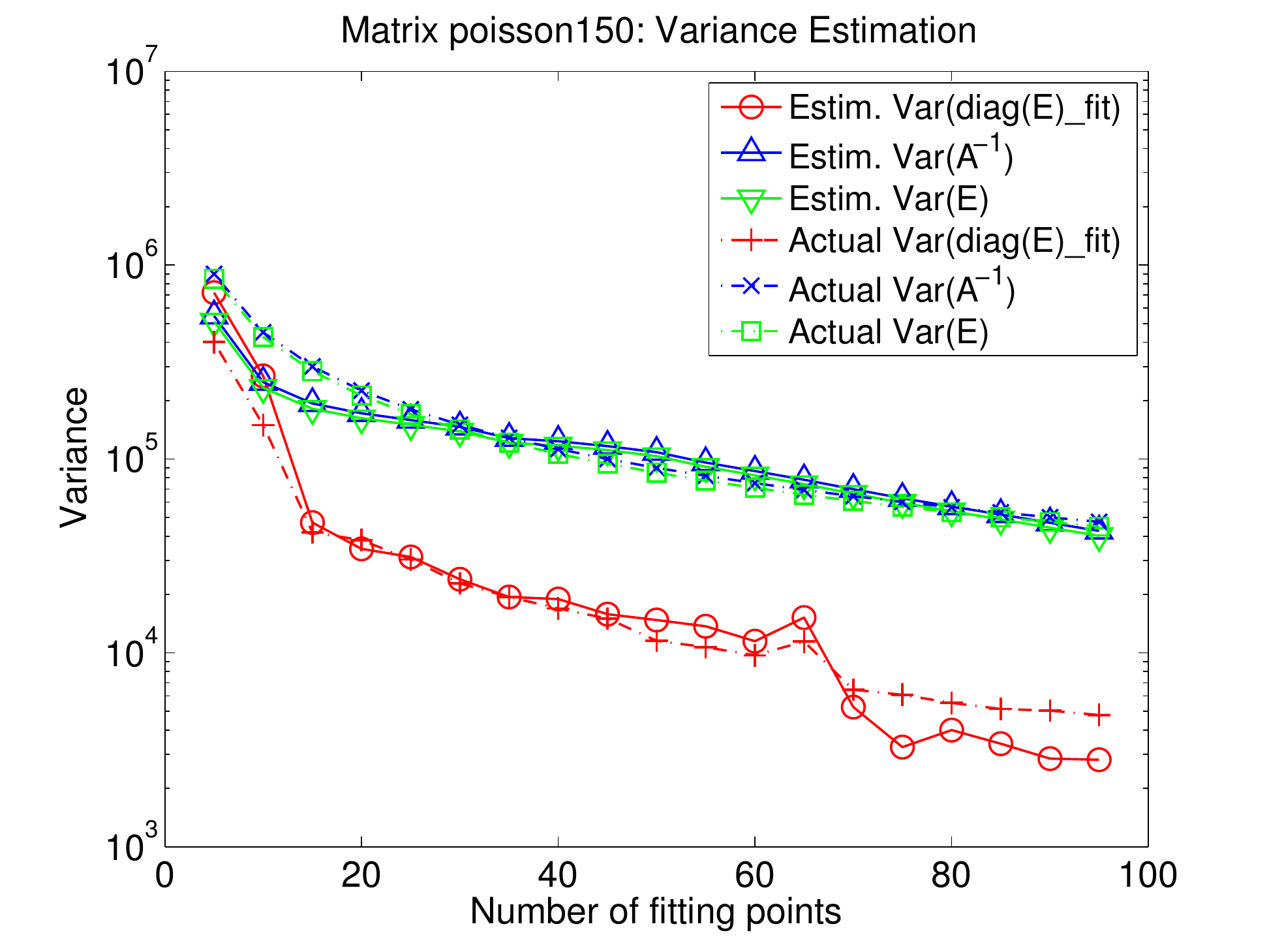}\label{fig:comparison of estimated variances and actual variances on Matrix poisson150ILU}}                
  \subfigure[Matrix VFH6]
  {\includegraphics[width=0.32\textwidth]{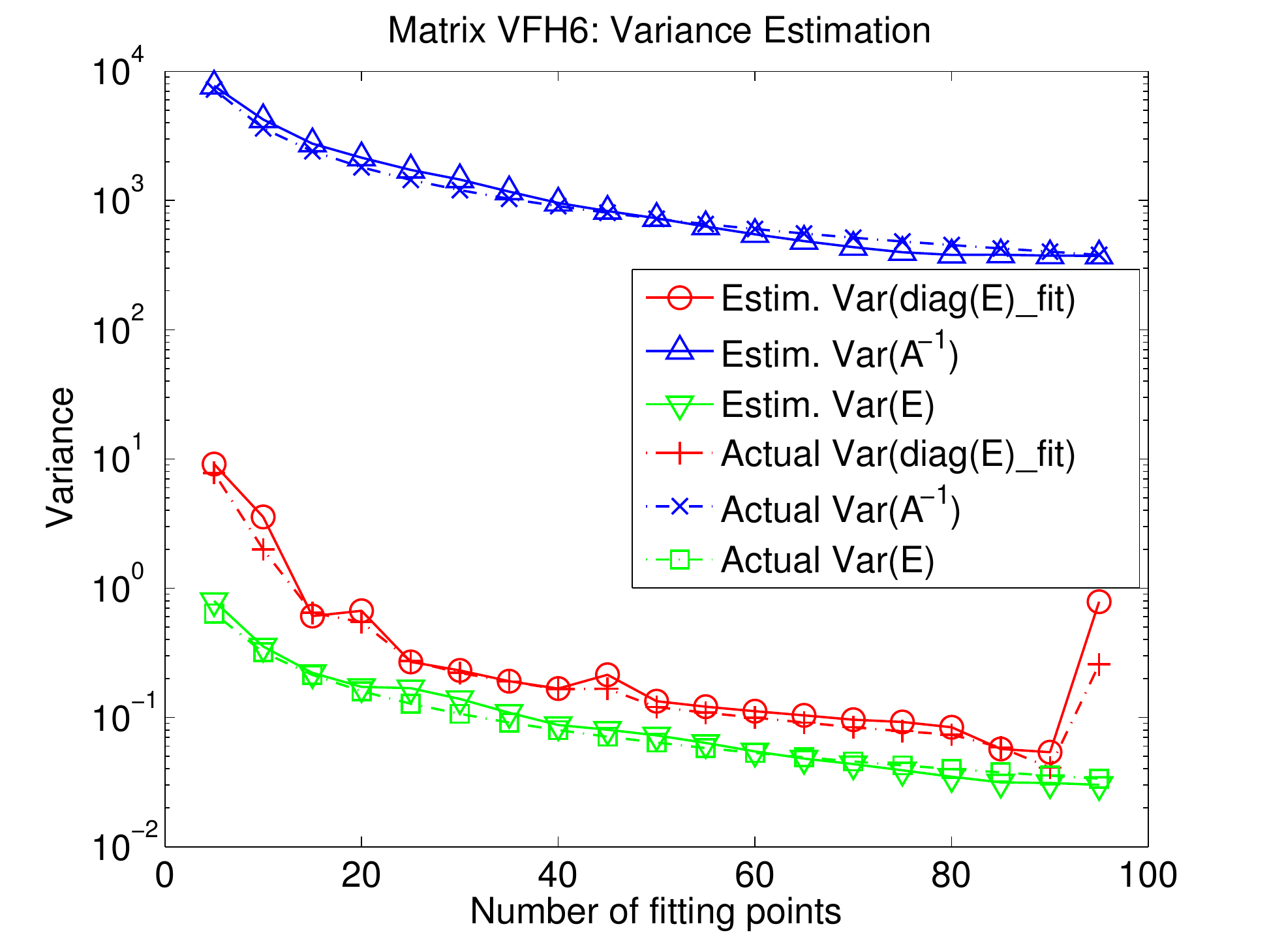}\label{fig:comparison of estimated variances and actual variances on Matrix VFH6ILU}}   
  \caption{Comparing estimated variances and actual variances of three MC methods with ILU.}
  \label{fig:comparison of the estimated variances with actual variances of MC in various cases with ILU}
\end{figure}

Figure \ref{fig:comparison of the estimated variances with actual variances of MC in various cases with ILU} compares the estimated variances with the actual variances of the three MC methods when increasing the number of fitting points from 5 to 100. The approximation $M$ is computed by using ILU.
We can see that the estimated values of $Var(T_{Z_2}(A^{-1}))$ and 
$Var(T_{Z_2}(E))$ converge to the actual variances after only a few sample points. 
The estimated value of $Var(T_{e_i}(E_{fit}))$ gets close to and captures the trend of the the actual variance as the fitting samples increase. Nevertheless, the relative differences between the variances of the various MC methods are apparent almost immediately. 

\begin{figure}
  \centering
  \subfigure[Matrix Heatflow160]
  {\includegraphics[width=0.32\textwidth]{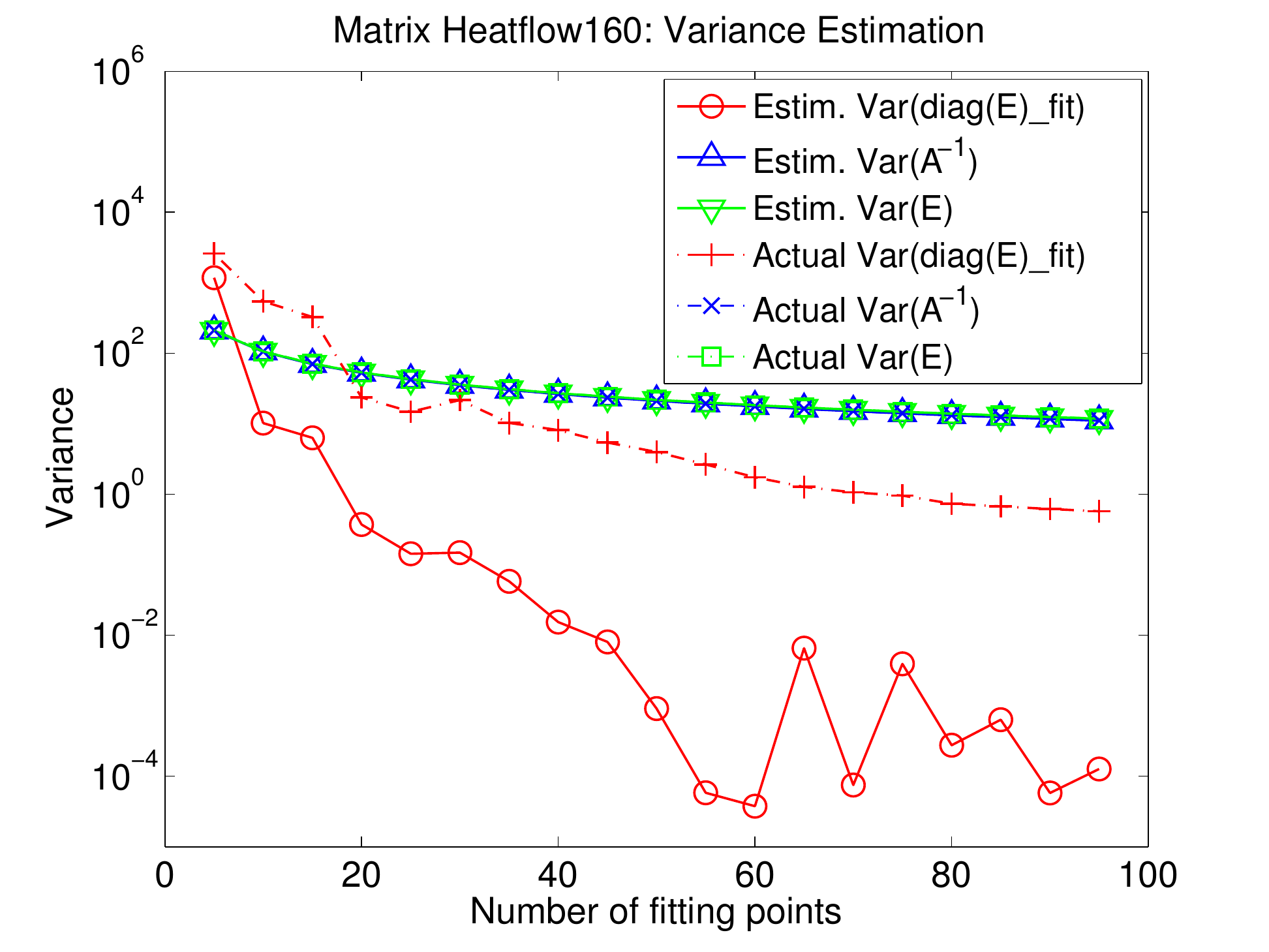}\label{fig:comparison of estimated variances and actual variances on Matrix Heatflow160SVD}} 
  \subfigure[Matrix Poisson150]
  {\includegraphics[width=0.32\textwidth]{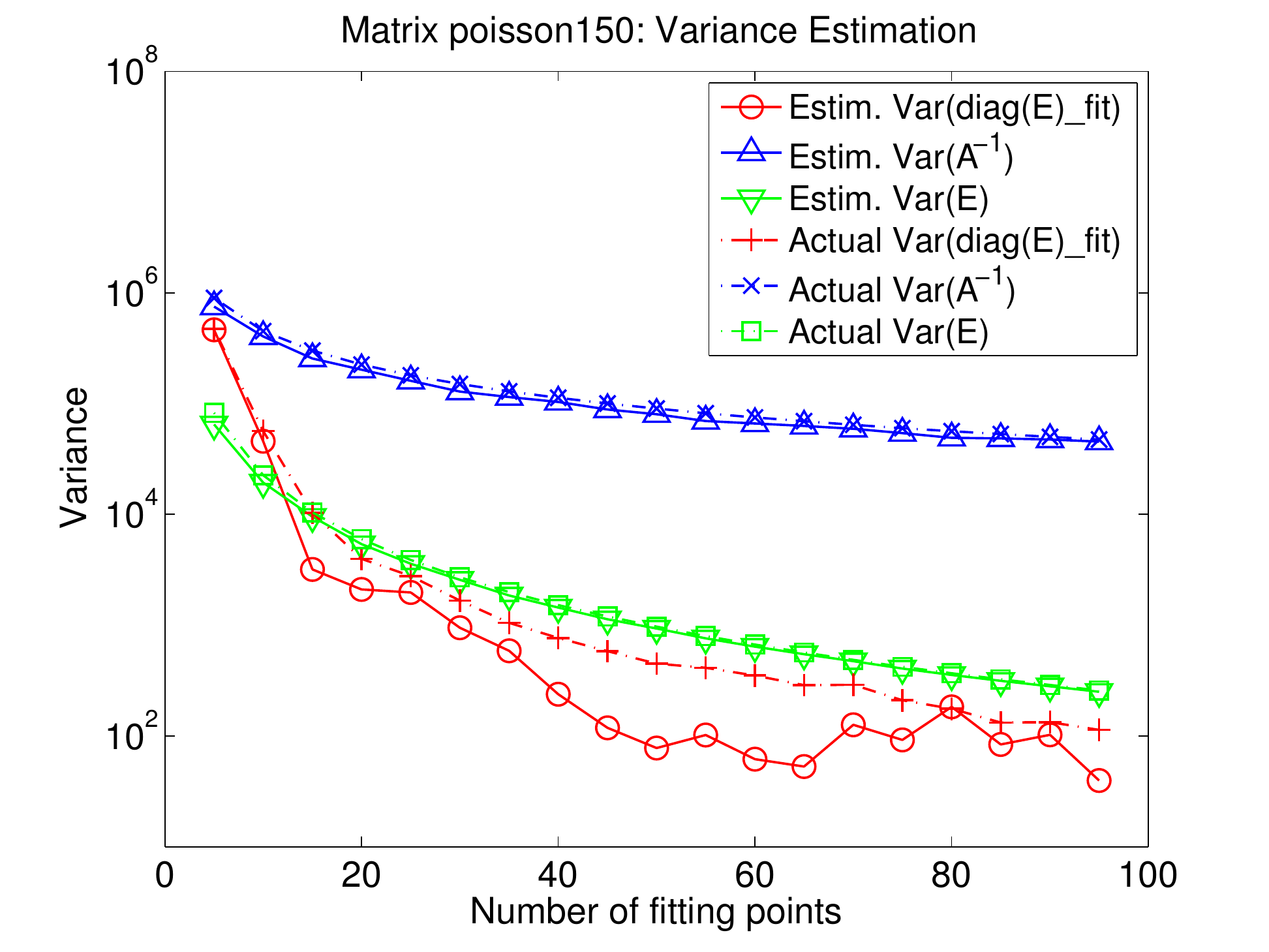}\label{fig:comparison of estimated variances and actual variances on Matrix poisson150SVD}}                
  \subfigure[Matrix VFH6]
  {\includegraphics[width=0.32\textwidth]{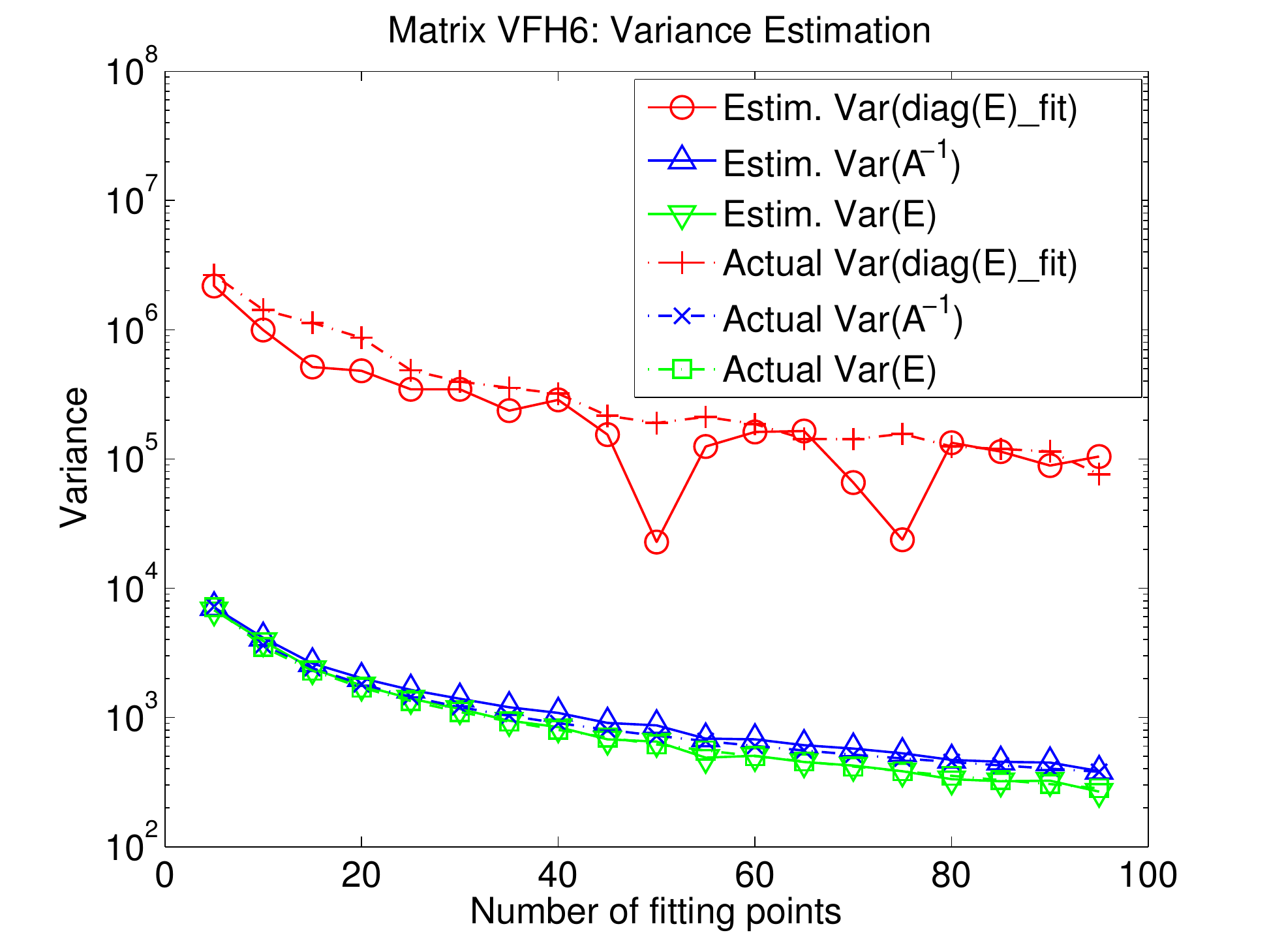}\label{fig:comparison of estimated variances and actual variances on Matrix VFH6SVD}}   
  \caption{Comparing estimated variances and actual variances of three MC methods with SVD.}
  \label{fig:comparison of the estimated variances with actual variances of MC in various cases with SVD}
\end{figure}

Figure \ref{fig:comparison of the estimated variances with actual variances of MC in various cases with SVD} shows the same experiments when the approximation $M$ is computed by SVD. 
Since $M$ is updated each step, $Var(T_{Z_2}(E))$ and $Var(T_{e_i}(E_{fit}))$ 
change accordingly.
As with ILU, $Var(T_{Z_2}(A^{-1}))$ and $Var(T_{Z_2}(E)$ can be estimated very well 
in a few steps. 
$Var(T_{e_i}(E_{fit}))$ could be underestimated but the relative variance difference between these MC methods becomes clear when the fitting points increase beyond 20. 
Thus we are able to determine whether the fitting process is beneficial as 
  a variance reduction preprocessing and which is the best MC method to proceed 
  with for the trace estimation. 

\begin{figure}
  \centering
  \subfigure[Matrix Heatflow160]
  {\includegraphics[width=0.32\textwidth]{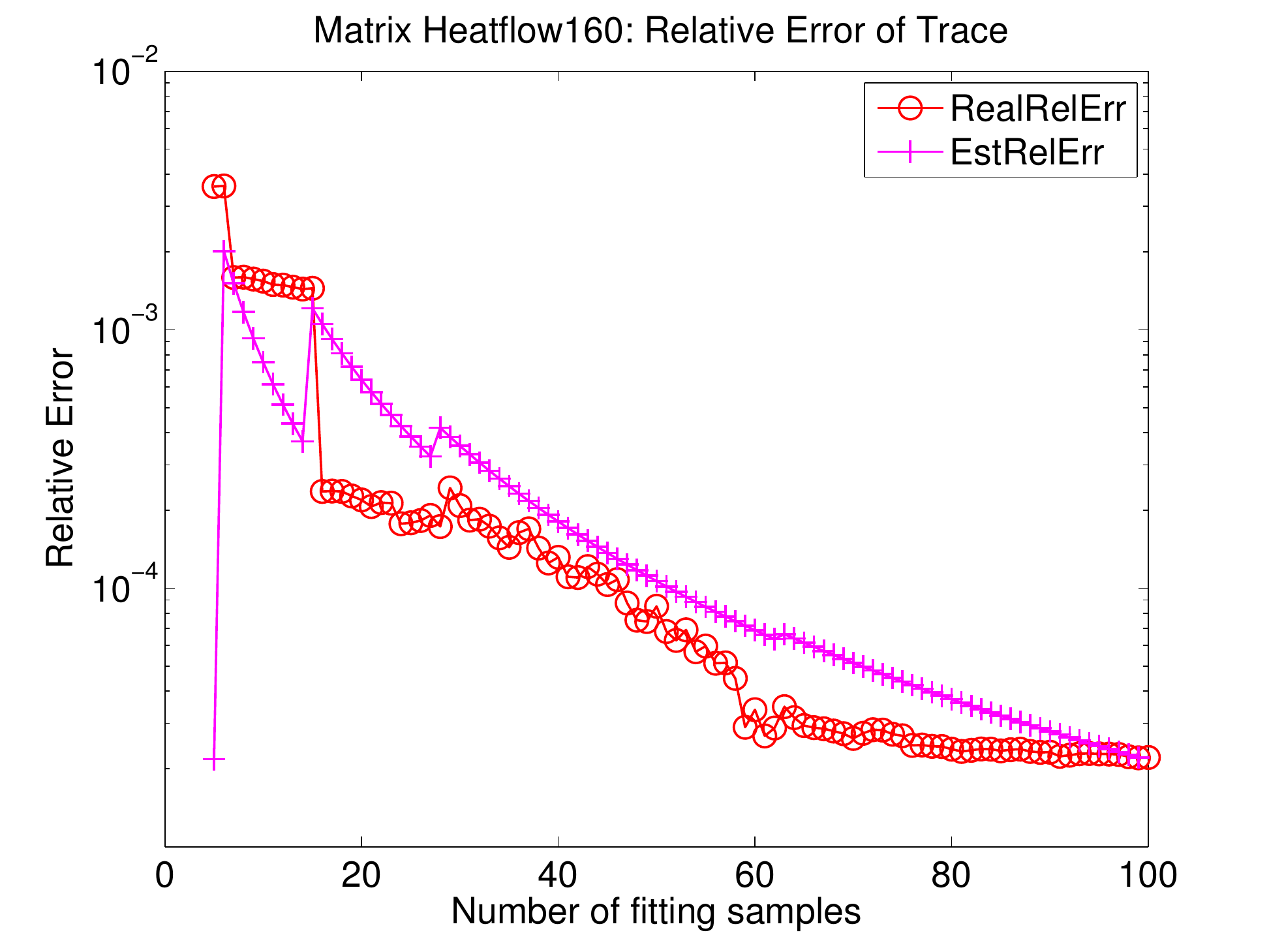}\label{fig:comparison of estimated and actual relative trace error on Matrix Heatflow160}}  
  \subfigure[Matrix Poisson150]
  {\includegraphics[width=0.32\textwidth]{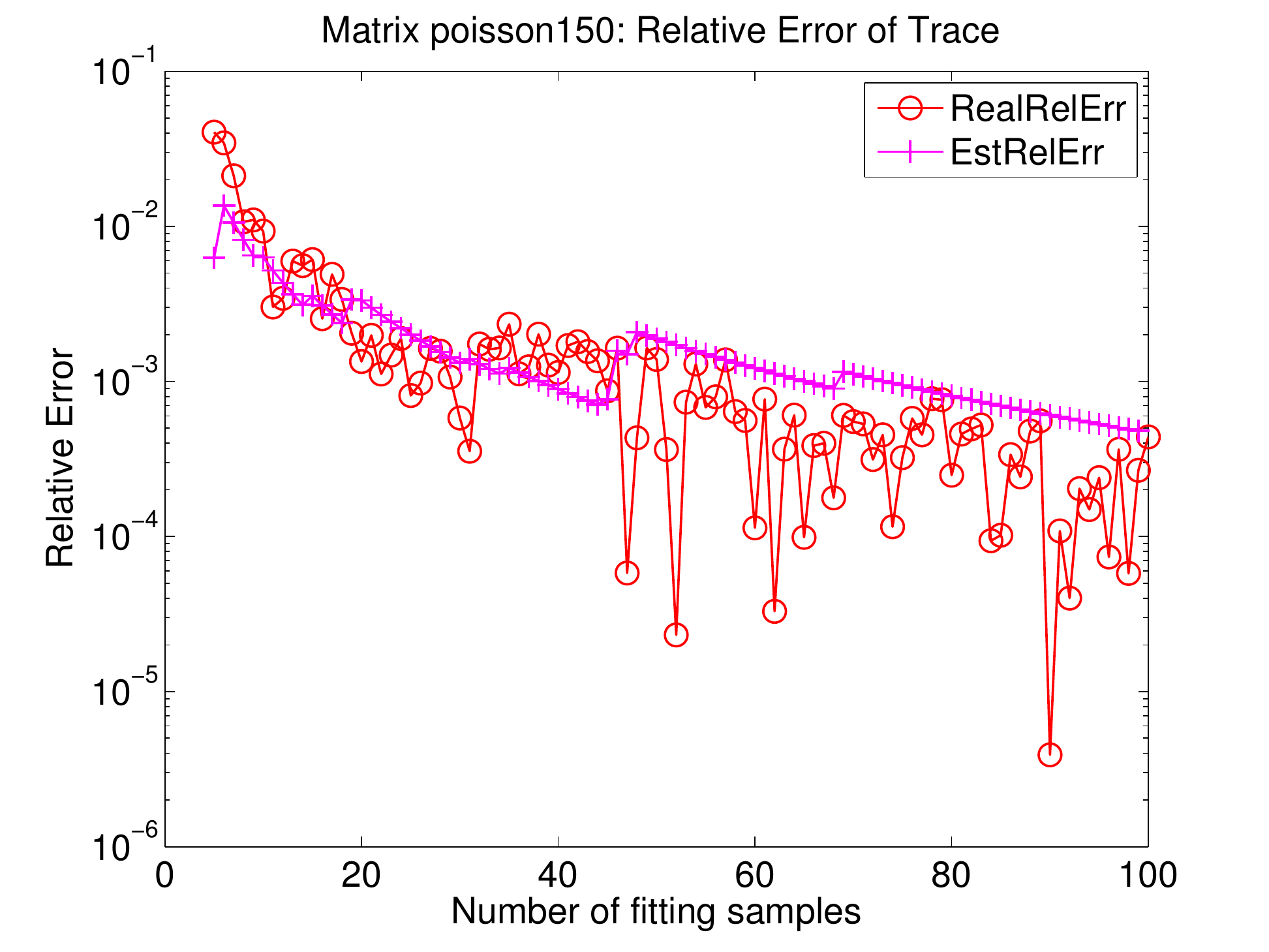}\label{fig:comparison of estimated and actual relative trace error on Matrix poisson150}}                
  \subfigure[Matrix VFH6]
  {\includegraphics[width=0.32\textwidth]{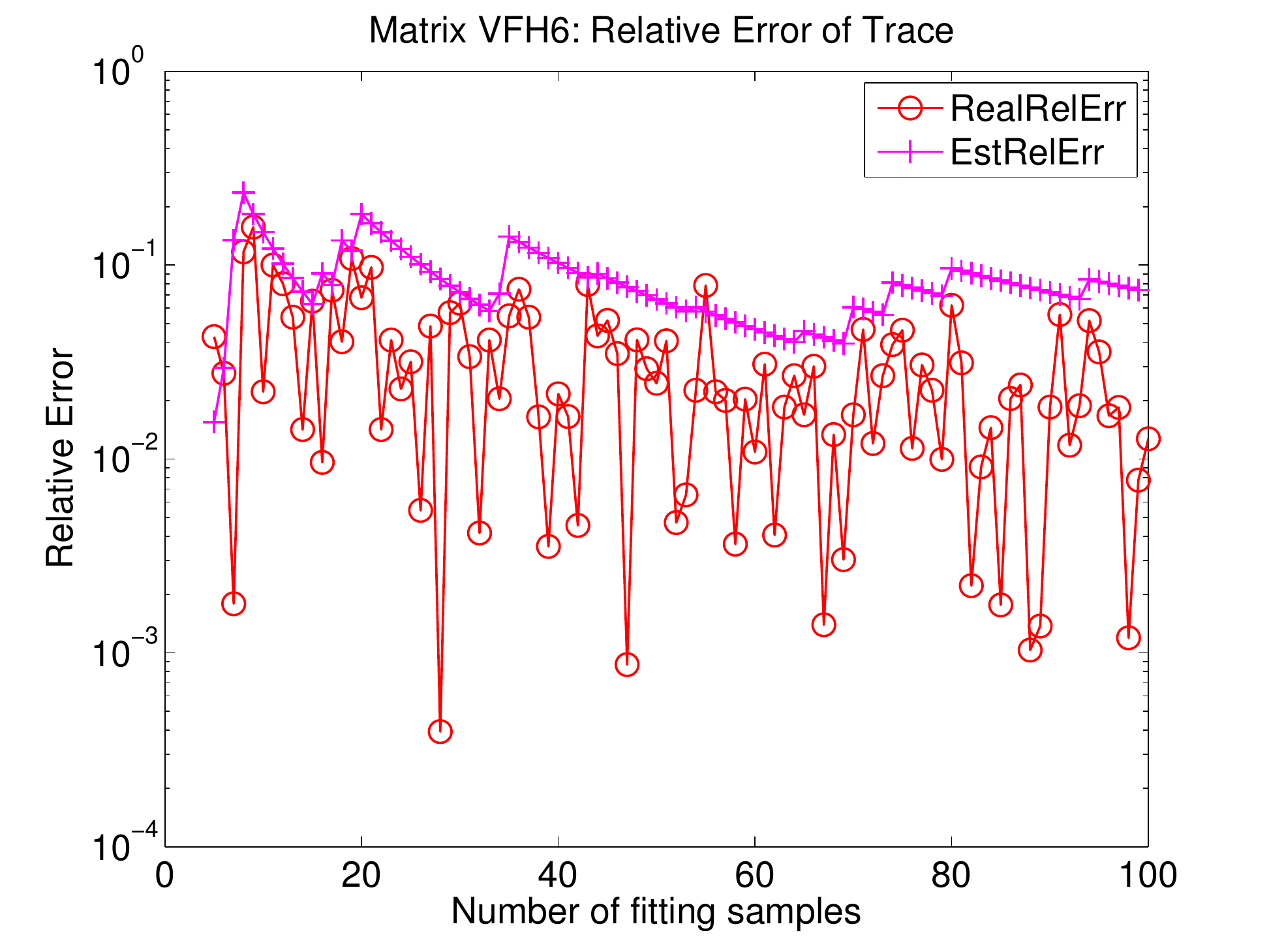}\label{fig:comparison of estimated and actual relative trace error on Matrix VFH6}}   
  \caption{Comparing estimated relative trace error with actual relative trace error with SVD.}
  \label{fig:comparison of the estimated relative trace error with actual relative trace error with SVD}
\end{figure}

Figure \ref{fig:comparison of the estimated relative trace error with actual relative trace error with SVD} compares the estimated relative trace error with the actual one in the cases of Figure \ref{fig:comparison of the estimated variances with actual variances of MC in various cases with SVD}.
We observe that the estimation is accurate as the fitting samples increase, even for cases such as VFH6 where the fitting process is not as successful.
Moreover, because our algorithm is based on upper bounds on the error of 
a piecewise cubic polynomial, the actual relative trace error could be lower than 
  predicted.

%---------------------------------------------------------------
\subsection{A large QCD problem}
The trace estimator presented in this paper has the potential of improving a number of LQCD calculations, where the trace of the Dirac matrix is related to an important property of QCD called spontaneous chiral symmetry breaking \citep{stathopoulos2010computing}. 
The authors in \citep{stathopoulos2013hierarchical}
presented the method of hierarchical probing that achieves almost optimal 
  variance reduction incrementally and inexpensively on regular lattices.
%We show that our fitting method can compute the trace much more accurately 
 %than hierarchical probing.

As shown in Figure \ref{fig:relative trace error, variance and PCHIP on matrix matb5 with SVD}(a), 
  a low rank approximation with 200 singular vectors yields a good approximation $M$ and 
  an excellent fit $p(M)$.
In Figure \ref{fig:relative trace error, variance and PCHIP on matrix matb5 with SVD}(b),
  we see that the actual relative trace error decreases very fast to the order of {\tt 1e-4}
  with increasing number of fitting points and singular vectors, 
  and can be monitored well by our dynamic trace estimation algorithm.
These singular vectors can be approximated while solving the linear systems with eigCG.
Interestingly, we can improve the relative trace error of hierarchical probing 
  by two orders of magnitude.
In addition, the variances of different MC methods can be estimated dynamically
  to allow us to continue with the best estimator if needed 
  (Figure \ref{fig:relative trace error, variance and PCHIP on matrix matb5 with SVD}(c)).

\begin{figure}
  \centering
  \subfigure[PCHIP results]
  {\includegraphics[width=0.32\textwidth]{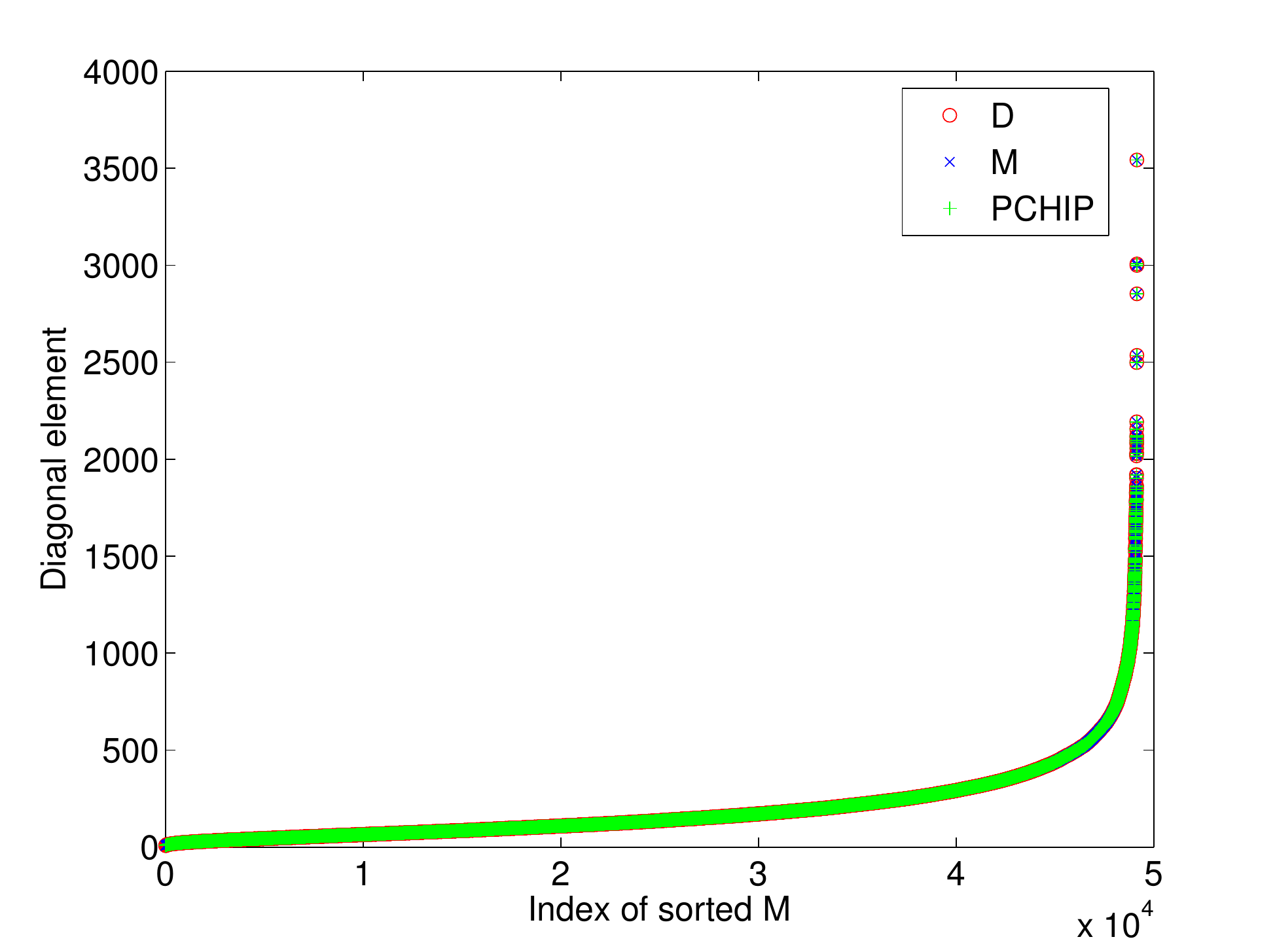}\label{fig:PCHIP on Matrix matb5}}               
  \subfigure[Relative trace error]
  {\includegraphics[width=0.32\textwidth]{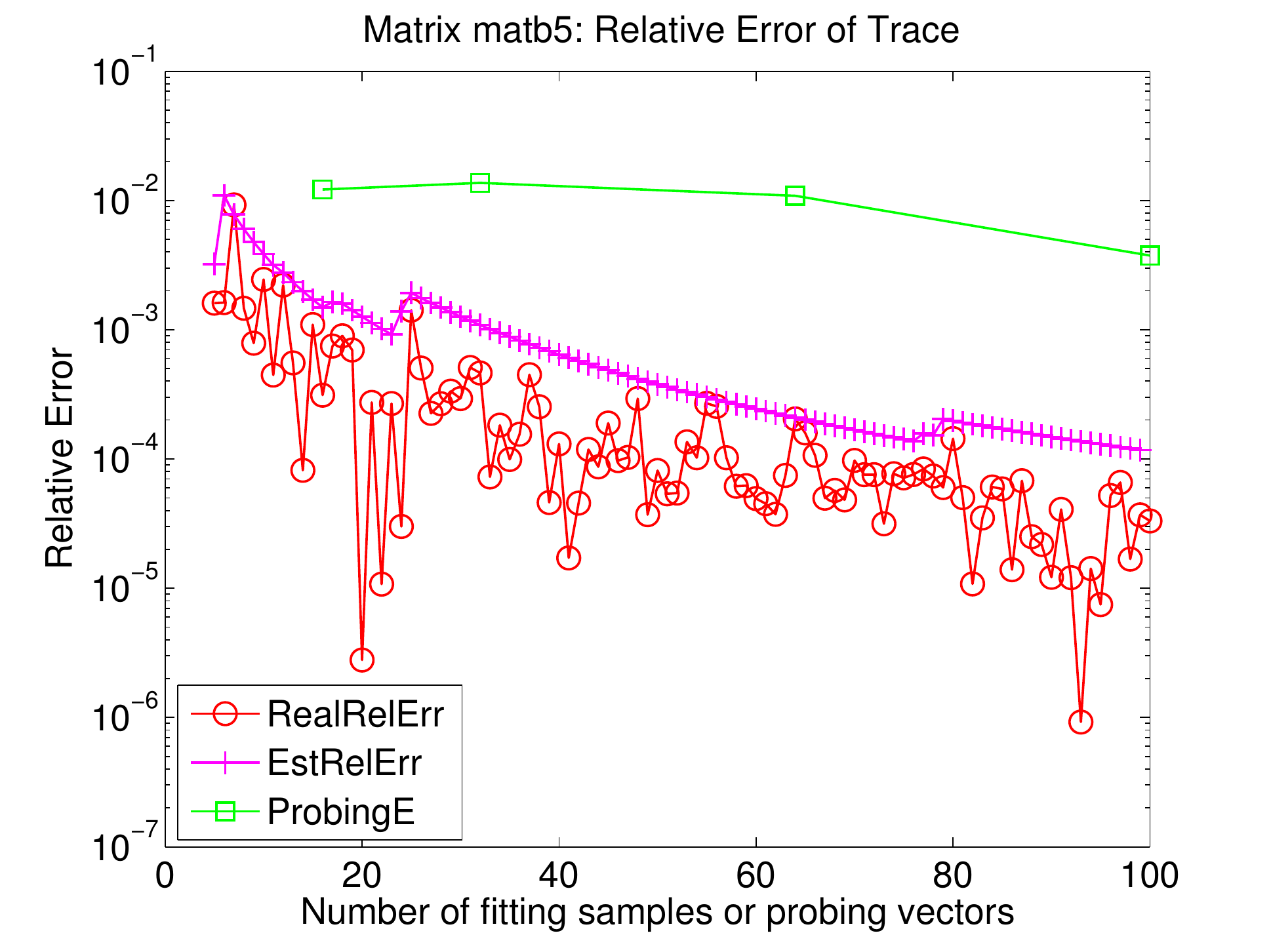}\label{fig:dynamic evaluation of relative trace error on Matrix matb5}}
  \subfigure[Variances]
  {\includegraphics[width=0.32\textwidth]{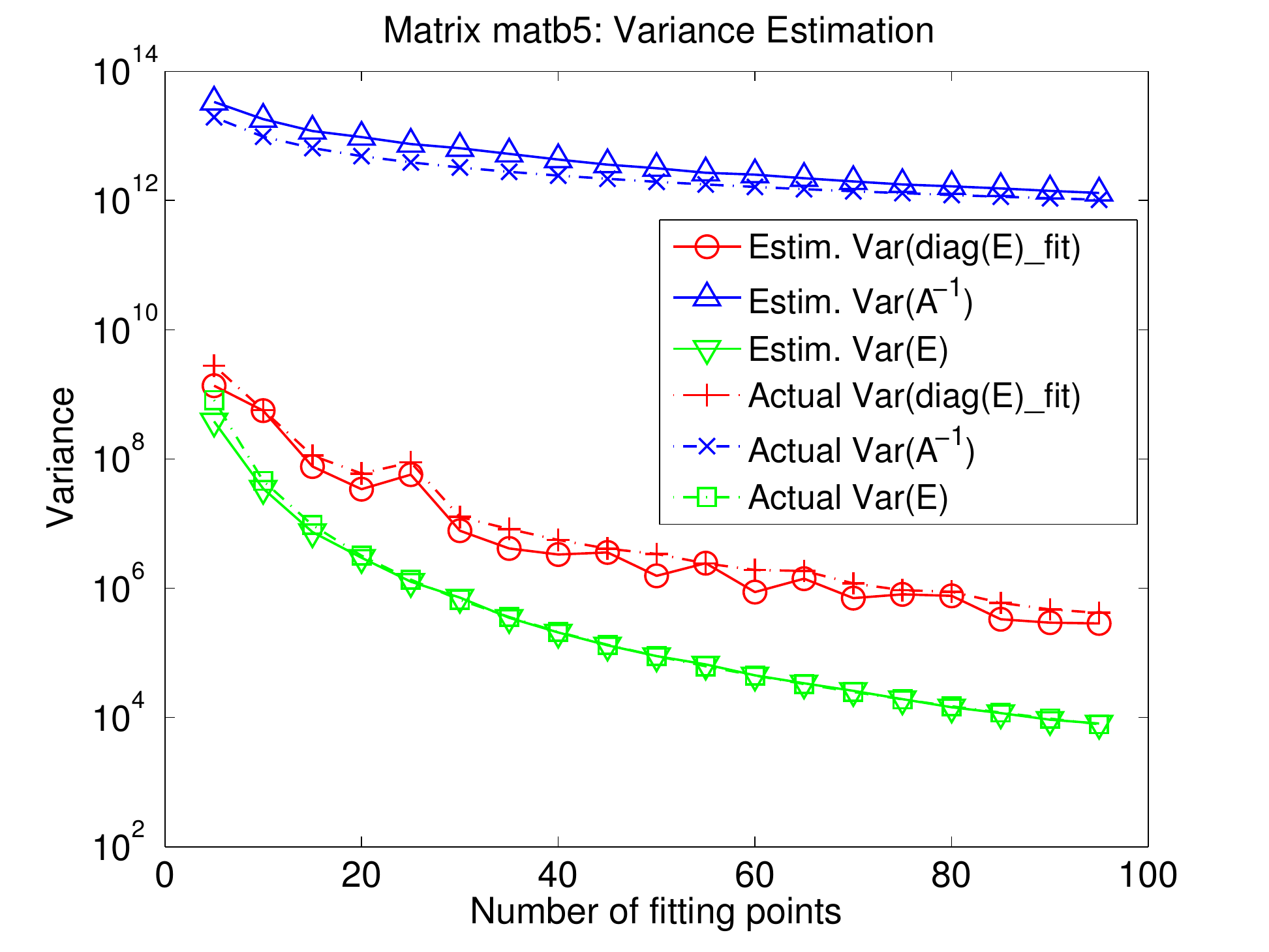}\label{fig:dynamic evaluation of variance on Matrix matb5}}    
  \caption{Fitting results, dynamic evaluation of relative trace error and variances with SVD on a large QCD matrix. In Figure \ref{fig:dynamic evaluation of relative trace error on Matrix matb5}, the green square denotes the relative trace error by applying hierarchical probing technique on deflated matrix $E$ in \citep{stathopoulos2013hierarchical}.}
  \label{fig:relative trace error, variance and PCHIP on matrix matb5 with SVD}
\end{figure}

%---------------------------------------------------------------
\section{Conclusion and future work}
A novel method has been presented to estimate the trace of the matrix inverse by exploiting the pattern correlation between the diagonal of the inverse of the matrix and some approximation. The key idea is to construct a good approximation $M\approx D$ through eigenvectors or some preconditioner, sample important patterns of $D$ by using the distribution of the elements of $M$, and use fitting techniques to obtain a better approximation 
$p(M)\approx D$ from where we obtain a trace estimate. The proposed method can provide a fast trace estimate with 2-3 digits relative accuracy given only a few samples while may or may not improve the variance of MC. When the variance is reduced sufficiently, our method can be also used as a diagonal estimator. We also propose an effective dynamic variance evaluation algorithm to determine the MC method with the smallest variance and a dynamic relative trace error estimation algorithm without any additional costs. 
We demonstrated the effectiveness of these methods through a set of experiments in some real applications.

\section*{Acknowledgement}
The authors thank Professor Yousef Saad for his helpful comments and discussions to improve the manuscript. The authors also thank reviewers' valuable comments to improve the quality and readability of the manuscript. This work is supported by NSF under grants No. CCF 1218349 and ACI SI2-SSE 1440700, and by DOE under a grant No. DE-FC02-12ER41890.

\section*{References}
  \bibliographystyle{elsarticle-num} 
  \bibliography{trace_JCP}

%% else use the following coding to input the bibitems directly in the TeX file.

%---------------------------------------------------------------
%\section*{References}
%\begin{center}
%\begin{thebibliography}{99}
%\parskip0pt 
%\itemsep0pt 
%
%\bibitem{avron2011randomized}
%Haim Avron and Sivan Toledo, \emph{Randomized algorithms for estimating the trace of an implicit symmetric positive semi-definite matrix}, Journal of the ACM, 58 (2001), p. Article 8.

\end{document}